\documentclass[11pt]{article}
\usepackage{latexsym}
\usepackage{amssymb}
\usepackage{amsmath}
\usepackage[francais]{babel}

\def\BBq{\mathbb Q}
\def\BBz{\mathbb Z}
\def\BBf{\mathbb F}

\def\BBp{\mathbb P}
\def\BBc{\mathbb C}

\def\BBn{\mathbb N}

\input xy
\xyoption{all}

\begingroup
\newtheorem{bigthm}{Th\'eor\`eme}
\newtheorem{thm}{Th\'eor\`eme}[section]
\newtheorem{prop}[thm]{Proposition}
\newtheorem{cor}[thm]{Corollaire}
\newtheorem{lem}[thm]{Lemme}
\newtheorem{rem}[thm]{Remarque}
\endgroup

\newenvironment{pf}{\noindent {\em Preuve.}}{\hfill $\Box$}

\begin{document}

\title{Les repr\'esentations $\ell$-adiques \\
associ\'ees aux courbes elliptiques sur ${\BBq}_p$}

\author{Maja Volkov}

\maketitle

\begin{abstract}
This paper is devoted to the study of the $\ell$-adic representations of the absolute Galois group $G$ of ${\BBq}_p$, $p\geq 5$, associated to an elliptic curve over ${\BBq}_p$, as $\ell$ runs through the set of all prime numbers (including $\ell =p$, in which case we use the theory of potentially semi-stable $p$-adic representations). \\
For each prime $\ell$, we give the complete list of isomorphism classes of ${\BBq}_{\ell}[G]$-modules coming from an elliptic curve over ${\BBq}_p$, that is, those which are isomorphic to the Tate module of an elliptic curve over ${\BBq}_p$. The $\ell =p$ case is the more delicate. It requires studying the liftings of a given elliptic curve over ${\BBf}_p$ to an elliptic scheme over the ring of integers of a totally ramified finite extension of ${\BBq}_p$, and combining it with a descent theorem providing a Galois criterion for an elliptic curve having good reduction over a $p$-adic field to be defined over a closed subfield. This enables us to state necessary and sufficient conditions for an $\ell$-adic representation of $G$ to come from an elliptic curve over ${\BBq}_p$, for each prime $\ell$.
\end{abstract}


\medskip

\bigskip
L'objet de cet article est l'\'etude des repr\'esentations $\ell$-adiques associ\'ees aux courbes elliptiques d\'efinies sur ${\BBq}_p$, o\`u $p$ est un nombre premier sup\'erieur ou \'egal \`a 5, lorsque $\ell$ parcourt l'ensemble de tous les nombres premiers, y compris $\ell =p$.

\medskip
Fixons un premier $p\geq 5$ et une cl\^oture alg\'ebrique $\overline{{\BBq}}_p$ de ${\BBq}_p$. Soit $E$ une courbe elliptique sur ${\BBq}_p$. Pour tout $\ell$ premier, soit $E[\ell^n]$, $n\geq 0$, le groupe des points de $\ell^n$-torsion de $E$ \`a valeurs dans $\overline{{\BBq}}_p$. Le module de Tate $\ell$-adique $\displaystyle T_{\ell}(E)=\varprojlim E[\ell^n]$ est un ${\BBz}_{\ell}$-module libre de rang 2, $V_{\ell}(E)={\BBq}_{\ell}\otimes_{{\BBz}_{\ell}}T_{\ell}(E)$ est un ${\BBq}_{\ell}$-espace vectoriel de dimension 2, et tous deux sont munis d'une action lin\'eaire et continue du groupe de Galois absolu $G=\mbox{Gal}(\overline{{\BBq}}_p/{\BBq}_p)$. On obtient ainsi des repr\'esentations pour tout $\ell$ 
 $$G\longrightarrow \mbox{Aut}_{{\BBz}_{\ell}}(T_{\ell}(E)) \makebox[1cm]{et} G\longrightarrow \mbox{Aut}_{{\BBq}_{\ell}}(V_{\ell}(E))$$ 
Ces repr\'esentations contiennent des informations concernant la courbe $E/{\BBq}_p$ et beaucoup de r\'esultats sur celles-ci sont devenus classiques (voir \cite{Se 1}, \cite{Se 2}, \cite{Se-Ta}, \cite{Ta}, \cite{Kr}, \cite{Ro}, et bien d'autres).

Soit maintenant $T_{\ell}$ un ${\BBz}_{\ell}$-module libre de rang 2 muni d'une action lin\'eaire et continue de $G$. On consid\`ere le probl\`eme suivant : quand $T_{\ell}$ provient-il d'une courbe elliptique sur ${\BBq}_p$, i.e. quand existe-t-il une courbe elliptique $E/{\BBq}_p$ telle que $T_{\ell}$ et $T_{\ell}(E)$ sont des ${\BBz}_{\ell}[G]$-modules isomorphes~? En fait, cette question se ram\`ene \`a celle obtenue en rempla\c cant $T_{\ell}$ par ${\BBq}_{\ell}\otimes_{{\BBz}_{\ell}}T_{\ell}$ et $T_{\ell}(E)$ par $V_{\ell}(E)$. On donne ici une r\'eponse \`a cette question, y compris dans le cas $\ell =p$, qui est le plus d\'elicat \`a traiter.

\medskip
L'objet de la section~\ref{sec:rappelsnots} est de pr\'esenter les outils utilis\'es. Pour $\ell \neq p$, \`a chaque repr\'esentation $\ell$-adique de $G$, on sait associer fonctoriellement une repr\'esentation de Weil-Deligne - sur laquelle les racines du polyn\^ome caract\'eristique d'un rel\`evement du Frobenius sont des unit\'es $\ell$-adiques - et vice versa (voir une d\'efinition ici en \ref{sec:repwd} ; il est important de noter que la construction que nous allons utiliser est contravariante). Une repr\'esentation de Weil-Deligne consiste en la donn\'ee d'une repr\'esentation du groupe de Weil et d'un op\'erateur nilpotent, les deux \'etant li\'es par une relation. Travailler avec de tels objets pr\'esente l'avantage de discr\'etiser l'action de $G$. Pour $\ell =p$, \`a chaque repr\'esentation $p$-adique de $G$ potentiellement semi-stable, on sait associer fonctoriellement un $(\varphi,N,G)$-module filtr\'e (voir \ref{sec:reppst} ; ici encore nous allons utiliser un foncteur contravariant), et r\'eciproquement si ce module filtr\'e est faiblement admissible (i.e. faiblement admissible \'equivaut \`a admissible, voir \cite{Co-Fo}). Un $(\varphi,N,G)$-module filtr\'e consiste en la donn\'ee d'un espace vectoriel muni d'un op\'erateur de Frobenius $\varphi$ semi-lin\'eaire, d'un op\'erateur nilpotent $N$ li\'e \`a $\varphi$ par une relation, d'une action semi-lin\'eaire de $G$ commutant avec ces deux op\'erateurs, ainsi que d'une filtration stable par l'action de $G$ ; la faible admissibilit\'e est une condition liant le Frobenius \`a la filtration. Travailler avec de tels objets pr\'esente l'avantage de pouvoir remplacer une repr\'esentation $p$-adique de $G$ par des donn\'ees de type alg\'ebrique. De plus, en oubliant la filtration, on sait associer fonctoriellement \`a chaque $(\varphi,N,G)$-module filtr\'e une repr\'esentation de Weil-Deligne d\'efinie sur une extension non ramifi\'ee de ${\BBq}_p$ ; via cette association, les classes d'isomorphisme des $(\varphi,N,G)$-modules (non filtr\'es) correspondent bijectivement \`a celles des repr\'esentations de Weil-Deligne. Cela permet de comparer les repr\'esentations $\ell$-adiques entre elles pour tout $\ell$ premier. En particulier, on sait que, pour une courbe elliptique $E/{\BBq}_p$ fix\'ee, les repr\'esentations de Weil-Deligne qui lui sont associ\'ees sont ind\'ependantes du nombre premier $\ell$. Ainsi, r\'epondre \`a la question pour $\ell \neq p$, c'est savoir quelles sont exactement les repr\'esentations de Weil-Deligne possibles~; et pour $\ell =p$, c'est savoir quelles sont exactement les filtrations faiblement admissibles possibles. 

\medskip
Donnons des conditions n\'ecessaires bien connues. Soit $E/{\BBq}_p$ une courbe elliptique.  \\
1) Pour tout $\ell$ premier, la repr\'esentation de Weil-Deligne associ\'ee \`a $V_{\ell}(E)$ v\'erifie~:
\begin{itemize}
\item[$(1^{\circ})$] le d\'eterminant sur $V_{\ell}(E)$ est le caract\`ere cyclotomique $\ell$-adique 
\item[$(2^{\circ})$] elle est d\'efinie sur ${\BBq}$ 
\item[$(3^{\circ})$] si $E$ a potentiellement bonne r\'eduction, les racines du polyn\^ome caract\'eristique d'un rel\`evement du Frobenius g\'eom\'etrique sont des $p$-nombres de Weil : $\mbox{Tr(Frob)} \in {\BBz}$ et $\mid\!\mbox{Tr(Frob)}\!\mid_{\infty} \leq 2\sqrt{p}\,$
\end{itemize}
2) Le ${\BBq}_p[G]$-module $V_p(E)$ est potentiellement semi-stable et de type Hodge-Tate $(0,1)$.

\medskip
En fait, on a un but double : d'une part, classifier pour tout $\ell$ les repr\'esentations $\ell$-adiques associ\'ees \`a une courbe elliptique sur ${\BBq}_p$ ; d'autre part, faire la liste pour tout $\ell$ de toutes les repr\'esentations $\ell$-adiques, \`a isomorphisme pr\`es, v\'erifiant ces conditions n\'ecessaires et d\'eterminer celles qui proviennent d'une courbe elliptique sur ${\BBq}_p$.

\medskip
La premi\`ere partie est l'objet de la section~\ref{sec:classif}; signalons que les r\'esultats \'enonc\'es dans cette section sont tous plus ou moins connus. On commence par construire une liste finie de repr\'esentations de Weil-Deligne deux \`a deux non isomorphes, liste que nous notons ${\bf WD^*}$ (\ref{sec:listl}). \`A chacune d'elles, on associe un $(\varphi,N,G)$-module et on construit un repr\'esentant de chaque classe d'isomorphisme de filtration faiblement admissible de type Hodge-Tate $(0,1)$ que l'on peut mettre sur ce module ; il y a, suivant les cas, un nombre fini ou infini de possibilit\'es. Cela nous donne une liste infinie de $(\varphi,N,G)$-modules filtr\'es faiblement admissibles deux \`a deux non isomorphes que nous notons ${\bf D^*}$ (\ref{sec:listp}). Puis on montre que si $E$ est une courbe elliptique sur ${\BBq}_p$, alors, pour $\ell \neq p$, la repr\'esentation de Weil-Deligne associ\'ee \`a $V_{\ell}(E)$ est isomorphe \`a un objet de la liste ${\bf WD^*}$ (\ref{sec:classl}), et le $(\varphi,N,G)$-module filtr\'e associ\'e \`a $V_p(E)$ est isomorphe \`a un objet de la liste ${\bf D^*}$ (\ref{sec:classp}). De plus, si l'on se donne la courbe $E$ sous la forme d'une \'equation de Weierstrass $y^2=x^3+Ax+B$, on donne une liste d'invariants de $E$ d\'efinis \`a partir de $A$ et de $B$ qui permettent de d\'eterminer l'objet de ${\bf WD^*}$ associ\'e et parfois aussi de ${\bf D^*}$.

\medskip
La deuxi\`eme partie est l'objet des sections~\ref{sec:lmodellip}, \ref{sec:construct} et \ref{sec:pmodellip}. Les r\'esultats principaux sont les deux th\'eor\`emes suivants :  

\begin{bigthm}
{\em (cf. thm. \ref{lmodellipthm} ci-dessous)} Soient $\ell \neq p$ et $V_{\ell}$ une repr\'esentation $\ell$-adique de $G$ de dimension $2$. Les assertions suivantes sont \'equivalentes : 
\begin{itemize}
\item[(1)] il existe une courbe elliptique $E$ sur ${\BBq}_p$ telle que $V_{\ell}(E)$ soit isomorphe \`a $V_{\ell}$,
\item[(2)] la repr\'esentation de Weil-Deligne associ\'ee \`a $V_{\ell}$ v\'erifie les conditions $(1^{\circ})$, $(2^{\circ})$ et $(3^{\circ})$ ci-dessus, 
\item[(3)] la repr\'esentation de Weil-Deligne associ\'ee \`a $V_{\ell}$ est isomorphe \`a un objet de la liste ${\bf WD^*}$.
\end{itemize}
\end{bigthm}

\begin{bigthm}
{\em (cf. thm. \ref{pmodellipthm} ci-dessous)} Soit $V_p$ une repr\'esentation $p$-adique de $G$ de dimension $2$. Les assertions suivantes sont \'equivalentes : 
\begin{itemize}
\item[(1)] il existe une courbe elliptique $E$ sur ${\BBq}_p$ telle que $V_p(E)$ soit isomorphe \`a $V_p$,
\item[(2)] la repr\'esentation $V_p$ est potentiellement semi-stable de type Hodge-Tate $(0,1)$ et la repr\'esentation de Weil-Deligne associ\'ee v\'erifie les conditions $(1^{\circ})$, $(2^{\circ})$ et $(3^{\circ})$ ci-dessus,
\item[(3)] la repr\'esentation $V_p$ est potentiellement semi-stable et le $(\varphi,N,G)$-module filtr\'e associ\'e est isomorphe \`a un objet de la liste ${\bf D^*}$.
\end{itemize}
\end{bigthm}

Le th\'eor\`eme 1 est facile \`a d\'emontrer (\ref{sec:resultl}). En effet, en utilisant la th\'eorie de Honda-Tate pour les cas de potentielle bonne r\'eduction, les op\'erations \'el\'ementaires sur les courbes elliptiques permettent de construire suffisament d'exemples pour obtenir toutes les classes ; on peut m\^eme produire pour chacune un exemple sous forme d'\'equation de Weierstrass (\ref{sec:exl}).

Le th\'eor\`eme 2 est plus d\'elicat et constitue le principal r\'esultat nouveau de cet article. Lorsque la repr\'esentation est potentiellement semi-stable mais non potentiellement cristalline on peut faire des calculs suffisament explicites, gr\^ace aux courbes de Tate ; on obtient une infinit\'e de classes param\'etr\'ees par ${\BBz}/2{\BBz} \times {\BBz}/2{\BBz} \times {\BBq}_p$. Les cas cristallins ou tordus-cristallins engendrent une famille finie de classes, de m\^eme pour les cas potentiellement ordinaires ; on peut donner des exemples sous forme d'\'equations de Weierstrass (\ref{sec:exp}). La principale difficult\'e se trouve dans les cas potentiellement supersinguliers qui ne sont pas des tordus de cas cristallins, ce qui arrive lorsque $12$ ne divise pas $p-1$ ; alors pour chaque entier dans $\{ 3,4,6\}$ divisant $p+1$ on obtient une infinit\'e de classes param\'etr\'ees par ${\BBp}^1({\BBq}_p)$.

Pour traiter ce dernier cas on construit dans la section~\ref{sec:construct} toutes les courbes elliptiques sur ${\BBq}_p$ ayant potentiellement bonne r\'eduction supersinguli\`ere, \`a ${\BBq}_p$-isomorphisme pr\`es, de la mani\`ere qui suit. Dans un premier temps, \'etant donn\'ee une courbe elliptique $\widetilde{E}$ sur ${\BBf}_p$ supersinguli\`ere, on d\'ecrit en \ref{sec:schellsup} les sch\'emas elliptiques la relevant sur l'anneau des entiers d'une extension totalement ramifi\'ee dont l'indice de ramification est un entier $e$ strictement inf\'erieur \`a $p-1$ (la r\'eponse \`a ce probl\`eme est bien connue quand $\widetilde{E}$ est ordinaire, cf. \cite{Me} ou \cite{Ka}). Pour cela on combine le th\'eor\`eme de Serre-Tate (\ref{sec:ST}) avec la description des groupes $p$-divisibles par les modules de Dieudonn\'e filtr\'es (\ref{sec:MD} et \ref{sec:MDfil}). Puis on d\'emontre en \ref{sec:galdesc} un th\'eor\`eme de descente qui fournit un crit\`ere galoisien pour qu'une courbe elliptique ayant bonne r\'eduction sur un corps $p$-adique puisse \^etre d\'efinie sur un sous-corps ferm\'e. Ce crit\`ere permet de d\'eterminer quels sont, parmi les sch\'emas elliptiques construits en \ref{sec:schellsup} et pour $e \in \{ 3,4,6\}$, ceux qui sont susceptibles d'\^etre d\'efinis sur ${\BBq}_p$ (\ref{sec:cpotsup}). En \ref{sec:cpotord} on donne les r\'esultats de ces m\'ethodes appliqu\'ees aux cas ordinaires. On r\'ecolte finalement les fruits de cette \'etude en \ref{sec:resultp} o\`u l'on d\'emontre le th\'eor\`eme 2 \'enonc\'e ci-dessus.

\medskip
Cet article est une version remani\'ee de la th\`ese de l'auteur sous la direction de J.-M. Fontaine. Au cours de ce travail l'auteur a b\'en\'efici\'e de pr\'ecieuses discussions avec lui. 

Je tiens donc \`a remercier chaleureusement J.-M. Fontaine, sans ses patientes explications ce travail n'aurait certainement pas pu \^etre men\'e \`a bien.

\bigskip

\tableofcontents

\bigskip

\section{Rappels et notations}
\label{sec:rappelsnots}

Soit ${\cal P}$ l'ensemble des nombres premiers. On fixe un $p \in {\cal P}$ tel que $p \geq 5$ et $\overline{{\BBq}}_p$ une cl\^oture alg\'ebrique de ${\BBq}_p$. On note $G=\mbox{Gal}(\overline{{\BBq}}_p/{\BBq}_p)$ et $I$ son sous-groupe d'inertie ; si $K\subset \overline{{\BBq}}_p$ est une extension de ${\BBq}_p$ on pose $G_K=\mbox{Gal}(\overline{{\BBq}}_p/K)$ et $I_K=G_K \cap I$. On note ${\BBq}_p^{nr}$ l'extension maximale non ramifi\'ee de ${\BBq}_p$ contenue dans $\overline{{\BBq}}_p$ et $\overline{{\BBf}}_p$ son corps r\'esiduel. Pour $\ell \in {\cal P}$, si $\mu_{\ell^n}(\overline{{\BBq}}_p)$ est le groupe des racines $\ell^n$-i\`emes de l'unit\'e contenues dans $\overline{{\BBq}}_p$, on \'ecrit $\displaystyle {\BBz}_{\ell}(1)= \varprojlim \mu_{\ell^n}(\overline{{\BBq}}_p)$ et ${\BBq}_{\ell}(1)={\BBq}_{\ell} \otimes_{{\BBz}_{\ell}} {\BBz}_{\ell}(1)$. La valuation $p$-adique $v_p$ sur ${\BBq}_p$ est normalis\'ee par $v_p(p)=1$ ; on note aussi $v_p$ la valuation qui l'\'etend sur $\overline{{\BBq}}_p$.

\subsection{Repr\'esentations $\ell$-adiques, $\ell \neq p$, et repr\'esentations de Weil-Deligne}
\label{sec:repwd}

Soit $\ell \in {\cal P}$ tel que $\ell \neq p$. On d\'esigne par ${\bf Rep}_{{\BBq}_{\ell}}(G)$ la cat\'egorie des repr\'esentations $\ell$-adiques de $G$, c'est-\`a-dire des ${\BBq}_{\ell}$-espaces vectoriels de dimension finie munis d'une action lin\'eaire et continue de $G$. Soit $L \subset \overline{{\BBq}}_p$ une extension finie de ${\BBq}_p$ ; une repr\'esentation $\ell$-adique de $G$ est semi-stable sur $L$ si $I_L$ op\`ere de fa\c con unipotente et a bonne r\'eduction sur $L$ si $I_L$ op\`ere trivialement. Le corps r\'esiduel de ${\BBq}_{\ell}$ \'etant fini, on sait que toutes les repr\'esentations $\ell$-adiques de $G$ sont potentiellement semi-stables.

\medskip
Le groupe de Weil $W$ de ${\BBq}_p$ est d\'efini par la suite exacte courte $ 1\rightarrow I\rightarrow W \stackrel{\upsilon}{\rightarrow} {\BBz}\rightarrow 1$ avec $\upsilon (\mbox{Frob.arithm.})=1$ ; c'est donc le sous-groupe de $G$ constitu\'e des \'el\'ements $g$ tels que $g \bmod I$ est une puissance enti\`ere du Frobenius.

Soit $K$ un corps de caract\'eristique $0$. On note ${\bf Rep}_{K}('W)$ la cat\'egorie des repr\'esentations $K$-lin\'eaires du groupe de Weil-Deligne $'W$ : les objets sont les triplets $({\Delta},{\rho}_0,N)$, o\`u ${\Delta}$ est un $K$-espace vectoriel de dimension finie, ${\rho}_0: W \rightarrow \mbox{Aut}_{K}({\Delta})$ est un morphisme dont le noyau contient un sous-groupe ouvert de $I$, et $N \in \mbox{End}_{K}({\Delta})$ v\'erifie ${\rho}_0(w)N=p^{\upsilon (w)}N{\rho}_0(w)$ pour tout $w\in W$ ; l'op\'erateur de monodromie $N$ est donc nilpotent. Soit $L\subset \overline{{\BBq}}_p$ une extension finie de ${\BBq}_p$ ; la repr\'esentation $({\Delta},{\rho}_0,N)$ est semi-stable sur $L$ si ${\rho}_0(I_L)=1$ et a bonne r\'eduction sur $L$ si $N=0$ et ${\rho}_0(I_L)=1$. Si l'action de $'W$ est $F$-semi-simple, i.e. si l'action de $W$ par ${\rho}_0$ est semi-simple, alors un objet de ${\bf Rep}_{K}('W)$ est d\'etermin\'e \`a isomorphisme pr\`es par les traces $\mbox{Tr}({\rho}_0) : W \rightarrow K$ ainsi que par le polyn\^ome minimal de $N$.

Soient $K$ et $K'$ deux corps de caract\'eristique $0$ munis de plongements $\iota : K \hookrightarrow {\BBc}$ et $\iota' : K' \hookrightarrow {\BBc}\,$ et soient $\Delta$, $\Delta'$ des objets de ${\bf Rep}_{K}('W)$ et ${\bf Rep}_{K'}('W)$ respectivement. On dit que $\Delta$ est d\'efini sur ${\BBq}$ si, \'etant donn\'es un ${\BBq}$-espace vectoriel $D$ tel que $\Delta = K\otimes_{\BBq}D$ et un corps alg\'ebriquement clos $\Omega$ contenant $K$, l'objet $\Omega\otimes_{\BBq}D$ de ${\bf Rep}_{\Omega}('W)$ est isomorphe \`a ses conjugu\'es sous $\mbox{Aut}(\Omega /{\BBq})$ (cette condition est ind\'ependante des choix de $D$ et de $\Omega$) ; dans ce cas $\Delta \otimes_{{K}^{\nearrow_{\!\iota}}}\!{\BBc}$ est un objet de ${\bf Rep}_{{\BBc}}('W)$ dont la classe d'isomorphisme ne d\'epend pas du choix de $\iota$. On dit que $\Delta$ et $\Delta'$ sont compatibles s'ils sont tous deux d\'efinis sur ${\BBq}$ et si $\Delta \otimes_{{K}^{\nearrow_{\!\iota}}}\!{\BBc}$ et $\Delta' \otimes_{{K'}^{\nearrow_{\!\iota'}}}\!{\BBc}$ sont isomorphes dans ${\bf Rep}_{{\BBc}}('W)$.

\medskip
Soit ${\bf Rep}^{\circ}_{{\BBq}_{\ell}}('W)$ la sous-cat\'egorie pleine de ${\bf Rep}_{{\BBq}_{\ell}}('W)$ form\'ee des objets sur lesquels les racines du polyn\^ome caract\'eristique d'un rel\`evement du Frobenius sont des unit\'es $\ell$-adiques. Il existe un foncteur \'etablissant une \'equivalence entre ${\bf Rep}^{\circ}_{{\BBq}_{\ell}}('W)$ et ${\bf Rep}_{{\BBq}_{\ell}}(G)$, voir \cite{Ro}, $\S$ 4 ou \cite{De 1}, $\S$ 8.

On renvoie \`a \cite{Fo 3} pour la d\'efinition de l'anneau $B_{st,\ell}$. On utilise ici le foncteur contravariant ${\bf W}^*_{\ell} : {\bf Rep}_{{\BBq}_{\ell}}(G) \rightarrow {\bf Rep}^{\circ}_{{\BBq}_{\ell}}('W)$ donn\'e par ${\bf W}^*_{\ell}(V) = \mbox{Hom}_{{\BBq}_{\ell}[I_L]}(V,B_{st,\ell})$ si $V$ est semi-stable sur l'extension finie $L\subset \overline{{\BBq}}_p$ de ${\BBq}_p$ ; si $V$ a bonne r\'eduction sur $L$ on a ${\bf W}^*_{\ell}(V) = \mbox{Hom}_{{\BBq}_{\ell}[I_L]}(V,{\BBq}_{\ell})$. Ce foncteur est \'equivalent \`a celui obtenu en appliquant le foncteur d\'ecrit dans \cite{Fo 3} \`a la repr\'esentation duale. Il \'etablit une anti-\'equivalence de cat\'egories via laquelle les notions d'objet semi-stable sur $L$ ou ayant bonne r\'eduction sur $L$ se correspondent. Comme $V_{\ell}(E)$ est le dual de $H^1_{\mbox{\scriptsize\em \'et}}(E\!\times_{{\BBq}_p}\!\overline{{\BBq}}_p\, ,{\BBq}_{\ell})$ pour une courbe elliptique $E/{\BBq}_p$, on peut voir ${\bf W}^*_{\ell}$ comme un foncteur covariant des $H^1_{\mbox{\scriptsize\em \'et}}$ dans ${\bf Rep}^{\circ}_{{\BBq}_{\ell}}('W)$.

\medskip
Choisissons pour chaque $\ell \neq p$ un plongement de corps $\iota_{\ell} : {\BBq}_{\ell} \hookrightarrow {\BBc}\,$. Un syst\`eme $(\Delta_{\ell})_{\ell \neq p}$ de repr\'esentations ${\BBq}_{\ell}$-lin\'eaires de $'W$ est compatible si les $\Delta_{\ell}$ sont deux \`a deux compatibles lorsque $\ell$ parcourt ${\cal P} \backslash \{ p\}$. Sur chaque $\Delta_{\ell}$ il existe une unique filtration finie croissante $\{ \mbox{Fil}_i \Delta_{\ell} \}_{i\in {\BBz}}$ telle que $N(\mbox{Fil}_i \Delta_{\ell}) \subset \mbox{Fil}_{i-2} \Delta_{\ell}$ et que $N$ induit un isomorphisme $N^i : \mbox{Gr}_i \Delta_{\ell} \stackrel{\sim}{\rightarrow} \mbox{Gr}_{-i} \Delta_{\ell}$ pour tout $i\in {\BBz}$ (\cite{De 2}). Supposons que chaque $\Delta_{\ell}$ est $F$-semi-simple ; alors la compatibilit\'e signifie que, pour tout $i\in {\BBz}$, les traces $\mbox{Tr}(\mbox{Gr}_i \Delta_{\ell}) : W \rightarrow {\BBq}_{\ell}$ sont \`a valeurs dans ${\BBq}$ et ind\'ependantes de $\ell$. Si $N=0$ sur tous les $\Delta_{\ell}$ cela signifie que les traces $\mbox{Tr}(\rho_{0,\ell}) : W \rightarrow {\BBq}_{\ell}$ sont \`a valeurs dans ${\BBq}$ et ind\'ependantes de $\ell$.

\medskip
Soit $E/{\BBq}_p$ une courbe elliptique. Pour tout $\ell \neq p$, si $L\subset \overline{{\BBq}}_p$ est une extension finie de ${\BBq}_p$, alors $E$ est semi-stable sur $L$ (resp. a bonne r\'eduction sur $L$) si et seulement si la repr\'esentation de Weil-Deligne ${\bf W}^*_{\ell}(V_{\ell}(E))$ associ\'ee \`a $V_{\ell}(E)$ l'est. De plus, on sait que ${\bf W}^*_{\ell}(V_{\ell}(E))$ est $F$-semi-simple, d\'efinie sur ${\BBq}\,$, et que le syst\`eme $({\bf W}^*_{\ell}(V_{\ell}(E)))_{\ell \neq p}$ est compatible (voir \cite{Ra} pour un \'enonc\'e dans un contexte bien plus g\'en\'eral ; voir aussi la rmq. \ref{lcomp}).

\subsection{Repr\'esentations $p$-adiques potentiellement semi-stables et $(\varphi,N,G)$-modules filtr\'es}
\label{sec:reppst}

On d\'esigne par ${\bf Rep}_{{\BBq}_p}(G)$ la cat\'egorie des repr\'esentations $p$-adiques de $G$, c'est-\`a-dire des ${\BBq}_p$-espaces vectoriels de dimension finie munis d'une action lin\'eaire et continue de $G$. Pour pouvoir disposer de notions similaires \`a celles du cas $\ell \neq p$, on a besoin de la th\'eorie de Fontaine (voir \cite{Fo 2}), en particulier des anneaux $B_{dR}$, $B_{cris}$ et $B_{st}$ ; on renvoie \`a \cite{Fo 1} pour les d\'efinitions de ceux-ci.

Soient $L\subset \overline{{\BBq}}_p$ une extension finie de ${\BBq}_p$ et $L_0$ l'extension maximale non ramifi\'ee contenue dans $L$~; alors $(B_{cris})^{G_L}=(B_{st})^{G_L}=L_0$. Donc, si $V$ est une repr\'esentation $p$-adique de $G$, les objets $\mbox{Hom}_{{\BBq}_p[G_L]}(V,B_{st})$ et $\mbox{Hom}_{{\BBq}_p[G_L]}(V,B_{cris})$ sont des $L_0$-espaces vectoriels ; on montre que leur dimension est toujours inf\'erieure ou \'egale \`a $\mbox{dim}_{{\BBq}_p}(V)$. On dit que $V$ est semi-stable sur $L$ si $\mbox{dim}_{L_0}(\mbox{Hom}_{{\BBq}_p[G_L]}(V,B_{st})) = \mbox{dim}_{{\BBq}_p}(V)$ et que $V$ est cristalline sur $L$ si $\mbox{dim}_{L_0}(\mbox{Hom}_{{\BBq}_p[G_L]}(V,B_{cris})) = \mbox{dim}_{{\BBq}_p}(V)$.  On note ${\bf Rep}_{cris}(G)$, ${\bf Rep}_{st}(G)$, ${\bf Rep}_{cris,L}(G)$, ${\bf Rep}_{st,L}(G)$, ${\bf Rep}_{pcris}(G)$ et ${\bf Rep}_{pst}(G)$ les sous-cat\'egories pleines de ${\bf Rep}_{{\BBq}_p}(G)$ constitu\'ees des objets qui sont respectivement cristallins sur ${\BBq}_p$, semi-stables sur ${\BBq}_p$, cristallins sur $L$, semi-stables sur $L$, potentiellement cristallins et potentiellement semi-stables. 

\medskip
Soient $K\subset \overline{{\BBq}}_p$ une extension galoisienne de ${\BBq}_p$ de groupe de Galois $G_{K/{\BBq}_p}$ et $K_0$ l'extension maximale non ramifi\'ee contenue dans $K$ ; le Frobenius absolu $\sigma$ agit sur $K_0$. La cat\'egorie des $(\varphi,N,G_{K/{\BBq}_p})$-modules filtr\'es est d\'efinie de la mani\`ere suivante : \\
{\bf -} les objets sont des $K_0$-espaces vectoriels $D$ munis : 
\begin{itemize}
\item[(i)] d'une action $\sigma$-semi-lin\'eaire de $G_{K/{\BBq}_p}$ (le sous-groupe d'inertie agit lin\'eairement)
\item[(ii)] d'un Frobenius $\varphi : D\rightarrow D$, injectif, $\sigma$-semi-lin\'eaire et $G_{K/{\BBq}_p}$-\'equivariant 
\item[(iii)] d'un endomorphisme $K_0$-lin\'eaire $G_{K/{\BBq}_p}$-\'equivariant $N:D\rightarrow D$ tel que $N\varphi=p\varphi N$  
\item[(iv)] d'une filtration index\'ee par ${\BBz}$, d\'ecroissante, exhaustive et s\'epar\'ee sur $D_K=K\otimes_{K_0}D$ par des sous-$K$-espaces vectoriels $\{ \mbox{Fil}^i D_K\: , \: i\in {\BBz} \}$ stables par $G_{K/{\BBq}_p}$, l'action de $G_{K/{\BBq}_p}$ \'etant \'etendue semi-lin\'eairement sur $D_K$ 
\end{itemize}
{\bf -} un morphisme $f : D_1 \rightarrow D_2$ est une application $K_0$-lin\'eaire commutant \`a l'action de $G_{K/{\BBq}_p}$, \`a $\varphi$ et \`a $N$, et telle que, si l'on note $f_K$ l'application $K$-lin\'eaire d\'eduite de $f$ par extension des scalaires, $f_K(\mbox{Fil}^i D_{1,K})\subset \mbox{Fil}^i D_{2,K}$ pour tout $i\in {\BBz}$. 

\medskip
On d\'esigne par ${\bf MF}_{K/{\BBq}_p}(\varphi,N)$ la sous-cat\'egorie pleine des $(\varphi,N,G_{K/{\BBq}_p})$-modules filtr\'es form\'ee des objets sur lesquels l'action de $G_{K/{\BBq}_p}$ est discr\`ete et qui sont de dimension finie en tant que $K_0$-espace vectoriel ; le Frobenius $\varphi$ est alors bijectif et l'op\'erateur de monodromie $N$ est nilpotent. C'est une cat\'egorie ${\BBq}_p$-lin\'eaire mais non ab\'elienne, qui est munie d'un produit tensoriel (\cite{Fo 2}, 4.3.4) ; un sous-objet est un sous-$(\varphi,N,G_{K/{\BBq}_p})$-module muni de la filtration induite. On note ${\bf MF}_{K/{\BBq}_p}(\varphi)$ la sous-cat\'egorie pleine form\'ee des objets sur lesquels $N=0$ ; si $K={\BBq}_p$ on \'ecrit ${\bf MF}_{{\BBq}_p}(\varphi,N)$ et ${\bf MF}_{{\BBq}_p}(\varphi)$.

Soit $D$ un objet de dimension $d$ dans ${\bf MF}_{K/{\BBq}_p}(\varphi,N)$. On pose $t_H(D)=t_H(\wedge^d D)=\mbox{Max}\{\, i \in {\BBz}\: / \: \mbox{Fil}^i (\wedge^d D_K) \neq 0 \,\}$ et $t_N(D)=v_p(\lambda)$, o\`u $\lambda \in K_0$ est tel que $\varphi x = \lambda x$ pour un $x$ non nul de $\wedge^d D$. On dit que $D$ est faiblement admissible si $t_H(D)=t_N(D)$ et $t_H(D')\leq t_N(D')$ pour tout sous-objet $D'$ de $D$ (\cite{Fo 2}, 4.4.1). Un objet de ${\bf MF}_{K/{\BBq}_p}(\varphi,N)$ est faiblement admissible si et seulement si l'objet de ${\bf MF}_{K}(\varphi,N)$ obtenu en oubliant l'action de $G_{K/{\BBq}_p}$ l'est (\cite{Fo 2}, prop.4.4.9). On note ${\bf MF}_{K/{\BBq}_p}^{fa}(\varphi,N)$ la sous-cat\'egorie pleine de ${\bf MF}_{K/{\BBq}_p}(\varphi,N)$ form\'ee des objets faiblement admissibles ; d\'efinition similaire pour ${\bf MF}_{K/{\BBq}_p}^{fa}(\varphi)$.

Le type de Hodge-Tate d'un objet $D$ de dimension 2 de ${\bf MF}_{K/{\BBq}_p}(\varphi,N)$ est le couple d'entiers $(r,s)$ tel que $\mbox{Fil}^i D_K= D_K \Leftrightarrow i \leq r$ et $\mbox{Fil}^i D_K= 0 \Leftrightarrow i >s$.

\medskip
L'anneau $B_{st}$ est un $(\varphi,N,G)$-module filtr\'e que l'on peut construire \`a partir de $B_{cris}$ de la mani\`ere suivante (pour l'influence de ce choix voir \cite{Fo 2}, 5.2). Soit $\mbox{{\boldmath $\pi$}} = (\pi^{(n)}) \in (\overline{{\BBq}}_p)^{{\BBn}}$ tel que $\pi^{(0)}=p$ et $(\pi^{(n+1)})^p = \pi^{(n)}$, et soit $\mbox{\bf u} = \log([\mbox{{\boldmath $\pi$}}]/p) \in B_{dR}$ ; alors on prend $B_{st}=B_{cris}[\mbox{\bf u}] \subset B_{dR}$ sur lequel le Frobenius est \'etendu par $\varphi \mbox{\bf u} =p \mbox{\bf u}$ et $N$ est l'unique $B_{cris}$-d\'erivation telle que $N \mbox{\bf u} =1$. Les foncteurs contravariants 
$${\bf D}_{cris,K/{\BBq}_p}^* : {\bf Rep}_{cris,K}(G) \rightarrow {\bf MF}_{K/{\BBq}_p}^{fa}(\varphi) \makebox[1cm]{et}  {\bf D}_{st,K/{\BBq}_p}^* : {\bf Rep}_{st,K}(G) \rightarrow {\bf MF}_{K/{\BBq}_p}^{fa}(\varphi,N)$$
donn\'es par ${\bf D}_{cris,K/{\BBq}_p}^*(V)= \mbox{Hom}_{{\BBq}_p[G_K]}(V,B_{cris})$ et ${\bf D}_{st,K/{\BBq}_p}^*(V)= \mbox{Hom}_{{\BBq}_p[G_K]}(V,B_{st})$ \'etablissent une anti-\'equivalence de cat\'egories (\cite{Fo 2} et \cite{Co-Fo}). Les quasi-inverses sont donn\'es par ${\bf V}_{st,K/{\BBq}_p}^*(D)=\mbox{Hom}_{(\varphi,N,G)-mf}(D,B_{st})$ et ${\bf V}_{cris,K/{\BBq}_p}^*(D)=\mbox{Hom}_{(\varphi,G)-mf}(D,B_{cris})$, o\`u l'indice ``$mf$'' signifie ``modules filtr\'es'', $G$ op\'erant sur $D$ via son quotient $G_{K/{\BBq}_p}$.

Si $K={\BBq}_p$ on \'ecrit ${\bf D}_{cris}^*$ et ${\bf D}_{st}^*$. On note ${\bf D}_{pcris}^*$ et ${\bf D}_{pst}^*$ les foncteurs obtenus comme limite inductive des ${\bf D}_{cris,K/{\BBq}_p}^*$ et ${\bf D}_{st,K/{\BBq}_p}^*$ lorsque $K$ parcourt l'ensemble des extensions finies galoisiennes de ${\BBq}_p$ contenues dans $\overline{\BBq}_p$ ; ce sont des foncteurs de ${\bf Rep}_{pcris}(G)$ et ${\bf Rep}_{pst}(G)$ dans la limite inductive des ${\bf MF}_{K/{\BBq}_p}(\varphi,N)$.

\medskip
Enfin, on a un foncteur ${\bf WD}_{K/{\BBq}_p} : {\bf MF}_{K/{\BBq}_p}(\varphi,N) \rightarrow {\bf Rep}_{K_0}('W)$ obtenu en oubliant la filtration et en faisant agir le groupe de Weil $K_0$-lin\'eairement (\cite{Fo 3}) : si $D$ est un objet de ${\bf MF}_{K/{\BBq}_p}(\varphi,N)$ et si $D^{(0)}$ est le $(\varphi,N,G_{K/{\BBq}_p})$-module obtenu en oubliant la filtration, alors ${\bf WD}_{K/{\BBq}_p}(D)$ s'identifie au $K_0$-espace vectoriel $D^{(0)}$ muni de l'op\'erateur $N$ avec $\rho_0 (w) = (w \bmod W_K) \cdot \varphi^{-\upsilon (w)}$ pour tout $w \in W$, o\`u $W_K$ est le groupe de Weil relatif \`a $K$. Si $D_1$ et $D_2$ sont des objets de ${\bf MF}_{K/{\BBq}_p}(\varphi,N)$ on a un isomorphisme 
 $$K_0\!\otimes_{{\BBq}_p}\!\mbox{Hom}_{(\varphi,N,G_{K/{\BBq}_p})-mod} (D_1^{(0)},D_2^{(0)}) \ \simeq \ \mbox{Hom}_{{\bf Rep}_{K_0}('W)} ({\bf WD}_{K/{\BBq}_p}(D_1),{\bf WD}_{K/{\BBq}_p}(D_2))$$
de sorte que les classes d'isomorphisme de $(\varphi,N,G_{K/{\BBq}_p})$-modules sont en bijection avec celles de ${\bf Rep}_{K_0}('W)$. En composant avec ${\bf D}^*_{st,K/{\BBq}_p}$ on obtient un foncteur ${\bf W}^*_p : {\bf Rep}_{st,K}(G) \rightarrow {\bf Rep}_{K_0}('W)$ ; si $V$ est un objet de ${\bf Rep}_{st,K}(G)$ on dira que ${\bf W}^*_p(V)$ est la repr\'esentation de Weil-Deligne associ\'ee \`a $V$.

Posons ${\BBq}_{\ell}'={\BBq}_{\ell}$ si $\ell \neq p$, ${\BBq}_p'= K_0$, et choisissons des plongements de corps $\iota_{\ell} : {\BBq}_{\ell}' \hookrightarrow {\BBc}$ pour tout $\ell \in {\cal P}$. Un syst\`eme $(\Delta_{\ell})_{\ell \in {\cal P}}$ de repr\'esentations ${\BBq}_{\ell}'$-lin\'eaires de $'W$ est compatible si les $\Delta_{\ell}$ sont deux \`a deux compatibles lorsque $\ell$ parcourt ${\cal P}$. On a alors les m\^emes notions et crit\`eres que pour les syst\`emes $(\Delta_{\ell})_{\ell \neq p}$ en rempla\c cant \`a chaque fois ``$\ell \neq p$'' par ``$\ell \in {\cal P}$''.

\medskip
Soit $E/{\BBq}_p$ une courbe elliptique. La repr\'esentation $V_p(E)$ est potentiellement semi-stable~; si $L\subset \overline{{\BBq}}_p$ est une extension finie de ${\BBq}_p$, alors $E$ est semi-stable sur $L$ si et seulement si $V_p(E)$ l'est et $E$ a bonne r\'eduction sur $L$ si et seulement si $V_p(E)$ est cristalline sur $L$. De plus, si $K$ est la cl\^oture galoisienne d'un corps sur lequel $E$ devient semi-stable, on sait que $D={\bf D}^*_{st,K/{\BBq}_p}(V_p(E))$ est de type Hodge-Tate $(0,1)$, i.e. $\mbox{Fil}^i(D_K)=D_K$ pour $i\leq 0$, $\mbox{Fil}^1(D_K)$ est une $K$-droite et $\mbox{Fil}^i(D_K)=0$ pour $i\geq 2$. Enfin, on sait que la repr\'esentation de Weil-Deligne ${\bf W}^*_p(V_p(E))$ associ\'ee \`a $V_p(E)$ est $F$-semi-simple, d\'efinie sur ${\BBq}\,$, et que le syst\`eme $({\bf W}^*_{\ell}(V_{\ell}(E)))_{\ell \in {\cal P}}$ est compatible (cf. \cite{C-D-T} prop. B.4.2. pour les cas de potentielle bonne r\'eduction ; sinon, la pr\'esence d'un op\'erateur de monodromie non nul rend ces assertions faciles \`a v\'erifier).

\subsection{Notations}
\label{sec:nots}

\subsubsection{Quelques invariants de courbes elliptiques sur ${\BBq}_p$}
\label{sec:notsgeom}

Pour tout ce qui concerne les courbes elliptiques on peut se r\'ef\'erer aux livres de J.H. Silverman \cite{Si 1} et \cite{Si 2}.

Soit $E$ une courbe elliptique sur ${\BBq}_p$, i.e. une vari\'et\'e ab\'elienne sur ${\BBq}_p$ de dimension relative 1. Elle admet un mod\`ele sous forme de cubique plane, dit mod\`ele de Weierstrass, qui est donn\'e par une \'equation de la forme 
 $$ E \ : \ y^2 = x^3 + Ax +B \makebox[2cm]{avec} A,B \in {\BBq}_p \makebox[1cm]{et} 4A^3+27B^2 \neq 0 $$
On dispose d'abord d'un invariant $j_E = j(E) = 1728 \frac{4A^3}{4A^3+27B^2}$ (avec $1728=12^3$), dit invariant modulaire de $E$ ; il caract\'erise la classe d'isomorphisme de $E$ sur $\overline{{\BBq}}_p$. Le discriminant de $E$ est $\Delta_E = -16 (4A^3+27B^2) \neq 0$, et l'on a $j_E= -12^3(4A)^3 \Delta_E^{-1}$ ; le quotient $\Delta_E \bmod ({\BBq}_p^{\times})^{12}$ est un invariant de la classe d'isomorphisme de $E$ sur ${\BBq}_p$. Le mod\`ele de Weierstrass est minimal si $A,B \in {\BBz}_p$ et $0 \leq v_p(\Delta_E)< 12$ ; on peut toujours se ramener \`a un tel mod\`ele pour $E$.

On note $E(\overline{{\BBq}}_p)$ le groupe des points de $E$ \`a valeurs dans $\overline{{\BBq}}_p$ et $E[m]$ son sous-groupe de $m$-torsion, pour tout entier $m\geq 1$.

\medskip
Si $v_p(j_E)<0$ alors $E$ a potentiellement r\'eduction multiplicative d\'eploy\'ee : il existe un unique $q=q(j_E) \in p{\BBz}_p \backslash \{0\}$ v\'erifiant  
$$j_E = \frac{1}{q} + 744 + 196844 \, q + \ldots  $$
tel que $E$ est la tordue par un caract\`ere d'ordre $1$ ou $2$ d'une courbe de Tate $E_q$ ; cette torsion correspond \`a l'extension ${\BBq}_p(\sqrt{\gamma_E})/{\BBq}_p$ o\`u $\gamma_E= - 2AB^{-1} \bmod ({\BBq}_p^{\times})^2$ (\cite{Si 2},V, lemme 5.2). Le groupe $E_q(\overline{{\BBq}}_p)$ est isomorphe, en tant que groupe analytique rigide, \`a $\overline{{\BBq}}_p^{\times} / q^{{\BBz}}$. Pour tout $\ell \in {\cal P}$, on a une suite exacte courte de ${\BBq}_{\ell}[G]$-modules (voir \cite{Se 1}, A.1.2.)
 $$(*_m) \makebox[1.5cm]{} 0 \longrightarrow {\BBq}_{\ell}(1) \longrightarrow V_{\ell}(E_q) \longrightarrow {\BBq}_{\ell} \longrightarrow 0 $$
o\`u l'indice ``$m$'' signifie ``multiplicatif''.

\medskip
Si $v_p(j_E)\geq 0$, alors $E$ a potentiellement bonne r\'eduction : elle acquiert bonne r\'eduction sur une extension finie de ${\BBq}_p$. Le {\em d\'efaut de semi-stabilit\'e} de $E$ est l'indice de ramification minimal d'un corps sur lequel elle acquiert bonne r\'eduction ; on le note $\mbox{dst}(E)$. On a 
 $$\mbox{dst}(E) \ =\ \frac{12}{\mbox{pgcd}(12,v_p(\Delta_E))}$$
Si l'on choisit une \'equation de Weierstrass minimale pour $E$, alors $0\leq v_p(\Delta_E)< 12$ et $v_p(j_E)\geq 0$ impliquent que $v_p(\Delta_E)$ n'est pas premier \`a $12$~; on voit donc que $\mbox{dst}(E)=e=1,2,3,4,$ ou $6$ suivant que $v_p(\Delta_E)=0$, $v_p(\Delta_E)=6$, $v_p(\Delta_E)\in \{ 4,8 \}$, $v_p(\Delta_E)\in \{ 3,9 \}$, ou $v_p(\Delta_E)\in \{ 2,10 \}$ respectivement. Les entiers $e$ qui interviennent sont ceux qui v\'erifient $\mbox{{\boldmath $\varphi$}}(e) \in \{ 1,2 \}$ o\`u {\boldmath $\varphi$} est la fonction arithm\'etique d'Euler.

Comme $p\geq 5$, l'entier $e=\mbox{dst}(E)$ est premier \`a $p$ et $E$ acquiert bonne r\'eduction sur une extension totalement ramifi\'ee $L$ de ${\BBq}_p$ de degr\'e $e$ ; on note alors $\widetilde{E}_L/{\BBf}_p$ la fibre sp\'eciale du mod\`ele de N\'eron de $E\!\times_{{\BBq}_p}\!L$ et $a_p(\widetilde{E}_L)=a_p(E)$ la trace du polyn\^ome caract\'eristique du Frobenius arithm\'etique agissant sur $V_{\ell}(\widetilde{E}_L)$, $\ell \neq p$. On sait que $a_p(\widetilde{E}_L)$ est un entier rationnel ind\'ependant de $\ell$ tel que $\mid\!a_p(\widetilde{E}_L)\!\mid_{\infty} \,\leq 2\sqrt{p}$ et qu'il caract\'erise la classe d'isog\'enie sur ${\BBf}_p$ de $\widetilde{E}_L$ (\cite{Ta}) ; on a la relation 
 $$a_p(\widetilde{E}_L) = p + 1 - \# \widetilde{E}_L({\BBf}_p)$$
De plus, la courbe $\widetilde{E}_L$ est ordinaire si $p$ ne divise pas $a_p(\widetilde{E}_L)$ ; supersinguli\`ere si $p$ divise $a_p(\widetilde{E}_L)$, ce qui \'equivaut \`a $a_p(\widetilde{E}_L)=0$. 

Si $O_K$ (resp. $k$) est l'anneau des entiers d'une extension $K$ de ${\BBq}_p$ contenue dans $\overline{{\BBq}}_p$ (resp. une extension de ${\BBf}_p$ contenue dans $\overline{{\BBf}}_p$), on a un foncteur ${\cal E} \mapsto {\cal E}(p)$ (resp. $\widetilde{E} \mapsto \widetilde{E}(p)$) de la cat\'egorie des sch\'emas elliptiques sur $O_K$ (resp. des courbes elliptiques sur $k$) dans celle des groupes $p$-divisibles sur $O_K$ (resp. sur $k$). Si $E/K$ est une courbe elliptique ayant bonne r\'eduction sur $K$, on note $E(p)$ le groupe $p$-divisible sur $O_K$ du mod\`ele de N\'eron de $E$.

Si $E$ acquiert bonne r\'eduction ordinaire sur $L$, la partie connexe $E_L(p)^0$ de $E_L(p)$ est de hauteur $1$, et l'on a une suite exacte de groupes $p$-divisibles sur l'anneau des entiers de $L$  
$$0\longrightarrow E_L(p)^0 \longrightarrow E_L(p) \longrightarrow \widetilde{E}_L(p) \longrightarrow 0$$
qui induit la suite exacte courte de ${\BBq}_p[G]$-modules 
$$(*_{ord}) \makebox[1.5cm]{} 0 \longrightarrow V_p(E_L(p)^0) \longrightarrow V_p(E) \longrightarrow V_p(\widetilde{E}_L) \longrightarrow 0$$

\subsubsection{Notations galoisiennes}
\label{sec:notsgal}

On note ${\BBq}_{p^2}$ l'extension non ramifi\'ee de degr\'e 2 de ${\BBq}_p$. On choisit $\pi_{12}\in \overline{{\BBq}}_p$ v\'erifiant $(\pi_{12})^{12}+p=0$ et $\zeta_{12}\in \overline{{\BBq}}_p$ une racine primitive $12$-i\`eme de l'unit\'e. Pour tout entier $e\in \{ 1,2,3,4,6 \}$ on pose $\pi_e = (\pi_{12})^{12/e}$ et $\zeta_e = (\zeta_{12})^{12/e}$.

\medskip
On consid\`ere pour tout $e\in \{ 1,2,3,4,6 \}$ le corps $L_e={\BBq}_p(\pi_e)$ : c'est une extension totalement ramifi\'ee de degr\'e $e$ de ${\BBq}_p$ ; l'entier $e$ \'etant premier \`a $p$, cette extension est mod\'er\'ement ramifi\'ee. On note $K_e$ la cl\^oture galoisienne de $L_e$ dans $\overline{{\BBq}}_p$, $G_{K_e/{\BBq}_p}=\mbox{Gal}(K_e/{\BBq}_p)$ et $I_e=I_{K_e}$. Comme $({\BBz}/e{\BBz})^{\times}$ est d'ordre 1 ou 2, on a $p \equiv 1 \bmod e{\BBz}$ ou bien $p \equiv -1 \bmod e{\BBz}$. On se trouve alors dans l'une des situations suivantes~:\\
$K_1={\BBq}_p$ et $G_{K_1/{\BBq}_p}=1$ ; $K_2={\BBq}_p(\pi_2)$ et $G_{K_2/{\BBq}_p}=<\tau_2>$, o\`u $\tau_2$ est d\'efini par $\tau_2\pi_2=-\pi_2$ ;\\
si $e\in \{ 3,4,6 \}$  et $e\mid p-1$, $K_e={\BBq}_p(\pi_e)$ et $G_{K_e/{\BBq}_p}=<\tau_e>$, o\`u $\tau_e$ est d\'efini par $\tau_e\pi_e=\zeta_e\pi_e$~;\\
si $e\in \{ 3,4,6 \}$  et $e\mid p+1$, $K_e={\BBq}_{p^2}(\pi_e)={\BBq}_p(\pi_e,\zeta_e) $ et $G_{K_e/{\BBq}_p}=<\tau_e>\rtimes <\omega>$, o\`u $\tau_e$ est d\'efini par $\tau_e\pi_e=\zeta_e\pi_e$, $\tau_e \zeta_e=\zeta_e$, et $\omega$ est le rel\`evement du Frobenius absolu qui fixe $\pi_e$ et tel que $\omega \zeta_e ={\zeta_e}^{-1}$ ; on a $\omega\tau_e={\tau_e}^{-1}\omega$.

\medskip 
Il y a 3 extensions quadratiques de ${\BBq}_p$, une non ramifi\'ee et deux totalement ramifi\'ees. On les note $M_1 = {\BBq}_{p^2}$, $M_2 = {\BBq}_p(\pi_2)$ et $M_3$.

\medskip
On pose ${\cal N}_p= \{\, a \in {\BBz}\, / \mid\!a\!\mid_{\infty} \leq 2\sqrt{p}\, \}$. L'ensemble ${\cal N}_p^{\times} \subset {\BBz} \cap {\BBz}_p^{\times}$ des \'el\'ements non nuls de ${\cal N}_p$ est de cardinal $2[2\sqrt{p}\,]$ (partie enti\`ere). \\
On pose $\gamma_e= \zeta_e + \zeta_e^{-1} = -1,0,1$ pour $e=3,4,6$ respectivement, i.e. $X^2-\gamma_eX+1$ est le $e$-i\`eme polyn\^ome cyclotomique. \\
Quand $e \in \{3,4,6 \}$ et $e \mid p-1$, on note ${\cal N}_{p,e}^{\times}$ l'ensemble des $a \in {\BBz}$ tels que $(\gamma_e ^2 -4)(a^2 -4p)$ est un carr\'e dans ${\BBq}$ ; c'est un sous-ensemble de ${\cal N}_p^{\times}$. Si $a \in {\cal N}_{p,e}^{\times}$ le polyn\^ome $X^2-aX+p$ est scind\'e dans ${\BBz}_p[X]$ et il existe un unique $u_a \in {\BBz}_p^{\times}$ tel que $a = u_a + u_a^{-1}p$ ; alors on a ${\cal N}_{p,e}^{\times} = \{\, \pm (u_a \zeta_e^i + u_a^{-1}p\zeta_e^{-i}) \, , \, i\in {\BBz}/e{\BBz}\, \}$. L'ensemble ${\cal N}_{p,3}^{\times}= {\cal N}_{p,6}^{\times}= \{\, a \in {\BBz}\, /\, a^2 -4p \equiv -3 \ \bmod ({\BBq}^{\times})^2 \}$ est de cardinal 6 et l'ensemble ${\cal N}_{p,4}^{\times}= \{\, a \in {\BBz}\, /\, a^2 -4p \equiv -1 \ \bmod ({\BBq}^{\times})^2 \}$ est de cardinal 4. Par exemple : \\
${\cal N}_{5,4}^{\times}= \{ \pm 2, \pm 4 \} \subset {\cal N}_5^{\times}= \{ \pm 1, \pm 2, \pm 3, \pm 4\}$ ; ${\cal N}_{7,3}^{\times}= \{ \pm 1, \pm 4, \pm 5 \} \subset {\cal N}_7^{\times}= \{ \pm 1, \ldots , \pm 5 \}$ ;\\
${\cal N}_{13,3}^{\times}= \{ \pm 2, \pm 5, \pm 7 \}$ et ${\cal N}_{13,4}^{\times}= \{ \pm 4, \pm 6 \} \subset {\cal N}_{13}^{\times}= \{ \pm 1, \ldots , \pm 7\}$.

\section{Classification des ${\BBq}_{\ell}[G]$-modules $V_{\ell}(E)$}
\label{sec:classif}

\subsection{Les cas $\ell \neq p$}
\label{sec:casl}

\subsubsection{La liste ${\bf WD^*}$}
\label{sec:listl}

Soit $\phi \in W$ un rel\`evement du Frobenius g\'eom\'etrique modulo $I(\overline{{\BBq}}_p /{\BBq}_p(\pi_{12}))$~: $\phi$ agit trivialement sur $L_e$ pour tout $e\in \{ 1,2,3,4,6 \}$ et $\phi \bmod I_e$ agit par $x\mapsto x^{1/p}$ sur $\overline{{\BBf}}_p$. Pour $e\in \{ 2,3,4,6 \}$ on note $\theta_e$ un rel\`evement dans $I$ de $\tau_e \in I/I_e=I(K_e/{\BBq}_p)$. \\
Soient $a \in {\BBz}_p^{\times}$ et $u_a \in {\BBz}_p^{\times}$ tels que $X^2-aX +p = (X-u_a)(X-u_a^{-1}p)$. Si $e\in \{3,4,6\}$ divise $p-1$ on pose pour $\epsilon \in \{ \pm 1 \}$
 $$t_{\epsilon}=t_{\epsilon}(a,e)= u_a \zeta_e^{\epsilon} + u_a^{-1}p\, \zeta_e^{-\epsilon} \in {\BBz}_p^{\times}$$
On a $(X-t_1)(X-t_{-1})= X^2- \gamma_e a X+ p{\gamma_e}^2 + {a}^2 -4p = T(X)$ et $t_{1}\neq t_{-1}$. La condition $a \in {\cal N}^{\times}_{p,e}$ signifie exactement que les racines de $T(X)$ sont dans ${\BBz}$ ; elles sont alors dans ${\cal N}^{\times}_{p,e}$.

\medskip
Soit $F$ un corps de caract\'eristique $0$. La liste ${\bf WD^*}$ suivante d\'efinit \`a isomorphisme pr\`es des objets de ${\bf Rep}_{F}('W)$ de dimension deux, qui sont dans ${\bf Rep}^{\circ}_{{\BBq}_{\ell}}('W)$ lorsque $F={\BBq}_{\ell}$ :

\medskip
${\bf WD^*_m(e;b)}$, $e \in \{ 1,2 \}$, $b\in \{ -1,1 \}$ : \\
${\rho}_0(I_2)=1$ ; ${\rho}_0(\theta_2)=(-1)^{e-1}$ ; ${\mbox{P}}_{min}({\rho}_0(\phi))=(X-b)(X-bp)$ ; ${\mbox{P}}_{min}(N)=X^2$ ; ${\rho}_0(\phi)N=p^{-1}N{\rho}_0(\phi)$.

\medskip
${\bf WD^*_c(e;a_p)}$, $e\in \{1,2\}$, $a_p \in {\cal N}_p$ : \\
${\rho}_0(I_2)=1$ ; ${\rho}_0(\theta_2)=(-1)^{e-1}$ ; ${\mbox{P}}_{min}({\rho}_0(\phi))=X^2-a_pX+p$ ; $N=0$.

\medskip
${\bf WD^*_{pc}(e;a_p;\epsilon)}$, $e\in \{3,4,6\}$ et $e\mid p-1$, $a_p \in {\cal N}^{\times}_{p,e}$, $\epsilon \in \{-1,1\}$ : \\
${\rho}_0(I_e)=1$ ; ${\mbox{P}}_{min}({\rho}_0(\theta_e))=X^2-\gamma_eX+1$ ; ${\mbox{P}}_{min}({\rho}_0(\phi))=X^2-a_pX+p$ ; ${\mbox{P}}_{min}({\rho}_0(\phi){\rho}_0(\theta_e)) =X^2-t_{\epsilon}X+p$ ; ${\rho}_0(\phi){\rho}_0(\theta_e)={\rho}_0(\theta_e){\rho}_0(\phi)$ ; $N=0$.

\medskip
${\bf WD^*_{pc}(e;0)}$, $e\in \{ 3,4,6 \}$ et $e\mid p+1$ :\\
${\rho}_0(I_e)=1$ ; ${\mbox{P}}_{min}({\rho}_0(\theta_e))=X^2-\gamma_eX+1$ ; ${\mbox{P}}_{min}({\rho}_0(\phi))=X^2+p$ ; ${\rho}_0(\phi){\rho}_0(\theta_e)={\rho}_0(\theta_e)^{-1}{\rho}_0(\phi)$~; $N=0$.

\begin{rem}
{\em Les classes d'isomorphisme de ces objets ne d\'ependent pas du choix des corps $K_e$ : si dans la description de l'un d'eux on remplace $K_e$ par une autre extension galoisienne de ${\BBq}_p$ d'indice de ramification $e$, on obtient un objet isomorphe.}
\end{rem}

Tous les objets $\Delta$ de la liste ${\bf WD^*}$ sont d\'efinis sur ${\BBq}\,$. De plus, la repr\'esentation $F$-lin\'eaire de $'W$ de dimension un $\wedge^2 \Delta$ est donn\'ee par : $\wedge^2{\rho}_0(I)=1$, $\wedge^2{\rho}_0(\phi)=p$, $\wedge^2N=0$ ; si $F={\BBq}_{\ell}$ l'objet obtenu en appliquant le foncteur quasi-inverse est ${\BBq}_{\ell}(1)$. \\
Les objets du type ${\bf WD^*_m}$ sont des tordus d'objets semi-stables sur ${\BBq}_p$ mais n'ont pas potentiellement bonne r\'eduction. Les objets du type ${\bf WD^*_c}$ sont des tordus d'objets ayant bonne r\'eduction sur ${\BBq}_p$. Les objets du type ${\bf WD^*_{pc}}$ ont potentiellement bonne r\'eduction mais ne sont pas des tordus d'objets ayant bonne r\'eduction sur ${\BBq}_p$.

\subsubsection{Description des twists quadratiques}
\label{sec:twistsl}
  
La liste ${\bf WD^*}$ est stable par twists quadratiques. En tordant un objet $\Delta$ de ${\bf WD^*}$ par l'un des caract\`eres quadratiques on obtient les objets $\Delta_1$, $\Delta_2$ et $\Delta_3$ correspondant respectivement aux extensions $M_1$, $M_2$ et $M_3$. En faisant varier $\Delta$ parmi les objets de la liste on obtient : 
 
\smallskip
\noindent $\Delta = {\bf WD^*_m(1;b)}$ ; $\Delta_1 = {\bf WD^*_m(1;-b)}$ ; $\Delta_2 = {\bf WD^*_m(2;b)}$ ; $\Delta_3 = {\bf WD^*_m(2;-b)}$ \\
$\Delta = {\bf WD^*_c(1;a_p)}$ ; $\Delta_1 = {\bf WD^*_c(1;-a_p)}$ ; $\Delta_2 = {\bf WD^*_c(2;a_p)}$ ; $\Delta_3 = {\bf WD^*_c(2;-a_p)}$ \\
$\Delta = {\bf WD^*_{pc}(4;a_p;\epsilon)}$ ; $\Delta_1 = {\bf WD^*_{pc}(4;-a_p;\epsilon)}$ ; $\Delta_2 = {\bf WD^*_{pc}(4;a_p;-\epsilon)}$ ; $\Delta_3 = {\bf WD^*_{pc}(4;-a_p;-\epsilon)}$ \\
$\Delta = {\bf WD^*_{pc}(3;a_p;\epsilon)}$ ; $\Delta_1 = {\bf WD^*_{pc}(3;-a_p;\epsilon)}$ ; $\Delta_2 = {\bf WD^*_{pc}(6;a_p;-\epsilon)}$ ; $\Delta_3 = {\bf WD^*_{pc}(6;-a_p;-\epsilon)}$  \\
$\Delta = {\bf WD^*_{pc}(4;0)} = \Delta_1 = \Delta_2 = \Delta_3$ \\
$\Delta = {\bf WD^*_{pc}(3;0)} =\Delta_1$ ; $\Delta_2 = {\bf WD^*_{pc}(6;0)} =\Delta_3$

\smallskip 
Si un objet $\Delta$ de la liste ${\bf WD^*}$ provient d'une courbe elliptique $E/{\BBq}_p$ alors les $\Delta_i$, $i \in \{ 1,2,3 \}$, proviennent des courbes elliptiques $E_i/{\BBq}_p$ obtenues en tordant $E$ sur $M_i$.

\subsubsection{Classification}
\label{sec:classl}

\begin{prop}
\label{classlprop}
Soit $\ell \in {\cal P}$ tel que $\ell \neq p$. Les repr\'esentations de la liste ${\bf WD^*}$ sont deux \`a deux non isomorphes. Si $E$ est une courbe elliptique sur ${\BBq}_p$ alors ${\bf W}^*_{\ell}(V_{\ell}(E))$ est isomorphe \`a l'un des objets de la liste ${\bf WD^*}$.
\end{prop}

\begin{pf}
Soit $E/{\BBq}_p$ une courbe elliptique et soit $(\Delta_{\ell},{\rho}_0,N)={\bf W}^*_{\ell}(V_{\ell}(E))$. 

Supposons $v_p(j_E)< 0$. Si $E$ est une courbe de Tate alors $\Delta_{\ell} \simeq {\bf WD^*_m(1;1)}$, et en tordant par les trois caract\`eres quadratiques (cf.~\ref{sec:twistsl}), on obtient les ${\bf WD^*_m(e;b)}$, voir \cite{Ro}, $\S$ 15.

Supposons $v_p(j_E)\geq 0$. Alors $E$ acquiert bonne r\'eduction sur $L_e$ avec $e=\mbox{dst}(E)$ et l'on a $N=0$, ${\rho}_0(I_e)=1$ et $\mbox{P}_{min}({\rho}_0(\phi))(X)=X^2-a_pX+p$ avec $a_p=a_p(\widetilde{E}_{L_e})\in {\cal N}_p$. La minimalit\'e de $e$ implique que ${\rho}_0 \bmod I_e : <\tau_e>=I/I_e \hookrightarrow \mbox{Aut}_{{\BBq}_{\ell}}({\Delta}_{\ell})$ est injective.

Si $e\in \{1,2\}$ alors $\Delta_{\ell} \simeq {\bf WD_c^*(e;a_p)}$. Supposons $e \in \{ 3,4,6 \}$, d'o\`u $({\BBz}/e{\BBz})^{\times}=\{ \pm 1 \}$. Alors $\mbox{P}_{min}({\rho}_0(\theta_e))(X) = (X-\zeta_e)(X-\zeta_e^{-1}) = X^2 -\gamma_e X +1$, puisque ${\rho}_0(\theta_e)$ est d'ordre $e$ et de d\'eterminant $1$. Comme $\theta_e \phi \equiv \phi \theta_e^p \bmod I_e$, on a ${\rho}_0(\phi){\rho}_0(\theta_e)= {\rho}_0(\theta_e){\rho}_0(\phi)$ si $p \equiv 1 \bmod e{\BBz}$ et ${\rho}_0(\phi){\rho}_0(\theta_e)= {\rho}_0(\theta_e)^{-1}{\rho}_0(\phi)$ si $p \equiv -1 \bmod e{\BBz}$. \\
Si $e\mid p+1$, en se pla\c cant dans ${\BBq}_{\ell}(\zeta_e)\!\otimes_{\!{\BBq}_{\ell}}\!\Delta_{\ell}$ on voit que $\mbox{Tr}({\rho}_0(\phi))= a_p =\mbox{Tr}({\rho}_0(\phi \theta_e))= 0$ et $\widetilde{E}_{L_e}$ est supersinguli\`ere. Dans ce cas la donn\'ee de $\mbox{Tr}({\rho}_0(\theta_e))$ et $\mbox{Tr}({\rho}_0(\phi))$ suffit pour d\'eterminer la classe de $\Delta_{\ell}$, qui correspond \`a ${\bf WD_{pc}^*(e;0)}$. \\
Si $e\mid p-1$ la repr\'esentation est ab\'elienne, donc sa classe est d\'etermin\'ee par la donn\'ee de $\mbox{Tr}({\rho}_0(\phi))$, $\mbox{Tr}({\rho}_0(\theta_e))$ et $\mbox{Tr}({\rho}_0(\phi \theta_e))$. Prenons une ${\BBq}_{\ell}(\zeta_e)$-base de ${\BBq}_{\ell}(\zeta_e)\!\otimes_{\!{\BBq}_{\ell}}\!\Delta_{\ell}$ dans laquelle la matrice de ${\rho}_0(\theta_e)$ s'\'ecrit $\mbox{Diag}(\zeta_e^{\epsilon},\zeta_e^{-\epsilon})$ avec $\epsilon \in \{ \pm 1 \}$ et celle de ${\rho}_0(\phi)$ s'\'ecrit $\mbox{Diag}(z_1 ,z_2)$. Alors $\mbox{Tr}({\rho}_0(\phi \theta_e)) = t_{\epsilon} = \zeta_e^{\epsilon}z_1 +\zeta_e^{-\epsilon}z_2$ est racine de $T(X) = X^2 - \gamma_e a_p X + p\gamma_e^2 + a_p^2 -4p$ dont le discriminant est $(a_p^2 -4p)(\gamma_e^2 -4) \neq 0$. Le fait que $\Delta_{\ell}$ est d\'efinie sur ${\BBq}$ implique $t_{\epsilon} \in {\BBq}$ ce qui \'equivaut \`a $a_p \in {\cal N}_{p,e}^{\times}$ ; en particulier $a_p\neq 0$ et $\widetilde{E}_{L_e}$ est ordinaire. Finalement $\Delta_{\ell} \simeq {\bf WD_{pc}^*(e;a_p;\epsilon)}$.
\end{pf}

\begin{rem}
{\em Soit $E/{\BBq}_p$ telle que $v_p(j_E)<0$, d'o\`u ${\bf W}^*_{\ell}(V_{\ell}(E))\simeq {\bf WD^*_m(e;b)}$ ; soit $\gamma_E \in {\BBq}_p^{\times}/({\BBq}_p^{\times})^2$ l'invariant d\'efini en \ref{sec:notsgeom}. Alors on a : $(e;b)=(1;1)\Leftrightarrow {\BBq}_p(\sqrt{\gamma_E})={\BBq}_p$ ; $(e;b)=(1;-1)\Leftrightarrow {\BBq}_p(\sqrt{\gamma_E})=M_1$ ; $(e;b)=(2;1)\Leftrightarrow {\BBq}_p(\sqrt{\gamma_E})=M_2$ ; $(e;b)=(2;-1)\Leftrightarrow {\BBq}_p(\sqrt{\gamma_E})=M_3$.}
\end{rem}

\begin{rem}
{\em Soit $E/{\BBq}_p$ telle que $v_p(j_E)\geq 0$ et $\mbox{dst}(E)=e \geq 3$. On voit que si $e \mid p-1$, alors $E$ est potentiellement ordinaire, et si $e \mid p+1$, alors $E$ est potentiellement supersinguli\`ere.}
\end{rem}

\begin{rem}
{\em Dans la m\^eme situation que ci-dessus, on d\'etermine l'invariant $\epsilon \in \{ \pm 1 \}$ qui intervient lorsque $e$ divise $p-1$ en \'etudiant le ${\BBf}_p[I]$-module $E[p]$ : si l'on prend une \'equation de Weierstrass minimale pour $E$ alors on a $\epsilon=1$ si $v_p(\Delta_E) < 6$ (i.e. $v_p(\Delta_E)\in \{2,3,4\}$) et $\epsilon=-1$ si $v_p(\Delta_E) > 6$ (i.e. $v_p(\Delta_E)\in \{8,9,10\}$) (\cite{Kr}, 2.3.1, prop.1).}
\end{rem}

\begin{rem}
\label{lcomp}
{\em Finalement, on constate que si l'on se donne une courbe elliptique $E/{\BBq}_p$ sous forme d'\'equation de Weierstrass, une liste d'invariants de la courbe (explicitement calculables) suffit pour retrouver la classe d'isomorphisme de $V_{\ell}(E)$, $\ell \neq p$. En particulier, lorsque $E/{\BBq}_p$ est fix\'ee et que $\ell$ parcourt ${\cal P} \backslash \{ p\}$, les classes des $V_{\ell}(E)$ sont ind\'ependantes de $\ell$ : ceci exprime la compatibilit\'e au sens de Weil-Deligne du syst\`eme de repr\'esentations $(V_{\ell}(E))_{\ell \neq p}$.}
\end{rem}

\subsection{Le cas $\ell = p$}
\label{sec:casp}

\subsubsection{La liste ${\bf D^*}$}
\label{sec:listp}

On d\'efinit la liste ${\bf D^*}$ d'objets de ${\bf MF}_{K_e/{\BBq}_p}(\varphi,N)$ et de type Hodge-Tate $(0,1)$ suivants :

\medskip
${\bf D^*_m(e;b;\alpha)}$, $e\in \{ 1,2 \} $, $b\in \{-1,1 \} $, $\alpha\in{\BBq}_p$ :\\
Pour $e=1$ : $D={\BBq}_pe_1\oplus{\BBq}_pe_2$ avec $\varphi e_1 = be_1$, $\varphi e_2 = pbe_2$ ; $Ne_1 = 0$, $N e_2 = e_1$ ; $\mbox{Fil}^1 D = (\alpha e_1 + e_2){\BBq}_p$.\\
Pour $e=2$ : $G_{K_e/{\BBq}_p}=<\tau_2>$, $D={\BBq}_pe_1\oplus{\BBq}_pe_2$ avec $\varphi e_1 = be_1$, $\varphi e_2 = pbe_2$ ; $Ne_1 = 0$, $N e_2 = e_1$ ; $\tau_2e_1=-e_1$, $\tau_2e_2=-e_2$ ; $\mbox{Fil}^1 D_{K_e} =(\alpha\cdot e_1\otimes 1 + e_2\otimes 1){\BBq}_p(\pi_2)$.

\medskip
${\bf D^*_c(e;a_p;\alpha)}$, $e\in \{ 1,2 \} $, $a_p \in {\cal N}_p^{\times}$, $\alpha \in \{ 0,1 \}$ :\\
Soit $u\in {\BBz}_p^{\times}$ tel que $u + u^{-1}p = a_p$.\\
Pour $e=1$ : $D={\BBq}_pe_1\oplus{\BBq}_pe_2$ avec $\varphi e_1=ue_1$, $\varphi e_2 =u^{-1} pe_2$ ; $N e_1 =N e_2=0$~; $\mbox{Fil}^1D=(\alpha e_1 + e_2){\BBq}_p$.\\
Pour $e=2$ : $G_{K_e/{\BBq}_p}=<\tau_2>$, $D={\BBq}_pe_1\oplus{\BBq}_pe_2$ avec $\varphi e_1 =ue_1$, $\varphi e_2=u^{-1}pe_2$ ; $Ne_1=Ne_2=0$ ; $\tau_2e_1=-e_1$, $\tau_2e_2 =-e_2$ ; $\mbox{Fil}^1D_{K_e}=(\alpha\cdot e_1\otimes 1 + e_2\otimes 1){\BBq}_p(\pi_2)$.

\medskip
${\bf D^*_c(e;0)}$, $e\in \{ 1,2 \}$ :\\
Pour $e=1$ : $D={\BBq}_pe_1\oplus{\BBq}_pe_2$ avec $\varphi e_1=e_2$, $\varphi e_2=- pe_1$ ; $Ne_1=Ne_2=0$~; $\mbox{Fil}^1D= e_1{\BBq}_p$. \\
Pour $e=2$ : $G_{K_e/{\BBq}_p}=<\tau_2>$, $D={\BBq}_pe_1\oplus{\BBq}_pe_2$ avec $\varphi e_1=e_2$, $\varphi e_2=-pe_1$ ; $Ne_1=Ne_2=0$~; $\tau_2e_1=-e_1$, $\tau_2e_2=-e_2$ ; $\mbox{Fil}^1D_{K_e}=(e_1\otimes 1){\BBq}_p(\pi_2)$.

\medskip
${\bf D^*_{pc}(e;a_p;\epsilon;\alpha)}$, $e\in \{ 3,4,6 \}$ et $e\mid p-1$, $a_p\in {\cal N}^{\times}_{p,e}$, $\epsilon \in \{ -1,1 \}$, $\alpha\in \{ 0,1 \}$ :\\
Soit $u\in {\BBz}_p^{\times}$ tel que $u + u^{-1}p = a_p$.\\
$G_{K_e/{\BBq}_p}=<\tau_e>$, $D={\BBq}_pe_1\oplus{\BBq}_pe_2$ avec $\varphi e_1=ue_1$, $\varphi e_2=u^{-1}pe_2$ ; $Ne_1=Ne_2= 0$ ; $\tau_e e_1= {\zeta_e}^{\epsilon}e_1$, $\tau_e e_2={\zeta_e}^{-\epsilon}e_2$ ; $\mbox{Fil}^1D_{K_e}=(\alpha\cdot e_1\otimes {\pi_e}^{-\epsilon}+e_2\otimes {\pi_e}^{\epsilon}){\BBq}_p(\pi_e)$. 

\medskip
${\bf D^*_{pc}(e;0;\alpha)}$, $e\in \{ 3,4,6 \}$ et $e\mid p+1$, $\alpha\in {\BBp}^1({\BBq}_p)$ :\\
$G_{K_e/{\BBq}_p}=<\tau_e>\rtimes <\omega>$, $D={\BBq}_{p^2}e_1\oplus{\BBq}_{p^2}e_2$ avec $\varphi e_1=e_2$, $\varphi e_2 =-pe_1$~; $Ne_1=Ne_2=0$ ; $\omega e_1=e_1$, $\omega e_2 =e_2$ ; $\tau_e e_1={\zeta_e}e_1$, $\tau_e e_2={\zeta_e}^{-1}e_2$ ; $\mbox{Fil}^1D_{K_e}=(\alpha\cdot e_1\otimes {\pi_e}^{-1}+ e_2\otimes {\pi_e}){\BBq}_{p^2}(\pi_e)$.

\medskip
La classe d'isomorphisme de ces objets est ind\'ependante du choix fait pour l'extension galoisienne $K_e$ (elle ne d\'epend que de l'indice de ramification $e$).

Tous les objets $D$ de la liste ${\bf D^*}$ sont faiblement admissibles. Pour chacun d'entre eux la repr\'esentation de Weil-Deligne associ\'ee est d\'efinie sur ${\BBq}\,$. De plus, $\wedge^2 D$ est un objet de ${\bf MF}_{{\BBq}_p}(\varphi)$ de dimension un d\'ecrit par : $\varphi = p$ et $\mbox{Fil}^1 (\wedge^2 D) = \wedge^2 D$, $\mbox{Fil}^2 (\wedge^2 D) = 0$ ; l'objet obtenu en appliquant le foncteur quasi-inverse est ${\BBq}_p(1)$. \\
Les objets du type ${\bf D^*_m}$ sont des tordus d'objets semi-stables sur ${\BBq}_p$ mais ne sont pas potentiellement cristallins. Les objets du type ${\bf D^*_c}$ sont des tordus d'objets cristallins sur ${\BBq}_p$. Les objets du type ${\bf D^*_{pc}}$ sont potentiellement cristallins mais ne sont pas des tordus d'objets cristallins sur ${\BBq}_p$.

\subsubsection{Description des twists quadratiques}
\label{sec:twistsp}
 
La liste ${\bf D^*}$ est stable par twists quadratiques. En tordant un objet $D$ de ${\bf D^*}$ par l'un des caract\`eres quadratiques on obtient les objets $D_1$, $D_2$ et $D_3$ correspondant respectivement aux extensions $M_1$, $M_2$ et $M_3$. En faisant varier $D$ parmi les objets de la liste on obtient : 

\smallskip
\noindent $D = {\bf D^*_m(1;b;\alpha)}$ ; $D_1 = {\bf D^*_m(1;-b;\alpha)}$ ; $D_2 = {\bf D^*_m(2;b;\alpha)}$ ; $D_3 = {\bf D^*_m(2;-b;\alpha)}$ \\
$D = {\bf D^*_c(1;a_p;\alpha)}$ ; $D_1 = {\bf D^*_c(1;-a_p;\alpha)}$ ; $D_2 = {\bf D^*_c(2;a_p;\alpha)}$ ; $D_3 = {\bf D^*_c(2;-a_p;\alpha)}$ \\
$D ={\bf D^*_c(1;0)}$ ; $D_1 = D $ ; $D_2 = D_3 = {\bf D^*_c(2;0)}$ \\
$D = {\bf D^*_{pc}(4;a_p;\epsilon;\alpha)}$ ; $D_1 = {\bf D^*_{pc}(4;-a_p;\epsilon;\alpha)}$ ; $D_2 = {\bf D^*_{pc}(4;a_p;-\epsilon;\alpha)}$ ; $D_3 = {\bf D^*_{pc}(4;-a_p;-\epsilon;\alpha)}$ \\
$D = {\bf D^*_{pc}(3;a_p;\epsilon;\alpha)}$ ; $D_1 = {\bf D^*_{pc}(3;-a_p;\epsilon;\alpha)}$ ; $D_2 = {\bf D^*_{pc}(6;a_p;-\epsilon;\alpha)}$ ; $D_3 = {\bf D^*_{pc}(6;-a_p;-\epsilon;\alpha)}$ \\
$D = {\bf D^*_{pc}(4;0;\alpha)}$ ; $D_1 = {\bf D^*_{pc}(4;0;-\alpha)}$ ; $D_2 = {\bf D^*_{pc}(4;0;p^2 \alpha^{-1})}$ ; $D_3 = {\bf D^*_{pc}(4;0;- p^2 \alpha^{-1})}$ \\
$D = {\bf D^*_{pc}(3;0;\alpha)}$ ; $D_1 = {\bf D^*_{pc}(3;0;-\alpha)}$ ; $D_2 = {\bf D^*_{pc}(6;0;p^2 \alpha^{-1})}$ ; $D_3 = {\bf D^*_{pc}(6;0;- p^2 \alpha^{-1})}$ 

\smallskip 
Si un objet $D$ de ${\bf D^*}$ provient d'une courbe elliptique $E/{\BBq}_p$ alors les $D_i$, $i \in \{ 1,2,3 \}$, proviennent des courbes elliptiques $E_i/{\BBq}_p$ obtenues en tordant $E$ sur $M_i$.

\subsubsection{L'image de Galois dans $\mbox{Aut}_{{\BBq}_p}(V)$}
\label{sec:imgal}

Les objets ${\bf D^*_c(e;0)}$ et ${\bf D^*_{pc}(e;0;\alpha)}$ sont irr\'eductibles. Tous les autres objets $D$ de la liste ${\bf D^*}$ sont r\'eductibles et il est possible de d\'ecrire l'action de $G$ sur le semi-simplifi\'e de $V \simeq {\bf V}^*_{pst}(D)$ (\`a comparer avec le lemme 7.1.3 de \cite{C-D-T}). Soit $\chi$ le caract\`ere cyclotomique $G \rightarrow {\BBz}_p^{\times}$. Pour tout $u \in {\BBz}_p^{\times}$ on note $\eta_u : G \rightarrow G/I \rightarrow {\BBz}_p^{\times}$ le caract\`ere non ramifi\'e qui envoie le Frobenius arithm\'etique sur $u$. Lorsque $e\geq 2$ et $e\mid p-1$ on note $\xi_e : G \rightarrow G_{K_e/{\BBq}_p} \rightarrow <\zeta_e> \subset {\BBz}_p^{\times}$ le caract\`ere ramifi\'e d\'efini par $\xi_e(g)= g\pi_e/ \pi_e$, $g\in G$ ; on a $\xi_e = [\chi \bmod p{\BBz}_p]^{\frac{p-1}{e}}$, o\`u $[-]$ est le repr\'esentant de Teichm\"uller.

\medskip
${\bf D^*_m(e;b;\alpha)}\simeq {\bf D}^*_{pst}(V)$, $e \in \{ 1,2 \}$, $b\in \{ -1,1 \}$, $\alpha\in{\BBq}_p$ : il existe une ${\BBq}_p$-base de $V$ telle que $G$ agit via 
$$ \left(  \begin {array}{cc}
            \eta_{-1}^{\frac{1-b}{2}} \xi_2^{e-1} \chi   &  * \\
                 0                                       &    \eta_{-1}^{\frac{b-1}{2}}\xi_2^{1-e}
            \end{array}
   \right)
\makebox[2cm]{avec} *\neq 0 $$

\smallskip
${\bf D^*_c(e;a_p;\alpha)}\simeq {\bf D}^*_{pst}(V)$, $e\in \{ 1,2 \}$, $a_p \in {\cal N}_p^{\times}$, $\alpha \in \{ 0,1 \}$~: soit $u \in {\BBz}_p^{\times}$ tel que $u+u^{-1}p=a_p$ ; il existe une ${\BBq}_p$-base de $V$ telle que $G$ agit via 
$$ \left(  \begin {array}{cc}
            \eta_{u}^{-1} \xi_2^{e-1} \chi       &  * \\
                 0                               &    \eta_{u} \xi_2^{e-1}
            \end{array}
   \right)
\makebox[2cm]{avec} *=0 \; \Leftrightarrow \; \alpha=0 $$

\smallskip
${\bf D^*_{pc}(e;a_p;\epsilon;\alpha)}\simeq {\bf D}^*_{pst}(V)$, $e\in \{ 3,4,6 \}$ et $e\mid p-1$, $a_p\in {\cal N}^{\times}_{p,e}$, $\epsilon \in \{ -1,1 \}$, $\alpha\in \{ 0,1 \}$~: soit $u \in {\BBz}_p^{\times}$ tel que $u+u^{-1}p=a_p$ ; il existe une ${\BBq}_p$-base de $V$ telle que $G$ agit via 
$$ \left(  \begin {array}{cc}
            \eta_{u}^{-1} \xi_e^{-\epsilon} \chi       &  * \\
                 0                               &    \eta_{u} \xi_e^{\epsilon}
            \end{array}
   \right)
\makebox[2cm]{avec} *=0 \; \Leftrightarrow \; \alpha=0 $$

\subsubsection{Classification}
\label{sec:classp}

\begin{prop}
\label{classpprop}
Les objets de la liste ${\bf D^*}$ sont deux \`a deux non isomorphes. Si $E$ est une courbe elliptique sur ${\BBq}_p$ alors ${\bf D}^*_{pst}(V_p(E))$ est isomorphe \`a l'un des objets de la liste ${\bf D^*}$.
\end{prop}

\begin{rem}
{\em Dans \cite{Fo-Ma}, J.-M. Fontaine et B. Mazur classifient les repr\'esentations potentiellement semi-stables et faiblement admissibles de $G$ sur un ${\BBq}_p$-espace vectoriel de dimension 2. Les objets de dimension 1 sont d\'ecrits au $\S$ 8, tandis que les objets de dimension 2 qui ne sont pas somme directe d'objets de dimension 1 sont d\'ecrits dans la liste du d\'ebut du $\S$ 11. Comme $12$ divise $p^2-1$ pour $p\geq 5$, les $K_e$, $e\in \{2,3,4,6\}$, sont des sous-corps du corps not\'e $F_2$ dans \cite{Fo-Ma}. Les correspondances dans les notations sont : 

\smallskip
\noindent ${\bf D^*_m(1;b;\alpha)}= D_{II}(0,1;b,\alpha \,;0)$, $b\in\{-1,1\}$, $\alpha\in{\BBq}_p$ \\
${\bf D^*_c(e;a_p;0)}=U_1\oplus U_2$, $e\in\{1,2\}$, $a_p \in {\cal N}_p^{\times}$ : cf.~\ref{sec:imgal} \\
${\bf D^*_c(1;a_p;1)}=D_I(0,1;a_p,p;0)$, $a_p \in {\cal N}_p^{\times}$  \\
${\bf D^*_c(2;a_p;1)}=D_I(0,1;a_p,p;\frac{p-1}{2})$, $a_p \in {\cal N}_p^{\times}$  \\
${\bf D^*_c(1;0)}=D_I(0,1;0,p;0)$  \\
${\bf D^*_c(2;0)}=D_I(0,1;0,p;\frac{p-1}{2})$  \\
${\bf D^*_{pc}(e;a_p;\epsilon;0)}=U_1\oplus U_2$, $e\in \{ 3,4,6 \}$, $e\mid p-1$, $a_p\in {\cal N}^{\times}_{p,e}$, $\epsilon \in \{ \pm 1\}$ : cf.~\ref{sec:imgal} \\
${\bf D^*_{pc}(e;a_p;\epsilon;1)}=D_{III}(0,1;u,u^{-1}p;\epsilon\, \frac{p-1}{e},-\epsilon\, \frac{p-1}{e})$, $e\in \{ 3,4,6 \}$ et $e\mid p-1$, $u\in {\BBz}_p^{\times}$ tel que $u + u^{-1}p = a_p \in {\cal N}^{\times}_{p,e}$, $\epsilon \in \{ \pm 1\}$ \\
${\bf D^*_{pc}(e;0;\alpha)}=D_{IV}(0,1;p; \frac{p+1}{e} -1 , p-\frac{p+1}{e};p^{-2}\alpha)$, $e\in \{ 3,4,6 \}$ et $e\mid p+1$, $\alpha\in {\BBp}^1({\BBq}_p)$.}
\end{rem}

\begin{pf}
Pour $E/{\BBq}_p$ fix\'ee le syst\`eme ${\bf W}^*_{\ell}(V_{\ell}(E))_{\ell \in {\cal P}}$ \'etant compatible, la classe d'isomorphisme du $(\varphi,N,G)$-module ${\bf D}^*_{pst}(V_p(E))^{(0)}$ obtenu en oubliant la filtration est d\'etermin\'ee par celle de ${\bf W}^*_{\ell}(V_{\ell}(E))$ pour un $\ell \neq p$. Il s'agit donc essentiellement de d\'eterminer quelles sont les filtrations possibles sur ces $(\varphi,N,G)$-modules, sachant que l'on veut obtenir des objets faiblement admissibles de type Hodge-Tate $(0,1)$. On laisse au lecteur le soin de se convaincre que les $(\varphi,N,G)$-modules d\'eduits de la liste ${\bf D^*}$ correspondent exactement aux objets de la liste ${\bf WD^*}$ (via le foncteur ${\bf WD}_{K_e/{\BBq}_p}$ d\'ecrit en~\ref{sec:reppst}) et que les filtrations sont bien celles que l'on veut. 

Signalons que l'on peut aussi faire des calculs directs, au cas par cas : ce sont les m\^emes que ceux effectu\'es dans \cite{Fo-Ma}, $\S$ A, mais avec des conditions. L'avantage de cette approche est qu'elle permet de retrouver le r\'esultat de compatibilit\'e (i.e. ${\bf W}^*_p(V_p(E))$ et ${\bf W}^*_{\ell}(V_{\ell}(E))$, $\ell \neq p$, sont compatibles). 
\end{pf}

\begin{rem}
\label{multrmq}
{\em Soit $E_q/{\BBq}_p$ une courbe de Tate. En utilisant la suite exacte $(*_m)$ (cf.~\ref{sec:notsgeom}) et en faisant des calculs explicites dans $B_{st}$ et $B_{dR}$, on trouve que ${\bf D}^*_{st}(V_p(E_q)) = {\bf D^*_m(1;1;}\alpha(q))$ avec  
$$\alpha(q)\ = \ - \, \frac{ \log(u_q)}{v_p(q)} \makebox[2cm]{o\`u} q = u_q p^{v_p(q)}\: ,\ v_p(q)\geq 1$$
($\log$ est le logarithme $p$-adique). On retrouve ainsi l'invariant not\'e ${\cal L}$ dans \cite{Ma}, $\S$ 3 (voir aussi \cite{LS}, $\S$ 9 et \cite{M-T-T}) ; la diff\'erence dans le signe provient du choix fait pour $B_{st}$, cf.~\ref{sec:reppst}).}
\end{rem}

\begin{rem}
\label{multisogrmq}
{\em Soient $E_q$ et $E_{q'}$ deux courbes de Tate sur ${\BBq}_p$. Avec les notations de la remarque pr\'ec\'edente, les ${\BBq}_p[G]$-modules $V_p(E_{q})$ et $V_p(E_{q'})$ sont isomorphes si et seulement si $\alpha(q)= \alpha(q')$, i.e. $\log(u_q^{v_p(q')})= \log(u_{q'}^{v_p(q)})$, ce qui \'equivaut \`a $q^{v_p(q')(p-1)}=(q')^{v_p(q)(p-1)}$. On retrouve ainsi le fait que $V_p(E_{q})$ et $V_p(E_{q'})$ sont isomorphes si et seulement si les courbes $E_q$ et $E_{q'}$ sont ${\BBq}_p$-isog\`enes (cf. \cite{LS}, $\S$ 9 et \cite{Se 1}, A.1.4).}
\end{rem}

\begin{rem}
{\em Soit $E/{\BBq}_p$ potentiellement ordinaire, i.e. telle que ${\bf D}_{pst}^*(V_p(E)) \simeq {\bf D_{c}^*(e;a_p;\alpha)}$ ou ${\bf D^*_{pc}(e;a_p;\epsilon;\alpha)}$. La suite exacte $(*_{ord})$ (cf.~\ref{sec:notsgeom}) est scind\'ee si $\alpha=0$ et non scind\'ee si $\alpha=1$. Avec les notations utilis\'ees en \ref{sec:imgal}, cela correspond au fait que le ${\BBq}_p$-espace vectoriel $\mbox{Ext}^1({\BBq}_p(\eta_u \xi_e^{\epsilon}),{\BBq}_p(\eta_u^{-1}\xi_e^{-\epsilon}\chi))\simeq H^1(G,{\BBq}_p(\eta_u^{-2}\xi_e^{-2\epsilon}\chi))$ est de dimension un. De plus, on a $\alpha =0$ si et seulement si $j_E = j(e)$ avec $j(3)=j(6)=0$ et $j(4)=1728$, i.e. $E_{L_e}$ est le rel\`evement canonique de $\widetilde{E}_{L_e}$ (cf.~\ref{sec:cpotord}).}
\end{rem}

\begin{rem}
{\em Soit $E/{\BBq}_p$ potentiellement supersinguli\`ere avec $\mbox{dst}(E)=e \geq 3$, i.e. telle que ${\bf D}_{pst}^*(V_p(E)) \simeq {\bf D^*_{pc}(e;0;\alpha)}$. L'invariant $\alpha$ est li\'e au logarithme du groupe formel (de hauteur 2) de $E_{K_e}$. D'apr\`es les r\'esultats de A. Kraus (\cite{Kr}, 2.3.2, prop.2 et lemme 2), les cas $v_p(\alpha)\neq 1$ et $v_p(\alpha)=1$ correspondent respectivement \`a $v_p(\tau_p)\geq \frac{p}{p+1}$ et $v_p(\tau_p)< \frac{p}{p+1}$ o\`u $\tau_p$ est le $p$-i\`eme terme de la s\'erie formelle donnant la multiplication par $p$ dans ce groupe formel ; si l'on choisit une \'equation minimale pour $E$, le cas $v_p(\alpha)\geq 2$ correspond \`a $v_p(\Delta_E)>6$ (i.e. $v_p(\Delta_E)\in \{8,9,10\}$) et le cas $v_p(\alpha)\leq 0$ correspond \`a $v_p(\Delta_E)<6$ (i.e. $v_p(\Delta_E)\in \{2,3,4\}$). De plus, on a $\alpha \in \{ 0, \infty \}$ si et seulement si $j_E=j(e)$ avec $j(3)=j(6)=0$ et $j(4)=1728$ (cf.~\ref{sec:cpotsup}, prop.~\ref{cpotsupprop}).}  
\end{rem}

\begin{rem}
{\em \'Etant donn\'e un objet $\Delta$ de la liste ${\bf WD^*}$, pour obtenir un objet de la liste ${\bf D^*}$ dont la repr\'esentation de Weil-Deligne associ\'ee est isomorphe sur ${\BBc}$ \`a $\Delta$, on a : \\
- une infinit\'e de possibilit\'es param\'etr\'ees par ${\BBq}_p$ dans les cas ${\bf WD^*_m(e;b)}$ \\
- deux possibilit\'es dans les cas ${\bf WD^*_c(e;a_p)}$, $a_p\neq0$ et ${\bf WD^*_{pc}(e;a_p;\epsilon)}$  \\
- une seule possibilit\'e dans les cas ${\bf WD^*_c(e;0)}$  \\
- une infinit\'e de possibilit\'es param\'etr\'ees par ${\BBp}^1({\BBq}_p)$ dans les cas ${\bf WD^*_{pc}(e;0)}$.} 
\end{rem}

\section{Les ${\BBq}_{\ell}[G]$-modules provenant d'une courbe elliptique sur ${\BBq}_p$, $\ell \neq p$} 
\label{sec:lmodellip}

\subsection{R\'esultat et cons\'equence}
\label{sec:resultl}

Soit $\ell \neq p$. Nous allons maintenant donner des conditions n\'ecessaires et suffisantes pour qu'une repr\'esentation $\ell$-adique $V_{\ell}$ de $G$ de dimension $2$ provienne d'une courbe elliptique sur ${\BBq}_p$, i.e. pour qu'il existe  $E/{\BBq}_p$ telle que $V_{\ell} \simeq V_{\ell}(E)$ en tant que ${\BBq}_{\ell}[G]$-modules.  

Pour $E/{\BBq}_p$ elliptique, on sait que ${\bf W}^*_{\ell}(V_{\ell}(E))$ est d\'efinie sur ${\BBq}$ et que $\wedge^2 V_{\ell}(E) = {\BBq}_{\ell}(1)$, i.e. le d\'eterminant sur $V_{\ell}(E)$ est le caract\`ere cyclotomique $\ell$-adique. De plus, si $E$ acquiert bonne r\'eduction sur une extension finie totalement ramifi\'ee $L$, on sait que $a_p(\widetilde{E}_L) \in {\BBz}$ avec $\mid\!a_p(\widetilde{E}_L)\!\mid_{\infty} \leq 2 \sqrt{p}\,$, i.e. les racines du $P_{car}(\mbox{Frob.arithm.}\widetilde{E}_L)$ sont des $p$-nombres de Weil. Cela nous donne trois conditions n\'ecessaires pour qu'une repr\'esentation $\ell$-adique de $G$ provienne d'une courbe elliptique sur ${\BBq}_p$. \'Ecrivons-les en termes de conditions sur les repr\'esentations de Weil-Deligne associ\'ees (i.e. les objets obtenus en appliquant le foncteur ${\bf W}^*_{\ell}$). Soit $(\Delta,\rho_0,N)$ un objet de ${\bf Rep}_{{\BBq}_{\ell}}('W)$ et soit $\phi$ un rel\`evement du Frobenius g\'eom\'etrique dans $W$. On consid\`ere les conditions suivantes : \\
 $(1^{\circ})\ $ $\wedge^2 \Delta$ est donn\'ee par $\wedge^2 \rho_0(I)=1$, $\wedge^2 \rho_0(\phi)=p$ et $\wedge^2 N=0$  \\
 $(2^{\circ})\ $ $\Delta$ est d\'efinie sur ${\BBq}$  \\
 $(3^{\circ})\ $ Si $N=0$ on a $\mbox{Tr}( \rho_0(\phi))\in {\BBz}$ et $\mid\!\mbox{Tr}( \rho_0(\phi))\!\mid_{\infty} \leq 2 \sqrt{p}\,$.

Maintenant le r\'esultat est que ces conditions n\'ecessaires sont \'egalement suffisantes : d'une part les repr\'esentations de Weil-Deligne v\'erifiant les conditions $(1^{\circ})$, $(2^{\circ})$ et $(3^{\circ})$ sont exactement celles de la liste ${\bf WD^*}$, et d'autre part chaque objet de cette liste provient effectivement d'une courbe elliptique sur ${\BBq}_p$.

\begin{thm}
\label{lmodellipthm}
Soient $\ell \neq p$ et $V_{\ell}$ une repr\'esentation $\ell$-adique de $G$ de dimension $2$. Les assertions suivantes sont \'equivalentes : 
\begin{itemize}
\item[(1)] il existe une courbe elliptique $E$ sur ${\BBq}_p$ telle que $V_{\ell}(E)$ soit isomorphe \`a $V_{\ell}$, 
\item[(2)] ${\bf W}^*_{\ell}(V_{\ell})$ v\'erifie les conditions $(1^{\circ})$, $(2^{\circ})$ et $(3^{\circ})$, 
\item[(3)] ${\bf W}^*_{\ell}(V_{\ell})$ est isomorphe \`a un objet de la liste ${\bf WD^*}$.
\end{itemize}
\end{thm}

\begin{pf}
Soient $\phi$ un rel\`evement dans $W$ du Frobenius g\'eom\'etrique et $(\Delta,\rho_0,N)$ un objet de ${\bf Rep}_{{\BBq}_{\ell}}('W)$ v\'erifiant les conditions $(1^{\circ})$, $(2^{\circ})$ et $(3^{\circ})$.

Si $N\neq 0$ les relations $N^2=0$ et $N\rho_0(\phi)=p\rho_0(\phi)N$ montrent que $\rho_0(\phi)$ est diagonalisable avec deux valeurs propres distinctes $(b,pb)$. La condition $(1^{\circ})$ donne $b \in \{ \pm 1\}$ ; en particulier, $\Delta$ est un objet de ${\bf Rep}^{\circ}_{{\BBq}_{\ell}}('W)$. Le sous-groupe d'inertie $I$ agissant \`a travers un quotient fini, donc par des racines de l'unit\'e, la relation $N\rho_0(g)=\rho_0(g)N$ pour $g \in I$ implique $\rho_0(g)=\pm 1$. Donc $\Delta \simeq {\bf WD^*_m(e;b)}$. Ces objets proviennent bien d'une courbe elliptique sur ${\BBq}_p$ : il suffit de prendre n'importe quelle courbe de Tate sur ${\BBq}_p$ et si n\'ecessaire de la tordre.

Si $N= 0$ les conditions $(1^{\circ})$ et $(3^{\circ})$ impliquent que $\Delta$ est un objet de ${\bf Rep}^{\circ}_{{\BBq}_{\ell}}('W)$. Soit $F$ le sous-corps de $\overline{{\BBq}}_p$ fixe par le noyau de la restriction de $\rho_0$ \`a $I$ : c'est une extension finie galoisienne de ${\BBq}_p^{nr}$ telle que ${\rho_0}_{\mid I}$ induit une injection $\rho_0 : I(F/{\BBq}_p^{nr}) \hookrightarrow \mbox{Aut}_{{\BBq}_{\ell}}(\Delta)$. Soit $\tau \in I(F/{\BBq}_p^{nr})$ ; la condition $(1^{\circ})$ impose $\mbox{d\'et}(\rho_0(\tau))=1$ et la condition $(2^{\circ})$ impose $\mbox{Tr}(\rho_0(\tau)) \in {\BBq}\,$. Comme $\mbox{dim}_{{\BBq}_{\ell}}(\Delta)=2$ et que $\rho_0(\tau)$ est d'ordre fini, le polyn\^ome minimal de $\rho_0(\tau)$ est le $e$-i\`eme polyn\^ome cyclotomique avec $e$ tel que $\mbox{{\boldmath $\varphi$}}(e)\in \{ 1,2\}$, o\`u $\mbox{{\boldmath $\varphi$}}$ est la fonction arithm\'etique d'Euler, c'est-\`a-dire $e \in \{ 1,2,3,4,6 \}$. Donc $F/{\BBq}_p^{nr}$ est mod\'er\'ee, cyclique d'ordre $e$, et $F={\BBq}_p^{nr}(\pi_e) = {\BBq}_p^{nr} K_e$, $I(F/{\BBq}_p^{nr})=I(K_e/{\BBq}_p)$.

Si $e\in \{1,2\}$, la condition $(3^{\circ})$ implique que $\Delta$ est isomorphe \`a l'un des ${\bf WD^*_c(e;a_p)}$. Ces objets proviennent tous de courbes elliptiques sur ${\BBq}_p$, puisque, par \cite{Ho-Ta}, pour tout $a_p \in {\cal N}_p$ il existe $\widetilde{E}/{\BBf}_p$ telle que $a_p(\widetilde{E})=a_p$ (laquelle se rel\`eve en un sch\'ema elliptique sur ${\BBz}_p$, que l'on tord sur $M_2$ si $e=2$).

Si $e\in \{ 3,4,6 \}$ et $e\mid p+1$, alors la trace de $\rho_0(\phi)$ doit \^etre nulle et $\Delta$ est isomorphe \`a ${\bf WD^*_{pc}(e;0)}$, voir~\ref{sec:classl}. Si $e\in \{3,4,6\}$ et $e\mid p-1$, alors la condition $(2^{\circ})$ implique que la trace de $\rho_0(\phi)$ est dans ${\cal N}^{\times}_{p,e}$ et $\Delta$ est isomorphe \`a l'un des ${\bf WD^*_{pc}(e;a_p;\epsilon)}$, voir~\ref{sec:classl}. Maintenant tous ces cas sont couverts par les exemples donn\'es en~\ref{sec:exsupl} et~\ref{sec:exordl} : une fois que l'on dispose du lemme~\ref{isogisomlem} qui suit, il s'agit d'un exercice facile sur les courbes elliptiques qui est laiss\'e au lecteur.  
\end{pf}

\medskip
Pour $u\in {\BBf}_p^{\times}$ et $e\in \{ 4,6\}$ on consid\`ere les courbes $\widetilde{E}_{e,u}/{\BBf}_p$ donn\'ees par les \'equations $y^2 = f_{e,u}(x)$ avec $f_{4,u}(x)=x^3 +ux$ et $f_{6,u}(x)=x^3 +u$. On a $j(\widetilde{E}_{e,u})=\tilde{\mbox{\j}}(e)$ avec $\tilde{\mbox{\j}}(4)=1728$ et $\tilde{\mbox{\j}}(6)=0$, donc $\widetilde{E}_{e,u}$ est ordinaire si et seulement si $e\mid p-1$. Toute courbe sur ${\BBf}_p$ d'invariant modulaire $\tilde{\mbox{\j}}(e)$ est ${\BBf}_p$-isomorphe \`a l'une des $\widetilde{E}_{e,u}$ ; de plus, $\widetilde{E}_{e,u}$ et $\widetilde{E}_{e,u'}$ sont ${\BBf}_p$-isomorphes si et seulement si $u \equiv u' \bmod ({\BBf}_p^{\times})^e$ (\cite{Si 1},X, prop.5.4). On en d\'eduit deux applications
$$\left\{ \begin{array}{rcl}
{\BBf}_p^{\times} / ({\BBf}_p^{\times})^e  &  \longrightarrow &  {\cal N}_p \\
u \bmod ({\BBf}_p^{\times})^e              & \longmapsto      &  a_p(\widetilde{E}_{e,u})
\end{array} \right.$$
de l'ensemble des classes de ${\BBf}_p$-isomorphisme de courbes elliptiques d'invariant modulaire $\tilde{\mbox{\j}}(e)$ dans l'ensemble des classes de ${\BBf}_p$-isog\'enie de courbes elliptiques sur ${\BBf}_p$.

\begin{lem}
\label{isogisomlem}
Si $e \mid p-1$ l'association $u \mapsto a_p(\widetilde{E}_{e,u})$ induit une bijection ${\BBf}_p^{\times} / ({\BBf}_p^{\times})^e \stackrel{\sim}{\rightarrow} {\cal N}_{p,e}^{\times}$. En particulier, la classe de ${\BBf}_p$-isog\'enie d'une courbe elliptique d'invariant modulaire $1728$ ou $0$ est aussi sa classe de ${\BBf}_p$-isomorphisme. 
\end{lem}

\begin{pf}
Si $e \mid p-1$, les $\widetilde{E}_{e,u}$ sont ordinaires, d'o\`u $a_p(\widetilde{E}_{e,u})\in {\cal N}_p^{\times}$, et ${\BBf}_p^{\times}/({\BBf}_p^{\times})^e \simeq \mu_e({\BBf}_p)$ est d'ordre $e$. Le fait que la repr\'esentation de Weil-Deligne associ\'ee \`a une tordue convenable d'un rel\`evement d'invariant modulaire $0$ ou $1728$ de $\widetilde{E}_{e,u}$ sur ${\BBq}_p$ est d\'efinie sur ${\BBq}$ implique $a_p(\widetilde{E}_{e,u})\in {\cal N}_{p,e}^{\times}$. Puis on a $a_p(\widetilde{E}_{e,u}) \equiv  1- \#(\widetilde{E}_{e,u}({\BBf}_p)) \bmod p{\BBz} = c_e\, u^{\frac{p-1}{e}} =$ coefficient de $x^{p-1}$ dans $f_{e,u}(x)^{\frac{p-1}{2}}$ (\cite{Si 1},V, preuve du thm.4.1(a)). Le coefficient bin\^omial $c_e$ est non nul et $u^{\frac{p-1}{e}}$ parcourt $\mu_e({\BBf}_p)$ lorsque $u$ parcourt ${\BBf}_p^{\times}$. Comme $\mbox{Card}({\cal N}_{p,e}^{\times})=e$, on en d\'eduit le r\'esultat.
\end{pf}

\begin{cor}
\label{nbclasslcor}
Soit $\ell \in {\cal P}$ tel que $\ell \neq p$. Le nombre de classes d'isomorphisme d'objets de ${\bf Rep}_{{\BBq}_{\ell}}(G)$ provenant d'une courbe elliptique sur ${\BBq}_p$ est fini et ind\'ependant de $\ell$ ; il vaut $4[2 \sqrt p\,] + \lambda (p)$ o\`u $\lambda(p)=38, 16, 31$ ou $9$ suivant que $p\equiv 1, 5, 7$ ou $11 \bmod 12$.
\end{cor}

Plus pr\'ecis\'ement, il y a : 4 classes dans ${\bf Rep}_{{\BBq}_{\ell}}(G)$ provenant de courbes elliptiques sur ${\BBq}_p$ n'ayant pas potentiellement bonne r\'eduction ; $\mbox{Card}({\cal N}_p^{\times})=2 [2\sqrt{p}\,]$ classes provenant de courbes ayant bonne r\'eduction ordinaire sur ${\BBq}_p$ et autant provenant d'un twist quadratique ramifi\'e de telles ; 1 classe provenant de courbes ayant bonne r\'eduction supersinguli\`ere sur ${\BBq}_p$ et 1 provenant d'un twist quadratique ramifi\'e de telles. Si $3 \mid p-1$ il y a $2\,\mbox{Card}({\cal N}_{p,3}^{\times})=12$ classes provenant de $E/{\BBq}_p$ potentiellement ordinaires avec $\mbox{dst}(E)=3$ et $2\,\mbox{Card}({\cal N}_{p,6}^{\times})=12$ classes provenant de telles courbes avec $\mbox{dst}(E)=6$ ; si $3 \mid p+1$ il y a 1 classe provenant de $E/{\BBq}_p$ potentiellement supersinguli\`eres avec $\mbox{dst}(E)=3$ et 1 classe provenant de telles courbes avec $\mbox{dst}(E)=6$. Si $4 \mid p-1$ il y a $2\,\mbox{Card}({\cal N}_{p,4}^{\times})=8$ classes provenant de $E/{\BBq}_p$ potentiellement ordinaires avec $\mbox{dst}(E)=4$ ; si $4 \mid p+1$ il y a 1 classe provenant de $E/{\BBq}_p$ potentiellement supersinguli\`eres avec $\mbox{dst}(E)=4$.

\begin{rem}
{\em Soit $E/{\BBq}_p$ une courbe elliptique et soit $T_{\ell}$ un r\'eseau $G$-stable de $V_{\ell}(E)$ ; \`a homoth\'etie pr\`es, on peut supposer que $T_{\ell}\subset T_{\ell}(E)$. Alors il existe une courbe elliptique $E'$ et une $\ell$-isog\'enie $\psi : E' \rightarrow E$, d\'efinies sur ${\BBq}_p$, telles que $\psi_{\ell}(T_{\ell}(E'))=T_{\ell}$. Donc, si $T_{\ell}$ est un ${\BBz}_{\ell}[G]$-module, pour qu'il existe une courbe elliptique $E'/{\BBq}_p$ telle que $T_{\ell} \simeq T_{\ell}(E')$, il faut et il suffit qu'il existe une courbe elliptique $E/{\BBq}_p$ telle que ${\BBq}_{\ell}\!\otimes_{{\BBz}_{\ell}}\!T_{\ell} \simeq V_{\ell}(E)$.} 
\end{rem}

\subsection{Exemples}
\label{sec:exl}

\subsubsection{Courbes elliptiques potentiellement ordinaires}
\label{sec:exordl}

\underline{Si $4 \mid p-1$} : pour chaque $a_{p,j} \in {\cal N}_{p,4}^{\times}$, $1\leq j \leq 4$, on choisit un $u_j \in {\BBf}_p^{\times}$ tel que $a_p(\widetilde{E}_j)=a_{p,j}$ avec $\widetilde{E}_j : y^2= x^3 +u_jx$ ; les $u_j$ sont un syst\`eme de repr\'esentants de ${\BBf}_p^{\times}/({\BBf}_p^{\times})^4$ (lemme~\ref{isogisomlem}). Ces courbes sont ordinaires d'invariant modulaire $1728$. Par exemple, si $p=5$, on peut prendre $u_1=1$, $u_2=2$, $u_3=-2$, $u_4=-1$, ce qui donne $a_{5,1}=2$, $a_{5,2}=4$, $a_{5,3}=-4$, $a_{5,4}=-2$. Alors $\{ \, [u_j](-p)^i,\, 1 \leq j \leq 4,\, 0\leq i \leq 3 \, \}$ est un syst\`eme de repr\'esentants de ${\BBq}_p^{\times}/({\BBq}_p^{\times})^4$, o\`u $[u_j] \in {\BBz}_p^{\times}$ est le repr\'esentant de Teichm\"{u}ller de $u_j$.

On pose $E_{i,j} : y^2 = x^3 + [u_j](-p)^i x$ pour $1 \leq j \leq 4$ et $0\leq i \leq 3$. Ce sont des courbes sur ${\BBq}_p$ d'invariant modulaire $1728$ repr\'esentant les \'el\'ements de $\mbox{Twist}((E_{0,1},{\bf 0}),{\BBq}_p) \simeq {\BBq}_p^{\times}/({\BBq}_p^{\times})^4$. On a alors pour chaque $j$ : 
$$\begin{array}{rclcrcl}
{\bf W}^*_{\ell}(V_{\ell}(E_{0,j}))  & \simeq  &  {\bf WD^*_c(1;a_{p,\,j})} & \hspace{0.5cm} & {\bf W}^*_{\ell}(V_{\ell}(E_{1,j})) & \simeq  &  {\bf WD^*_{pc}(4;a_{p,\,j};1)}   \\
{\bf W}^*_{\ell}(V_{\ell}(E_{2,j})) & \simeq  &  {\bf WD^*_c(2;a_{p,\,j})} & \hspace{0.5cm} & {\bf W}^*_{\ell}(V_{\ell}(E_{3,j})) & \simeq  &  {\bf WD^*_{pc}(4;a_{p,\,j};-1)}
\end{array}$$

\medskip
\noindent \underline{Si $3 \mid p-1$} : pour chaque $a_{p,j} \in {\cal N}_{p,3}^{\times}$, $1\leq j \leq 6$, on choisit un $v_j \in {\BBf}_p^{\times}$ tel que $a_p(\widetilde{\cal E}_j)=a_{p,j}$ avec $\widetilde{\cal E}_j : y^2= x^3 +v_j$ ; les $v_j$ sont un syst\`eme de repr\'esentants de ${\BBf}_p^{\times}/({\BBf}_p^{\times})^6$ (lemme~\ref{isogisomlem}). Ces courbes sont ordinaires d'invariant modulaire $0$. Par exemple, si $p=7$, on peut prendre $v_1=1$, $v_2=2$, $v_3=3$, $v_4=-3$, $v_5=-2$, $v_6=-1$, ce qui donne $a_{7,1}=-4$, $a_{7,2}=-1$, $a_{7,3}=-5$, $a_{7,4}=5$, $a_{7,5}=1$, $a_{7,6}=4$. Alors $\{\, [v_j](-p)^i,\, 1 \leq j \leq 6,\, 0\leq i \leq 5\, \}$ est un syst\`eme de repr\'esentants de ${\BBq}_p^{\times}/({\BBq}_p^{\times})^6$.

On pose ${\cal E}_{i,j} : y^2 = x^3 + [v_j](-p)^i$ pour $1 \leq j \leq 6$ et $0\leq i \leq 5$. Ce sont des courbes sur ${\BBq}_p$ d'invariant modulaire $0$ repr\'esentant les \'el\'ements de $\mbox{Twist}(({\cal E}_{0,1},{\bf 0}),{\BBq}_p) \simeq {\BBq}_p^{\times}/({\BBq}_p^{\times})^6$. On a alors pour chaque $j$ :
$$\begin{array}{rclcrcl}
{\bf W}^*_{\ell}(V_{\ell}({\cal E}_{0,j}))  & \simeq  &  {\bf WD^*_c(1;a_{p,\,j})} & \hspace{0.5cm} & {\bf W}^*_{\ell}(V_{\ell}({\cal E}_{1,j}))  & \simeq  &  {\bf WD^*_{pc}(6;a_{p,\,j};1)}  \\
{\bf W}^*_{\ell}(V_{\ell}({\cal E}_{2,j}))  & \simeq  &  {\bf WD^*_{pc}(3;a_{p,\,j};1)} & \hspace{0.5cm} & {\bf W}^*_{\ell}(V_{\ell}({\cal E}_{3,j}))  & \simeq  &  {\bf WD^*_c(2;a_{p,\,j})}      \\
{\bf W}^*_{\ell}(V_{\ell}({\cal E}_{4,j}))  & \simeq  &  {\bf WD^*_{pc}(3;a_{p,\,j};-1)} & \hspace{0.5cm} & {\bf W}^*_{\ell}(V_{\ell}({\cal E}_{5,j}))  & \simeq  &  {\bf WD^*_{pc}(6;a_{p,\,j};-1)}
\end{array}$$

\subsubsection{Courbes elliptiques potentiellement supersinguli\`eres}
\label{sec:exsupl}

\underline{Si $4 \mid p+1$} : on pose $E_i : y^2 = x^3 + (-p)^i x$ pour $0\leq i \leq 3$. Ce sont des courbes sur ${\BBq}_p$ d'invariant modulaire $1728$ repr\'esentant les \'el\'ements de $\mbox{Twist}((E_0,{\bf 0}),{\BBq}_p) \simeq {\BBq}_p^{\times}/({\BBq}_p^{\times})^4$, la courbe r\'eduite de $E_0$ \'etant d'\'equation $y^2 = x^3 + x$. On a alors :
$$\begin{array}{rclcrcl}
{\bf W}^*_{\ell}(V_{\ell}(E_0)) & \simeq & {\bf WD^*_c(1;0)} & \hspace{0.5cm} & {\bf W}^*_{\ell}(V_{\ell}(E_1)) & \simeq & {\bf WD^*_{pc}(4;0)}  \\
{\bf W}^*_{\ell}(V_{\ell}(E_2)) & \simeq & {\bf WD^*_c(2;0)} & \hspace{0.5cm} & {\bf W}^*_{\ell}(V_{\ell}(E_3)) & \simeq & {\bf WD^*_{pc}(4;0)}
\end{array}$$

\medskip
\noindent \underline{Si $3 \mid p+1$} : on pose ${\cal E}_i : y^2 = x^3 + (-p)^i$ pour $0\leq i \leq 5$. Ce sont des courbes sur ${\BBq}_p$ d'invariant modulaire $0$ repr\'esentant les \'el\'ements de $\mbox{Twist}(({\cal E}_0,{\bf 0}),{\BBq}_p) \simeq {\BBq}_p^{\times}/({\BBq}_p^{\times})^6$, la courbe r\'eduite de ${\cal E}_0$  \'etant d'\'equation $y^2 = x^3+1$. On a alors :
$$\begin{array}{rclcrcl}
{\bf W}^*_{\ell}(V_{\ell}({\cal E}_0)) & \simeq & {\bf WD^*_c(1;0)} & \hspace{0.5cm} & {\bf W}^*_{\ell}(V_{\ell}({\cal E}_1)) & \simeq & {\bf WD^*_{pc}(6;0)}     \\
{\bf W}^*_{\ell}(V_{\ell}({\cal E}_2)) & \simeq & {\bf WD^*_{pc}(3;0)} & \hspace{0.5cm} & {\bf W}^*_{\ell}(V_{\ell}({\cal E}_3)) & \simeq & {\bf WD^*_c(2;0)}      \\
{\bf W}^*_{\ell}(V_{\ell}({\cal E}_4)) & \simeq & {\bf WD^*_{pc}(3;0)} & \hspace{0.5cm} & {\bf W}^*_{\ell}(V_{\ell}({\cal E}_5)) & \simeq & {\bf WD^*_{pc}(6;0)}
\end{array}$$

\section{Construction de courbes potentiellement supersinguli\`eres}
\label{sec:construct}

Pour d\'eterminer tous les ${\BBq}_p[G]$-modules provenant d'une courbe elliptique sur ${\BBq}_p$ nous allons suivre la m\^eme strat\'egie que pour les cas $\ell \neq p$. Mais cette fois les op\'erations \'el\'ementaires sur les courbes elliptiques ne suffisent plus : il faut construire toutes les courbes sur ${\BBq}_p$ potentiellement supersinguli\`eres de d\'efaut de semi-stabilit\'e sup\'erieur ou \'egal \`a $3$.

Soit $O_{L_e}={\BBz}_p[\pi_e]$ l'anneau des entiers de $L_e={\BBq}_p(\pi_e)$. Dans un premier temps, pour $e<p-1$, on construit \`a partir d'une courbe elliptique $\widetilde{E}/{\BBf}_p$ supersinguli\`ere fix\'ee tous les sch\'emas elliptiques ${\cal E}/O_{L_e}$ relevant $\widetilde{E}$, \`a $O_{L_e}$-isomorphisme pr\`es (la restriction sur l'indice de ramification provient de la th\'eorie des modules de Dieudonn\'e filtr\'es). Puis, pour $e\in \{ 3,4,6\}$ et $e<p-1$, on d\'etermine parmi ces sch\'emas ceux qui sont susceptibles d'\^etre d\'efinis sur ${\BBq}_p$ avec un d\'efaut de semi-stabilit\'e $e$, c'est-\`a-dire ceux pour lesquels il existe une courbe elliptique $E/{\BBq}_p$ telle que $E\!\times_{\!{\BBq}_p}\!L_e \simeq {\cal E}\!\times_{\!O_{L_e}}\!L_e$ et $\mbox{dst}(E)=e$. On obtient finalement toutes les courbes elliptiques sur ${\BBq}_p$ potentiellement supersinguli\`eres, \`a ${\BBq}_p$-isomorphisme pr\`es.

\subsection{Pr\'eliminaires}
\label{sec:prelim}

\subsubsection{Le foncteur de Serre-Tate}
\label{sec:ST}

Soit ${\bf {\cal SE}}_{O_{L_e}}$ la cat\'egorie des sch\'emas elliptiques sur $O_{L_e}$. On note ${\bf {\cal C}}_{O_{L_e}}$ la cat\'egorie dont les objets sont les triplets $(\widetilde{B},\Gamma,\nu)$ o\`u $\widetilde{B}/{\BBf}_p$ est une courbe elliptique, $\Gamma$ un groupe $p$-divisible sur $O_{L_e}$ et $\nu : \widetilde{B}(p) \stackrel{\sim}{\rightarrow} \widetilde{\Gamma}= \Gamma\!\times_{O_{L_e}}\!{\BBf}_p$ un isomorphisme de groupes $p$-divisibles sur ${\BBf}_p$ ; les morphismes $(\widetilde{B},\Gamma,\nu) \rightarrow (\widetilde{B}',\Gamma',\nu')$ sont les couples $(\gamma,\psi)$ o\`u $\gamma : \widetilde{B} \rightarrow \widetilde{B}'$ est un morphisme de courbes elliptiques sur ${\BBf}_p$ et $\psi : \Gamma \rightarrow \Gamma'$ est un morphisme de groupes $p$-divisibles sur $O_{L_e}$ tels que $\nu' \circ \gamma(p) = \widetilde{\psi} \circ \nu $. Le th\'eor\`eme de Serre-Tate implique que le foncteur ${\bf ST}$ de ${\bf {\cal SE}}_{O_{L_e}}$ dans ${\bf {\cal C}}_{O_{L_e}}$ d\'efini par ${\bf ST}({\cal A}) = (\widetilde{\cal A}, {\cal A}(p),\nu_{can})$, o\`u $\widetilde{\cal A}= {\cal A}\!\times_{O_{L_e}}\!{\BBf}_p$, \'etablit une \'equivalence de cat\'egories (voir \cite{Ka} plus le fait que tout sch\'ema alg\'ebro\"{\i}de de dimension relative 1 est alg\'ebrisable).

Pour $\widetilde{E}/{\BBf}_p$ fix\'ee, on note ${\bf {\cal SE}}_{O_{L_e}}(\widetilde{E})$ la sous-cat\'egorie pleine de ${\bf {\cal SE}}_{O_{L_e}}$ form\'ee des objets ${\cal A}$ tels qu'il existe un ${\BBf}_p$-isomorphisme $\widetilde{\cal A} \simeq \widetilde{E}$ (i.e. ${\cal A}$ est un sch\'ema elliptique sur $O_{L_e}$ relevant $\widetilde{E}$) et ${\bf {\cal C}}_{O_{L_e}}(\widetilde{E})$ la sous-cat\'egorie pleine de ${\bf {\cal C}}_{O_{L_e}}$ form\'ee des triplets du type $(\widetilde{E},\Gamma ,\nu)$. On voit que ${\bf {\cal SE}}_{O_{L_e}}(\widetilde{E})$ est la sous-cat\'egorie pleine de ${\bf {\cal SE}}_{O_{L_e}}$ form\'ee des objets ${\cal A}$ tels qu'il existe un objet $X$ de ${\bf {\cal C}}_{O_{L_e}}(\widetilde{E})$ avec ${\bf ST}({\cal A}) \simeq X$ dans ${\bf {\cal C}}_{O_{L_e}}$ ; on a, avec des notations \'evidentes,
 $$\mbox{Hom}_{{\bf {\cal SE}}_{O_{L_e}}(\widetilde{E})}({\cal A},{\cal A}') \ \simeq \ \mbox{Hom}_{{\bf {\cal C}}_{O_{L_e}}(\widetilde{E})} (X,X') $$ 
Donc le foncteur ${\bf ST}$ induit une bijection entre les classes d'isomorphisme de ${\bf {\cal SE}}_{O_{L_e}}(\widetilde{E})$ et celles de ${\bf {\cal C}}_{O_{L_e}}(\widetilde{E})$. Ainsi, \'etudier les sch\'emas elliptiques sur $O_{L_e}$ relevant $\widetilde{E}$ revient \`a \'etudier les rel\`evements du groupe $p$-divisible $\widetilde{E}(p)$ sur ${\BBf}_p$ en un groupe $p$-divisible sur $O_{L_e}$.

\subsubsection{Modules de Dieudonn\'e}
\label{sec:MD}

Soient $k\subset \overline{{\BBf}}_p$ un corps fini et $\sigma$ le Frobenius absolu agissant sur $k$ (par $x \mapsto x^p$) et sur $W(k)$, l'anneau des vecteurs de Witt \`a coefficients dans $k$. Soit $\widetilde{\Gamma}$ un groupe $p$-divisible sur $k$. On note $M= {\bf M}_k(\widetilde{\Gamma})= \mbox{Hom}_{{\bf D}_k}(\widetilde{\Gamma}, C\widehat{W}(k))$ son module de Dieudonn\'e sur $k$ (voir \cite{Fo 4})~: c'est un $W(k)$-module libre de rang fini muni d'un op\'erateur de Frobenius $\sigma$-semi-lin\'eaire $\varphi$ v\'erifiant $pM \subset \varphi M$. Le foncteur ${\bf M}_k$ induit une anti-\'equivalence entre la cat\'egorie des groupes $p$-divisibles sur $k$ et celle des $W(k)$-modules libres de rang fini munis d'un op\'erateur $\varphi$ comme ci-dessus (\cite{Fo 4},III, prop.6.1 et rmq.3 qui suit). On note ${\bf MD}_{k}$ cette derni\`ere cat\'egorie ; si $k={\BBf}_p$ on \'ecrit ${\bf M}_{{\BBf}_p}={\bf M}$.

Si $\widetilde{E}/{\BBf}_p$ est supersinguli\`ere alors $M= {\bf M}(\widetilde{E}(p))$ est un ${\BBz}_p$-module libre de rang 2 muni d'un Frobenius ${\BBz}_p$-lin\'eaire $\varphi= {\bf M}(\mbox{Frob}_{\widetilde{E}}(p))$ v\'erifiant ${\varphi}^2 +p=0$. Si $x \in M$ est tel que $x \not\in \varphi M$ alors $\varphi x \not\in pM$ et $(x,\varphi x)$ est une base de $M$ ; un tel $x$ est donc un g\'en\'erateur du ${\BBz}_p[\varphi]$-module $M$. En particulier, on en d\'eduit que deux courbes elliptiques sur ${\BBf}_p$ supersinguli\`eres ont des groupes $p$-divisibles isomorphes, et l'isomorphisme canonique ${\BBz}_p\!\otimes_{{\BBz}}\mbox{Hom}_{{\BBf}_p}(\widetilde{E},\widetilde{E}') \simeq \mbox{Hom}_{{\BBf}_p} (\widetilde{E}(p), \widetilde{E}'(p))$ (cf. \cite{Wa-Mi},II) montre qu'elles sont li\'ees par une ${\BBf}_p$-isog\'enie de degr\'e premier \`a $p$.

\medskip
Il y a, \`a ${\BBf}_p$-isomorphisme pr\`es, deux courbes elliptiques sur ${\BBf}_p$ supersinguli\`eres ayant un invariant modulaire donn\'e, l'une est une tordue sur ${\BBf}_{p^2}$ de l'autre ; de plus, $\mbox{Aut}_{{\BBf}_p}(\widetilde{E})$ est toujours d'ordre 2 lorsque $\widetilde{E}$ est supersinguli\`ere (\cite{Si 1}, prop.5.4 et cor.5.4.1).

Soit $\widetilde{E}/{\BBf}_p$ supersinguli\`ere d'invariant modulaire $\tilde{ \mbox{\j}}(e)$ avec $e \in \{ 3,4,6 \}$ et $\tilde{ \mbox{\j}}(3)=\tilde{ \mbox{\j}}(6)=0$, $\tilde{ \mbox{\j}}(4)=1728$~; donc $e \mid p+1$ et $[\zeta_e] \in \mbox{Aut}_{{\BBf}_{\!p^2}}(\widetilde{E})$. Soit $f$ le Frobenius arithm\'etique agissant sur $\widetilde{E}$~; comme $p\equiv -1 \bmod e{\BBz}$, on a $[\zeta_e]f=f[\zeta_e]^{-1}$. Le Frobenius $\varphi$ agissant sur $M={\bf M}(\widetilde{E}(p))$ s'\'etend $\sigma$-semi-lin\'eairement sur $R=M\!\otimes_{{\BBz}_p}\!{\BBz}_{p^2}$ o\`u ${\BBz}_{p^2}=W({\BBf}_{p^2})$ est l'anneau des entiers de ${\BBq}_{p^2}$. Notons $\xi_e = {\bf M}_{{\BBf}_{p^2}}([\zeta_e](p))$ : c'est un automorphisme ${\BBz}_{p^2}$-lin\'eaire de $R$ d'ordre $e$ et de d\'eterminant 1. L'objet $D=M\!\otimes_{{\BBz}_p}\!{\BBq}_{p^2}$ est un $(\varphi,G_{K_e/L_e})$-module de dimension $2$ dans lequel la relation $[\zeta_e]f=f[\zeta_e]^{-1}$ se traduit par les relations $\xi_e\varphi=\varphi\xi_e$ et $\omega \xi_e = \xi_e^{-1} \omega$ (avec $<\omega>=G_{K_e/L_e}$, cf.~\ref{sec:notsgal}). Il existe alors une ${\BBz}_{p^2}$-base $(e_1 ,e_2)$ de $R \subset D$ et un $\eta \in ({\BBz}/e{\BBz})^{\times}=\{ \pm 1 \}$ tels que 
$$\varphi e_1 = e_2 \  ,\  \varphi e_2 = -pe_1 \makebox[1cm]{;}  \omega e_1=e_1 \ ,\  \omega e_2 =e_2 \makebox[1cm]{;} \xi_e e_1 = \zeta_e^{\eta}e_1 \ ,\  \xi_e e_2 = \zeta_e^{-\eta} e_2 $$
En effet, on montre d'abord l'existence d'une telle base pour $D$ ; puis, \`a homoth\'etie pr\`es, tout ${\BBz}_{p^2}$-r\'eseau de $D$ stable par $\varphi$ est de la forme ${\BBz}_{p^2}e_1 \oplus {\BBz}_{p^2}e_2$ ou $\varphi({\BBz}_{p^2}e_1 \oplus {\BBz}_{p^2}e_2)$. Dans cette situation, on dira que le g\'en\'erateur $e_1$ du ${\BBz}_p[\varphi]$-module $M$ est adapt\'e au groupe d'automorphismes de $\widetilde{E}$. Le choix d'un tel g\'en\'erateur est unique \`a un \'el\'ement de ${\BBz}_p^{\times}$ pr\`es.

\subsubsection{Modules de Dieudonn\'e filtr\'es}
\label{sec:MDfil}

Pour \'etudier les rel\`evements sur $O_{L_e}$ de groupes $p$-divisibles sur ${\BBf}_p$ on utilise la th\'eorie des modules de Dieudonn\'e filtr\'es sur $O_{L_e}$ telle qu'elle est d\'ecrite dans \cite{Fo 4},IV,$\S$2 \`a 5. Pour $e \leq p-1$ on d\'efinit la cat\'egorie ${\bf MD}_{O_{L_e}}$ suivante : \\
{\bf -} les objets sont les couples $(M,{\cal L})$, o\`u $M$ est un ${\BBz}_p$-module libre de rang fini muni d'un op\'erateur de Frobenius ${\BBz}_p$-lin\'eaire $\varphi$ tel que $pM \subset \varphi M$ (i.e. $M$ est un objet de ${\bf MD}_{{\BBf}_p}$), et ${\cal L}$ est un sous-$O_{L_e}$-module de 
$$ {\cal M}= M\!\otimes_{{\BBz}_p}\!O_{L_e} + \varphi M\!\otimes_{{\BBz}_p}\!\pi_e^{1-e}O_{L_e} \subset M\!\otimes_{{\BBz}_p}\!L_e $$
tel que l'inclusion ${\cal L} \hookrightarrow {\cal M}$ induit un isomorphisme de ${\BBf}_p$-espaces vectoriels 
$$ {\cal L} / \pi_e {\cal L} \; \simeq \; {\cal M} / (\varphi M \!\otimes_{{\BBz}_p}\!\pi_e^{1-e}O_{L_e}) \hspace{2cm} (*) $$
{\bf -} un morphisme $(M,{\cal L}) \rightarrow (M',{\cal L'})$ est une application ${\BBz}_p$-lin\'eaire $\psi :\, M \rightarrow M'$ qui commute aux Frobenius et qui, apr\`es extension des scalaires, envoie ${\cal L}$ dans ${\cal L'}$.

Soit $\Gamma$ un groupe $p$-divisible sur $O_{L_e}$ et $\widetilde{\Gamma}= \Gamma\!\times_{O_{L_e}}\!{\BBf}_p$ sa fibre sp\'eciale. Dans \cite{Fo 4} J.-M. Fontaine construit un sous-$O_{L_e}$-module ${\cal L}(\Gamma)$ de ${\cal M}(\widetilde{\Gamma}) = {\bf M}(\widetilde{\Gamma})\!\otimes_{{\BBz}_p}\!O_{L_e} + \varphi {\bf M}(\widetilde{\Gamma})\!\otimes_{{\BBz}_p}\!\pi_e^{1-e}O_{L_e}$ qui v\'erifie $(*)$, de sorte que le couple $({\bf M}(\widetilde{\Gamma}),{\cal L}(\Gamma))$ est un objet de ${\bf MD}_{O_{L_e}}$. L'association $\Gamma \mapsto {\bf M}_{O_{L_e}}(\Gamma) = ({\bf M}(\widetilde{\Gamma}),{\cal L}(\Gamma))$ est fonctorielle et lorsque $e<p-1$ elle induit une anti-\'equivalence entre la cat\'egorie des groupes $p$-divisibles sur $O_{L_e}$ et ${\bf MD}_{O_{L_e}}$ (\cite{Fo 4},IV, prop.5.1). De plus, on a un isomorphisme canonique de $\varphi$-modules filtr\'es ${\BBq}_p\otimes_{{\BBz}_p}{\bf M}_{O_{L_e}}(\Gamma) \simeq {\bf D}^*_{cris,L_e}(V_p(\Gamma))$ qui fait le lien entre le Module de Dieudonn\'e filtr\'e de $\Gamma$ et la th\'eorie cristalline (voir \cite{Fo 5}).

Dans toute la suite on suppose $e<p-1$. Soit $M$ un objet de ${\bf MD}_{{\BBf}_p}$ ; les rel\`evements de $M$ en un objet de ${\bf MD}_{O_{L_e}}$ sont alors en correspondance bijective avec les sous-$O_{L_e}$-modules ${\cal L}$ de ${\cal M}$ tels que l'inclusion ${\cal L} \hookrightarrow {\cal M}$ induit un isomorphisme de ${\BBf}_p$-espaces vectoriels ${\cal L} / \pi_e {\cal L}  \simeq  {\cal M} / (\varphi M \!\otimes_{{\BBz}_p}\!\pi_e^{1-e}O_{L_e})$.

\subsection{Sch\'emas elliptiques supersinguliers}
\label{sec:schellsup}

On fixe $\widetilde{E}/{\BBf}_p$ supersinguli\`ere et l'on note $M= {\bf M}(\widetilde{E}(p))$ et ${\cal M} = {\cal M}(\widetilde{E}(p))$. On choisit un g\'en\'erateur $e_1$ du ${\BBz}_p[\varphi]$-module $M$, d'o\`u $M= {\BBz}_p e_1 \oplus {\BBz}_p e_2$ avec $\varphi e_1 = e_2$ et $\varphi e_2 = -pe_1$. Alors on a ${\cal M} = O_{L_e}(e_1\otimes 1)\oplus O_{L_e}(e_2\otimes \pi_e^{1-e})$ et 
$${\cal M}/ (\varphi M\!\otimes_{{\BBz}_p}\!\pi_e^{1-e}O_{L_e})\; \simeq \; {\BBf}_p ((e_1 \otimes 1) \bmod \pi_e {\cal M})$$
On en d\'eduit que les rel\`evements de $M$ en un module de Dieudonn\'e sur $O_{L_e}$ correspondent bijectivement aux
$${\cal L}(\beta) = ( e_1 \otimes 1 + \beta \cdot e_2 \otimes \pi_e^{1-e} )O_{L_e} \makebox[1cm]{,}\beta \in O_{L_e}$$

Soit $(\widetilde{E},\Gamma,\nu)$ un objet de ${\bf {\cal C}}_{O_{L_e}}(\widetilde{E})$. Alors ${\bf M}_{O_{L_e}}(\Gamma) = ({\bf M}(\widetilde{\Gamma}),{\cal L}(\Gamma))$ est un objet de ${\bf MD}_{O_{L_e}}$ et ${\bf M}(\nu) : {\bf M}(\widetilde{\Gamma}) \stackrel{\sim}{\rightarrow} M$ est un isomorphisme dans ${\bf MD}_{{\BBf}_p}$ ; on a donc ${\bf M}(\nu)(\varphi {\bf M}(\widetilde{\Gamma}))= \varphi M$. On note ${\bf M}(\nu)_{L_e}$ l'application $L_e$-lin\'eaire de ${\cal M}(\widetilde{\Gamma})$ dans ${\cal M}$ d\'eduite par extension des scalaires ; on a ${\bf M}(\nu)_{L_e}({\cal M}(\widetilde{\Gamma})) = {\cal M}$. Alors ${\bf M}(\nu)_{L_e}({\cal L}(\Gamma))$ est un sous-$O_{L_e}$-module de rang 1 de ${\cal M}$ tel que l'inclusion induit un isomorphisme de ${\BBf}_p$-espaces vectoriels
 $$ {\bf M}(\nu)_{L_e}({\cal L}(\Gamma))\, / \, \pi_e ({\bf M}(\nu)_{L_e}({\cal L}(\Gamma))) \; \simeq \; {\cal M} / (\varphi M\!\otimes_{{\BBz}_p}\!\pi_e^{1-e}O_{L_e}) $$
Il existe donc un unique $\beta \in O_{L_e}$ tel que ${\bf M}(\nu)_{L_e}({\cal L}(\Gamma)) = {\cal L}(\beta) \subset {\cal M}$.

\begin{lem}
\label{homlem}
Soient $(\widetilde{E},\Gamma,\nu)$, $(\widetilde{E},\Gamma',\nu')$ deux objets de ${\bf {\cal C}}_{O_{L_e}}(\widetilde{E})$ et $\beta ,\beta' \in O_{L_e}$ tels que ${\bf M}(\nu)_{L_e}({\cal L}(\Gamma))={\cal L}(\beta)$ et ${\bf M}(\nu')_{L_e}({\cal L}(\Gamma'))={\cal L}(\beta')$. Alors on a un isomorphisme
 $$\mbox{\em Hom}_{{\bf {\cal C}}_{O_{L_e}}(\widetilde{E})} \Bigl( (\widetilde{E} , \Gamma , \nu ) , (\widetilde{E} , \Gamma' , \nu' ) \Bigr) \   \simeq \  \{\; \gamma \in \mbox{\em End}_{{\BBf}_p}(\widetilde{E})\  / \ {\bf M}(\gamma(p))_{L_e}({\cal L}(\beta')) \subset {\cal L}(\beta)\; \} $$
\end{lem}

\begin{pf}
C'est imm\'ediat, en utilisant les d\'efinitions des cat\'egories ${\bf {\cal C}}_{O_{L_e}}(\widetilde{E})$ et ${\bf MD}_{O_{L_e}}$ ainsi que la pleine fid\'elit\'e des foncteurs ${\bf M}$ et ${\bf M}_{O_{L_e}}$.
\end{pf}

\begin{prop}
\label{schellsupprop}
Soit $e<p-1$. Soit $\widetilde{E}/{\BBf}_p$ une courbe elliptique supersinguli\`ere. Via le choix d'un g\'en\'erateur du ${\BBz}_p[\varphi]$-module ${\bf M}(\widetilde{E}(p))$, l'association 
$$ \left \{
   \begin{array}{rll}
   {\bf {\cal C}}_{O_{L_e}}(\widetilde{E})  &  \rightarrow   &  O_{L_e}   \\
   (\widetilde{E} , \Gamma , \nu )           &  \mapsto       &  \beta \makebox[1.5cm]{tel que} {\cal L}(\beta)= {\bf M}(\nu)_{L_e}({\cal L}(\Gamma))
   \end{array}
   \right . $$
induit une bijection entre les classes d'isomorphisme dans ${\bf {\cal C}}_{O_{L_e}}(\widetilde{E})$ et $O_{L_e}$. 
\end{prop}

En composant avec le foncteur ${\bf ST}$ on obtient une bijection entre les classes d'isomorphisme dans ${\bf {\cal SE}}_{O_{L_e}}(\widetilde{E})$ et $O_{L_e}$. Notons que le choix d'un autre g\'en\'erateur du ${\BBz}_p[\varphi]$-module $M$ change l'invariant $\beta$ en $(a+b\pi_e \beta)^{-1}(b\pi_e^{e-1}+a\beta)$ avec $a \in {\BBz}_p^{\times}$ et $b\in {\BBz}_p$.

\medskip
\begin{pf}
Montrons d'abord la surjectivit\'e. Soit $\beta \in O_{L_e}$. Il existe un groupe $p$-divisible $J_{\beta}$ sur $O_{L_e}$ et un isomorphisme $\xi_{\beta} : {\bf M}_{O_{L_e}}(J_{\beta}) \stackrel{\sim}{\rightarrow} (M,{\cal L}(\beta))$ dans ${\bf MD}_{O_{L_e}}$. Il induit un isomorphisme $\xi_{\beta} : {\bf M}(\widetilde{J}_{\beta}) \stackrel{\sim}{\rightarrow} M$ de modules de Dieudonn\'e sur ${\BBf}_p$ et il existe un unique isomorphisme $\nu_{\beta} : \widetilde{E}(p) \stackrel{\sim}{\rightarrow} \widetilde{J}_{\beta}$ de groupes $p$-divisibles sur ${\BBf}_p$ tel que ${\bf M}(\nu_{\beta})=\xi_{\beta}$. Alors le triplet $(\widetilde{E},J_{\beta}, \nu_{\beta})$ est un objet de ${\cal C}_{O_{L_e}}(\widetilde{E})$ tel que ${\bf M}(\nu_{\beta})_{L_e}({\cal L}(J_{\beta}))= \xi_{\beta}({\cal L}(J_{\beta}))={\cal L}(\beta)$. \\
Montrons maintenant l'injectivit\'e. Soient $(\widetilde{E},\Gamma,\nu)$ et $(\widetilde{E},\Gamma',\nu')$ deux objets de ${\bf {\cal C}}_{O_{L_e}}(\widetilde{E})$, avec ${\bf M}(\nu)_{L_e}({\cal L}(\Gamma))={\cal L}(\beta)$ et ${\bf M}(\nu')_{L_e}({\cal L}(\Gamma'))={\cal L}(\beta')$, $\beta ,\beta' \in O_{L_e}$. D'apr\`es le lemme~\ref{homlem} ils sont isomorphes si et seulement si il existe $\gamma \in \mbox{Aut}_{{\BBf}_p}(\widetilde{E})$ tel que ${\bf M}(\gamma(p))_{L_e}({\cal L}(\beta')) \subset {\cal L}(\beta)$. Or, $\widetilde{E}/{\BBf}_p$ \'etant supersinguli\`ere, on a $\mbox{Aut}_{{\BBf}_p}(\widetilde{E})=\{ \pm 1 \}$. La multiplication par $(-1)$ sur $\widetilde{E}$ induit $-\mbox{Id}_{M}$ sur $M$, d'o\`u ${\bf M}([-1]_{\widetilde{E}}(p))({\cal L}(\beta'))={\cal L}(\beta')$. Donc $(\widetilde{E},\Gamma,\nu)$ et $(\widetilde{E},\Gamma',\nu')$ sont isomorphes si et seulement si ${\cal L}(\beta') \subset {\cal L}(\beta)$, c'est-\`a-dire $\beta = \beta'$.
\end{pf}

\medskip
Pour tout $\beta \in O_{L_e}$ on note ${\cal E}_{\beta}$ le sch\'ema elliptique sur $O_{L_e}$, unique \`a $O_{L_e}$-isomorphisme pr\`es, qui correspond par ${\bf ST}$ \`a un objet isomorphe dans ${\cal C}_{O_{L_e}}$ \`a un triplet $(\widetilde{E},J_{\beta}, \nu_{\beta})$, avec ${\bf M}(\nu_{\beta})_{L_e}({\cal L}(J_{\beta}))={\cal L}(\beta) \subset {\cal M}$. Soient $\beta , \beta' \in O_{L_e}$ ; on a un isomorphisme 
 $$\mbox{Hom}_{{\bf {\cal SE}}_{O_{L_e}}(\widetilde{E})} ({\cal E}_{\beta},{\cal E}_{\beta'}) \ \simeq \ \{\, \gamma \in \mbox{End}_{{\BBf}_p}(\widetilde{E})\  / \ {\bf M}(\gamma(p))_{L_e}({\cal L}(\beta')) \subset {\cal L}(\beta)\, \} $$

\begin{rem}
\label{supnonisogrmq}
{\em Lorsque $\gamma \in \mbox{End}_{{\BBf}_p}(\widetilde{E})$ le polyn\^ome caract\'eristique de ${\bf M}(\gamma(p))$ doit \^etre dans ${\BBq}[X]$. Dans la base $(e_1,e_2)$ de $M$ sa matrice doit s'\'ecrire sous la forme 
 $$\left ( 
\begin{array}{cc}
a & -pc \\
c & a 
\end{array}
\right )  \hspace{0.5cm} a,c \in {\BBz}_p$$
pour commuter avec $\varphi$, d'o\`u $a\in {\BBq}$ et $c^2\in {\BBq}\,$. Prenons $e=1$ et $\beta \in {\BBz}_p$ ; alors on a ${\bf M}(\gamma(p))({\cal L}(\beta)) \subset {\cal L}(0)$ si et seulement si $a\beta +c=0$. On en d\'eduit que si $\beta$ n'appartient \`a aucune extension quadratique de ${\BBq}$ les sch\'emas ${\cal E}_{0}$ et ${\cal E}_{\beta}$ ne sont pas ${\BBz}_p$-isog\`enes.}
\end{rem}

\begin{rem}
{\em Si $j(\widetilde{E})=0$ ou $1728$ on peut choisir un g\'en\'erateur du ${\BBz}_p[\varphi]$-module ${\bf M}(\widetilde{E}(p))$ qui est adapt\'e au groupe d'automorphismes de $\widetilde{E}$ (cf.~\ref{sec:MD}). Avec ce choix on a $j({\cal E}_{\beta})=0$ ou $1728$ si et seulement si $\beta=0$ (voir prop.~\ref{cpotsupprop}). Par analogie avec le cas ordinaire on appelle ${\cal E}_0$ le rel\`evement canonique de $\widetilde{E}$ sur $O_{L_e}$.} 
\end{rem}

\medskip
Dans la suite, pour chaque $e \in \{ 3,4,6 \}$, on s'int\'eresse au probl\`eme suivant : parmi les sch\'emas elliptiques sur $O_{L_e}$ relevant $\widetilde{E}$ d\'ecrits ci-dessus, quels sont ceux qui proviennent d'une courbe elliptique d\'efinie sur ${\BBq}_p$ ? Nous allons chercher ceux qui sont d\'efinis sur ${\BBq}_p$ et dont le d\'efaut de semi-stabilit\'e est $e$. Soit $\beta \in O_{L_e}$. Par le th\'eor\`eme de pleine fid\'elit\'e de Tate, le module de Dieudonn\'e filtr\'e $(M,{\cal L}(\beta))$ caract\'erise le ${\BBz}_p[G_{L_e}]$-module $T_p({\cal E}_{\beta})$ ; on a un isomorphisme canonique de $(\varphi, G_{K_e/L_e})$-modules filtr\'es
 $$(M\!\otimes_{{\BBz}_p}\!{\BBq}_{p^2}\, , \, {\cal L}(\beta)\!\otimes_{O_{L_e}}\!K_e) \ \simeq \ {\bf D}^*_{cris, K_e} (V_p({\cal E}_{\beta}))$$
Si ${\cal E}_{\beta}$ est d\'efinie sur ${\BBq}_p$ avec un d\'efaut de semi-stabilit\'e $e$ alors ${\bf D}^*_{cris, K_e} (V_p({\cal E}_{\beta}))$ devient un $(\varphi, G_{K_e/{\BBq}_p})$-module filtr\'e, c'est-\`a-dire qu'il est en plus muni d'une action de $G_{K_e/{\BBq}_{p^2}} = I(K_e/{\BBq}_p)$ compatible avec toutes les autres structures. On sait qu'alors cette action doit \^etre comme dans les objets ${\bf D^*_{pc}(e;0;\alpha)}$ de la liste ${\bf D^*}$ d\'ecrits en~\ref{sec:listp}, ce qui am\`ene \`a consid\'erer un crit\`ere galoisien de descente.

\subsection{Un crit\`ere galoisien de descente}
\label{sec:galdesc}

Soient $F\subset \overline{\BBq}_p$ une extension finie de ${\BBq}_p$ et $K\subset \overline{\BBq}_p$ une extension finie galoisienne de $F$ totalement ramifi\'ee, d'indice de ramification absolu $e(K)$, d'anneau des entiers $O_K$ et de corps r\'esiduel $k$. On note $G_{K/F}=\mbox{Gal}(K/F)$.

Soit $E/F$ une courbe elliptique qui acquiert bonne r\'eduction sur $K$. Alors le groupe $G_{K/F}$ op\`ere sur $E_K = E\times_{F}K$ via son action sur $K$ et cette action s'\'etend par fonctorialit\'e au mod\`ele de N\'eron de $E_K$. Comme $G_{K/F}$ agit trivialement sur $k$, son action sur la fibre sp\'eciale $\widetilde{E}_K$ s'effectue par des $k$-automorphismes de celle-ci, c'est-\`a-dire par un morphisme $G_{K/F} \rightarrow \mbox{Aut}_k(\widetilde{E}_K)$. En particulier, l'action de $G_K$ sur $T_p(E_K)$ s'\'etend en une action de $G_F$ de telle sorte que, pour $e(K) <p-1$, sur le module de Dieudonn\'e filtr\'e associ\'e au groupe $p$-divisible $E_K(p)$, on a une action de $G_{K/F}$ sur ${\bf M}_k(\widetilde{E}_K(p))$ qui pr\'eserve la filtration et qui provient d'un morphisme 
 $$G_{K/F} \rightarrow \mbox{Aut}_k(\widetilde{E}_K) \hookrightarrow \mbox{Aut}_k(\widetilde{E}_K(p)) \simeq \mbox{Aut}_{{\bf MD}_k} ({\bf M}_k(\widetilde{E}_K(p))) $$
R\'eciproquement, on a le th\'eor\`eme suivant :

\begin{thm}
\label{descthm}
Soit $E/K$ une courbe elliptique ayant bonne r\'eduction sur $K$. Alors $E$ est d\'efinie sur $F$ si et seulement si l'action de $G_K$ sur $T_p(E)$ s'\'etend en une action de $G_F$ qui provient de $k$-automorphismes de la fibre sp\'eciale $\widetilde{E}$ de $E$. \\
Plus pr\'ecis\'ement, il existe alors une courbe elliptique $E_0$ sur $F$ et un $K$-isomorphisme $\psi : E_0\times_{F}K \stackrel{\sim}{\rightarrow} E$ induisant un isomorphisme $G_F$-\'equivariant $\psi_p : T_p(E_0) \stackrel{\sim}{\rightarrow} T_p(E)$, o\`u $G_F$ agit naturellement sur $T_p(E_0)$ et sur $T_p(E)$ par l'action prolong\'ee que l'on s'est donn\'ee ; un tel couple $(E_0,\psi)$ est unique \`a $F$-isomorphisme pr\`es. De plus, le d\'efaut de semi-stabilit\'e de $E_0$ est \'egal \`a l'indice du noyau de $G_{K/F} \rightarrow \mbox{\em Aut}_k(\widetilde{E})$ dans $G_{K/F}$. 
\end{thm}

\begin{pf}
Pour tout $\omega \in G_{K/F}$ on note $E^{\omega}$ la courbe elliptique sur $K$ d\'eduite de $E$ par le changement de base $\mbox{Spec}(\omega ^{-1}) : \mbox{Spec}(K) \rightarrow \mbox{Spec}(K)$. L'association $E \mapsto E^{\omega}$ d\'efinit un foncteur ${\cal F}_{\omega}$ de la cat\'egorie des courbes elliptiques sur $K$ dans elle-m\^eme ; on a ${\cal F}_{\omega \tau} = {\cal F}_{\omega} \circ {\cal F}_{\tau}$ pour tous $\omega , \tau \in G_{K/F}$ et ${\cal F}_{1}=\mbox{Id}$ (le groupe $G_{K/F}$ agit sur la cat\'egorie). Rappelons le crit\`ere de Weil~: la courbe $E$ est d\'efinie sur $F$ si et seulement si il existe un ensemble de $K$-isomorphismes $f_{\omega} : E \rightarrow E^{\omega}$, $\omega \in G_{K/F}$, v\'erifiant 
 $$\makebox[2cm]{$(*)$}  f_{\omega\tau}= (f_{\tau})^{\omega}\circ f_{\omega} \makebox[1cm]{ } \forall \: \omega, \tau \in G_{K/F}$$
Il existe alors une courbe elliptique $E_0/F$ et un $K$-isomorphisme $\psi : E_0\times_{F}K \stackrel{\sim}{\rightarrow} E$ tel que $f_{\omega} = \psi^{\omega} \circ \psi^{-1}$ pour tout $\omega \in G_{K/F}$ ; le couple $(E_0, \psi)$ est unique \`a $K$-isomorphisme pr\`es (\cite{We} et \cite{La} thm.2G). Un ensemble $\{ f_{\omega} : E \rightarrow E^{\omega},\, \omega \in G_{K/F} \}$ v\'erifiant la condition $(*)$ est appel\'e un syst\`eme coh\'erent d'isomorphismes.

\medskip
Soit $\widehat{\omega}$ rel\`evement de $\omega \in G_{K/F}$ dans $G_F$ ; on a un isomorphisme de groupes de $E(\overline{\BBq}_p)$ dans $E^{\omega}(\overline{\BBq}_p)$ donn\'e par $\eta \mapsto \widehat{\omega} \circ \eta$ (il d\'epend du rel\`evement choisi) induisant une bijection ${\BBz}_p$-lin\'eaire $\alpha_{\widehat{\omega}} : T_p(E) \rightarrow T_p(E^{\omega})$. On se donne un morphisme $\rho : G_F \rightarrow \mbox{Aut}_{{\BBz}_p}(T_p(E))$ dont la restriction \`a $G_K$ est l'action naturelle. Alors, pour tout $\omega \in G_{K/F}$, l'isomorphisme ${\BBz}_p$-lin\'eaire 
 $$f_{\omega ,p} :  \left\{ 
\begin{array}{rcl}
T_p(E) & \rightarrow & T_p(E^{\omega}) \\
x     &   \mapsto    & \alpha_{\widehat{\omega}}(\rho(\widehat{\omega}^{-1})(x)) 
\end{array}
\right.$$
est $G_K$-\'equivariant et ne d\'epend pas du rel\`evement de $\omega$ dans $G_F$. Supposons qu'il existe un syst\`eme coh\'erent d'isomorphismes $\{ f_{\omega} : E \rightarrow E^{\omega},\, \omega \in G_{K/F} \}$ tel que $T_p(f_{\omega})=f_{\omega ,p}$ pour tout $\omega \in G_{K/F}$ ; soit $(E_0 ,\psi)$ le couple obtenu par le crit\`ere de Weil. Alors $\psi$ induit un isomorphisme $G_F$-\'equivariant $\displaystyle T_p(\psi)=\psi_p  : \underbrace{ T_p(E_0) }_{\mbox{\scriptsize action naturelle}} \stackrel{\sim}{\longrightarrow} \underbrace{ T_p(E) }_{\mbox{\scriptsize action \'etendue}}$.

\smallskip
Le choix d'une cl\^oture alg\'ebrique $\overline{\BBq}_p$ d\'efinit un foncteur contravariant $\Phi$ de la cat\'egorie ${\cal C}$ des sch\'emas en groupes finis et \'etales sur $\mbox{Spec}(K)$ dans la cat\'egorie ${\cal T}$ des groupes ab\'eliens finis munis d'une action de $G_K$, par $\Phi : X \mapsto X(\overline{\BBq}_p)=\mbox{Hom}_{K-alg}(B,\overline{\BBq}_p)$ o\`u $X=\mbox{Spec}(B)$~; le foncteur $\Psi : T \mapsto \mbox{Spec}((\mbox{Fcts}(T,\overline{\BBq}_p))^{G_K})$ de ${\cal T}$ dans ${\cal C}$ est un quasi-inverse. Alors, pour tout $\omega \in G_{K/F}$ et pour tout $T\in \mbox{Ob}({\cal T})$, on pose 
 $$ T^{\omega} = \Phi (\Psi(T)^{\omega}) = \mbox{Hom}_{K-alg} \Bigl( (\mbox{Fcts}(T,\overline{\BBq}_p))^{G_K}\otimes_{K^{\nearrow_{\!\omega^{-1}}}}\!\!K \, , \, \overline{\BBq}_p \Bigr)$$
L'objet $T^{\omega}$ est bien d\'efini et l'on obtient ainsi une action de $G_{K/F}$ sur la cat\'egorie ${\cal T}$ par ``transport de structure''. Par passage \`a la limite on d\'efinit $T^{\omega}$ o\`u $T$ est un module de Tate~; si $T= T_p(E)$ alors on a $(T_p(E))^{\omega} = T_p(E^{\omega})$. On a ainsi une notion de syst\`eme coh\'erent d'isomorphismes de ${\BBz}_p[G_K]$-modules. Alors le fait que l'action de $G_K$ sur $T_p(E)$ s'\'etend en une action de $G_F$ implique que le syst\`eme $\{ f_{\omega,p} : T_p(E) \rightarrow T_p(E^{\omega}),\, \omega \in G_{K/F} \}$ est coh\'erent (c'est purement formel).

\medskip
L'unicit\'e du mod\`ele de N\'eron ainsi que la pleine fid\'elit\'e du foncteur de Serre-Tate impliquent que l'existence d'un syst\`eme coh\'erent de $K$-isomorphismes $\{ f_{\omega} : E \rightarrow E^{\omega},\, \omega \in G_{K/F} \}$ \'equivaut aux donn\'ees suivantes :

\begin{itemize}
\item[(1)] un ensemble d'isomorphismes $\{ f_{\omega}(p) : E(p) \rightarrow E^{\omega}(p),\, \omega \in G_{K/F} \}$ de groupes $p$-divisibles sur $O_K$ tels que $f_{\omega\tau}(p)= (f_{\tau}(p))^{\omega}\circ f_{\omega}(p)$ pour tous $\omega, \tau \in G_{K/F}$ (coh\'erence)~;
\item[(2)] un ensemble d'isomorphismes $\{ \tilde{f}_{\omega} : \widetilde{E} \rightarrow \widetilde{E^{\omega}},\, \omega \in G_{K/F} \}$ de courbes elliptiques sur $k$ tels que $\tilde{f}_{\omega}(p) = \widetilde{f_{\omega}(p)}$ pour tout $\omega\in G_{K/F}$ (recollement). 
\end{itemize}

\noindent En effet, la condition de coh\'erence sur les $\tilde{f}_{\omega}$ est v\'erifi\'ee gr\^ace aux injections canoniques 
 $$\mbox{Hom}_{k}(\widetilde{E},\widetilde{E^{\omega}}) \hookrightarrow  \mbox{Hom}_{k}(\widetilde{E}(p),\widetilde{E^{\omega}}(p))$$
Le syst\`eme coh\'erent $\{ f_{\omega,p},\, \omega \in G_{K/F} \}$ fournit par le th\'eor\`eme de pleine fid\'elit\'e de Tate un syst\`eme coh\'erent $\{ f_{\omega}(p) : E(p) \rightarrow E^{\omega}(p), \, \omega \in G_{K/F} \}$ d'isomorphismes de groupes $p$-divisibles sur $O_K$, ce qui permet de satisfaire (1). Soit $\{ \widetilde{f_{\omega}(p)} : \widetilde{E}(p) \rightarrow \widetilde{E^{\omega}}(p),\, \omega \in G_{K/F} \}$ le syst\`eme coh\'erent d\'eduit sur les groupes $p$-divisibles des fibres sp\'eciales. L'extension $K/F$ \'etant totalement ramifi\'ee, on a $\widetilde{E^{\tau}}=\widetilde{E}$ pour tout $\tau \in G_{K/F}$, et la coh\'erence signifie que l'association $\tau \mapsto \widetilde{f_{\tau}(p)}$ est un morphisme $r_p : G_{K/F} \rightarrow \mbox{Aut}_k(\widetilde{E}(p)) \simeq ({\BBz}_p \otimes_{\BBz} \mbox{End}_{k}(\widetilde{E}))^{\times}$. Alors $E$ est d\'efinie sur $F$ si et seulement si $r_p$ provient par extension des scalaires d'un morphisme $r : G_{K/F} \rightarrow \mbox{Aut}_k(\widetilde{E})$. En effet, il suffit de poser $\tilde{f}_{\tau}=r(\tau)$ pour tout $\tau \in G_{K/F}$ : ce sont des $k$-isomorphismes de $\widetilde{E}$ qui v\'erifient, par construction, la condition (2).

\medskip
Soit $K_0 = \mbox{Frac}(W(k))$ l'extension maximale non ramifi\'ee contenue dans $K$. Prolonger l'action de $G_K$ \`a $G_F$ sur $V_p(E)$ revient \`a munir le $\varphi$-module filtr\'e $D= {\bf D}^*_{cris,K}(V_p(E))$ d'une structure d'objet de ${\bf MF}_{K/F}(\varphi)$, c'est-\`a-dire faire agir $K_0$-lin\'eairement $G_{K/F}$ sur $D$ de sorte que cette action commute avec $\varphi$ et respecte la filtration sur $D\otimes_{\!K_0}\!K$ ; l'objet $D$ est alors isomorphe dans ${\bf MF}_{K/F}(\varphi)$ \`a ${\bf D}^*_{cris,K/F}(V_p(E_0))$. En oubliant la filtration on obtient un morphisme $\nu : G_{K/F} \rightarrow \mbox{Aut}_{K_0[\varphi]}(D)$ et le diagramme suivant est commutatif :
 $$\xymatrix{
                     &  \mbox{Aut}_k(\widetilde{E}) \ar@{^{(}->} [d]  \\
G_{K/F} \ar[ur]^r \ar[r]^{r_p} \ar[dd]_{\nu}   &  \mbox{Aut}_k(\widetilde{E}(p)) \ar[d]^{\wr}  \\
                     &  \mbox{Aut}_{{\bf MD}_k} ({\bf M}_k(\widetilde{E}(p))) \ar@{^{(}->} [d]  \\
\mbox{Aut}_{K_0[\varphi]}(D) \ar[r]^-{\sim}_-{can}  & \mbox{Aut}_{K_0[\varphi]}(K_0\otimes_{\!W(k)}{\bf M}_k(\widetilde{E}(p)))  } $$
Soit $F'$ le sous-corps de $K$ fixe par $\mbox{Ker}(r) = \mbox{Ker}(\nu)$. La courbe $E_0/F$ acquiert bonne r\'eduction sur $F'$ puisque $V_p(E_0)$ est cristalline sur $F'$. Elle n'acquiert bonne r\'eduction sur aucun sous-corps strict de $F'$ contenu dans $F$ car $V_p(E_0)$ ne peut \^etre cristalline sur un tel corps : $\nu$ induit une injection $\mbox{Gal}(F'/F) \hookrightarrow \mbox{Aut}_{K_0[\varphi]}(D)$. Donc $\mbox{dst}(E_0)$ est \'egal au degr\'e de l'extension $F'/F$ qui est aussi l'indice de $\mbox{Ker}(r)$ dans $G_{K/F}$.
\end{pf}

\medskip
Le lemme qui suit permet de traiter des situations o\`u l'extension $K/F$ est totalement ramifi\'ee mais pas n\'ecessairement galoisienne. Soient $F_1,F_2$ deux corps tels que $F \subset F_i \subset K$, $F=F_1 \cap F_2$ et $K= F_1 F_2$ ; posons $G_{K/F_i}= \mbox{Gal}(K/F_i)$. Supposons l'extension $F_1/F$ galoisienne : $G_{K/F}$ est un produit semi-direct de $G_{K/F_2}$ par $G_{K/F_1}$.

\begin{lem}
\label{desclem}
Soient $F_1,F_2$ comme ci-dessus. Supposons $E$ d\'efinie sur $F_1$ et sur $F_2$, ce qui prolonge l'action de $G_K$ sur $T_p(E)$ en une action de $G_{F_1}$ et une de $G_{F_2}$. Si l'action de $G_K$ s'\'etend sur $T_p(E)$ en une action de $G_F$ qui co\"{\i}ncide avec celles de $G_{F_1}$ et $G_{F_2}$, alors $E$ est d\'efinie sur $F$.
\end{lem}

\begin{pf}
On dispose de deux syst\`emes coh\'erents d'isomorphismes $\{ f_{\omega_i} : E \rightarrow E^{\omega_i},\, \omega_i \in G_{K/F_i} \}$, $i=1,2$, et il s'agit de montrer l'existence d'un syst\`eme coh\'erent $\{ f_{\omega} : E \rightarrow E^{\omega},\, \omega \in G_{K/F} \}$. Tout \'el\'ement $\omega$ de $G_{K/F}$ s'\'ecrivant de mani\`ere unique $\omega= \omega_1 \omega_2$ avec $\omega_i \in G_{K/F_i}$, on pose $f_{\omega} = f_{\omega_1 \omega_2} = (f_{\omega_2})^{\omega_1} \circ f_{\omega_1}$. Comme l'action de $G_K$ sur $T_p(E)$ s'\'etend \`a $G_F$ de sorte qu'elle co\"{\i}ncide avec celles de $G_{F_1}$ et de $G_{F_2}$, le syst\`eme de ${\BBz}_p[G_K]$-isomorphismes $\{ T_p(f_{\omega}) : T_p(E) \rightarrow T_p(E^{\omega}),\, \omega \in G_{K/F} \}$ est coh\'erent. Alors l'injection canonique $\mbox{Hom}_K (E,E') \hookrightarrow  \mbox{Hom}_{{\BBz}_p[G_K]} (T_p(E),T_p(E'))$ pour deux courbes elliptiques $E$ et $E'$ sur $K$ permet d'obtenir la coh\'erence du syst\`eme $\{ f_{\omega},\, \omega \in G_{K/F} \}$ \`a partir de celle de $\{ T_p(f_{\omega}),\, \omega \in G_{K/F} \}$.
\end{pf}

\begin{rem}
{\em Le th\'eor\`eme \ref{descthm} et le lemme \ref{desclem} sont \'egalement valables pour des vari\'et\'es ab\'eliennes de dimension relative quelconque (les preuves sont les m\^emes).} 
\end{rem}

\subsection{Courbes potentiellement supersinguli\`eres}
\label{sec:cpotsup}

Nous allons appliquer le th\'eor\`eme de descente~\ref{descthm} \`a notre situation. Comme $e \in \{ 3,4,6\}$ on voit que l'on doit partir d'une courbe $\widetilde{E}/{\BBf}_p$ ayant suffisament d'automorphismes (d\'efinis sur ${\BBf}_{p^2}$), c'est-\`a-dire telle que $[\zeta_e] \in \mbox{Aut}_{{\BBf}_{p^2}}(\widetilde{E})$. On fixe $\widetilde{E}/{\BBf}_p$ supersinguli\`ere d'invariant modulaire $\tilde{ \mbox{\j}}(e)$ avec $\tilde{ \mbox{\j}}(3)=\tilde{ \mbox{\j}}(6)=0$ et $\tilde{ \mbox{\j}}(4)=1728$ ; alors $e\mid p+1$ et pour satisfaire la condition $e<p-1$ il faut exclure le cas $(e,p)=(6,5)$ (voir cependant la rmq.~\ref{allsupsingrmq}(ii)). De plus, l'extension $L_e/{\BBq}_p$ n'\'etant pas galoisienne, il faut monter sur ${\BBq}_{p^2}$ et travailler avec l'extension $K_e/{\BBq}_{p^2}$~; on utilise alors le lemme~\ref{desclem}.

\medskip
On note $O_{L_e}'$ le sous-ensemble de $O_{L_e}$ correspondant aux invariants $\beta \in O_{L_e}$ pour lesquels ${\cal E}_{\beta}$ peut \^etre d\'efinie sur ${\BBq}_p$ avec un d\'efaut de semi-stabilit\'e $e$ et $O_{L_e}''$ l'ensemble des classes d'isomorphisme de courbes elliptiques sur ${\BBq}_p$ qui prolongent un sch\'ema ${\cal E}_{\beta}$ avec un d\'efaut de semi-stabilit\'e $e$. L'ensemble $O_{L_e}''$ consiste en la donn\'ee d'un \'el\'ement $\beta \in O_{L_e}'$ avec en plus celle d'une action prolong\'ee de $G$ sur $T_p({\cal E}_{\beta})$ ; on a bien s\^ur une fl\`eche naturelle surjective $\lambda_e : O_{L_e}'' \rightarrow O_{L_e}'$.

\begin{prop}
\label{cpotsupprop}
Soit $e \in \{ 3,4,6 \}$ tel que $e\mid p+1$ et $e < p-1$. Soit $\widetilde{E}/{\BBf}_p$ supersinguli\`ere d'invariant modulaire $\tilde{ \mbox{\j}}(e)$ avec $\tilde{ \mbox{\j}}(3)=\tilde{ \mbox{\j}}(6)=0$ et $\tilde{ \mbox{\j}}(4)=1728$. On choisit un g\'en\'erateur du ${\BBz}_p[\varphi]$-module ${\bf M}(\widetilde{E}(p))$ adapt\'e au groupe d'automorphismes de $\widetilde{E}$ pour param\'etrer les rel\`evements de $\widetilde{E}$ en un sch\'ema elliptique sur $O_{L_e}$. Alors : \\
1) $j({\cal E}_{\beta}) = 0$ ou $1728$ si et seulement si $\beta =0$. \\
2) $O_{L_e}' = {\BBz}_p \pi_e \cup {\BBz}_p \pi_e^{e-3}$, de sorte que $O_{L_e}'$ s'identifie \`a un sous-${\BBz}_p$-module de rang $1$ de $O_{L_e}$ si $e=4$ et \`a ${\BBz}_p \cup {\BBz}_p$ si $e=3$ ou $6$. \\
3) $O_{L_e}''$ s'identifie \`a ${\BBz}_p \sqcup {\BBz}_p$ (r\'eunion disjointe) et les fibres de $\lambda_e$ sont de cardinal $2$ si $\beta =0$ ou $e=4$, de cardinal $1$ sinon. 
\end{prop}

\begin{pf}
Fixons un g\'en\'erateur $e_1$ du ${\BBz}_p[\varphi]$-module $M={\bf M}(\widetilde{E}(p))$ adapt\'e au groupe d'automorphismes de $\widetilde{E}$ ; rappelons qu'alors $(e_1,\varphi(e_1)=e_2)$ est une base de $M$ qui, apr\`es extension des scalaires \`a ${\BBz}_{p^2}$, diagonalise l'action de $\xi_e = {\bf M}_{{\BBf}_{\!p^2}}([\zeta_e](p))$ (\ref{sec:MD}). Soit $\beta \in O_{L_e}$~; notons encore ${\cal E}_{\beta}$ le sch\'ema ${\cal E}_{\beta}\!\times_{O_{L_e}}\!O_{K_e}$ o\`u $O_{K_e}$ est l'anneau des entiers de $K_e$. \\
1) On a $j({\cal E}_{\beta}) =0$ (resp. $1728$) si et seulement si $[\zeta_e] \in \mbox{Aut}_{{\BBf}_{\!p^2}}(\widetilde{E})$ se rel\`eve dans $\mbox{Aut}_{O_{K_e}}({\cal E}_{\beta})$ avec $e=3$ ou $6$ (resp. $e=4$), i.e. si et seulement si $\xi_e$ stabilise apr\`es extension des scalaires la filtration ${\cal L}(\beta)\!\otimes_{O_{L_e}}\!O_{K_e}=(e_1 \otimes 1 + \beta \cdot e_2 \otimes \pi_e^{1-e})O_{K_e}$. Cette condition s'\'ecrit $e_1\otimes \zeta_e^{\eta} +\beta \cdot e_2\otimes \zeta_e^{-\eta} \pi_e^{1-e} \in {\cal L}(\beta)\!\otimes_{O_{L_e}}\!O_{K_e}$ avec $\eta =\pm 1$, ce qui \'equivaut \`a $\beta =0$. \\
2) D'apr\`es le th\'eor\`eme~\ref{descthm} et le lemme~\ref{desclem}, la courbe ${\cal E}_{\beta}$ est d\'efinie sur ${\BBq}_p$ avec un d\'efaut de semi-stabilit\'e $e$ si et seulement si l'action de $G_{K_e}$ sur $T_p({\cal E}_{\beta})$ s'\'etend en une action de $G$, dont la restriction \`a $G_{{\BBq}_{p^2}}$ induit une injection $<\tau_e> \hookrightarrow \mbox{Aut}_{{\bf MD}_{{\BBf}_{\!p^2}}}(M\!\otimes_{{\BBz}_p}\!{\BBz}_{\!p^2})$ pr\'eservant la filtration et provenant d'une injection $<\tau_e> \hookrightarrow \mbox{Aut}_{{\BBf}_{\!p^2}}(\widetilde{E})$. Cette derni\`ere condition \'equivaut \`a $\tau_e=\xi_e$ ou $\xi_e^{-1}$, d'o\`u $\tau_e e_1 = \zeta_e^{\epsilon} e_1$ et $\tau_e e_2 = \zeta_e^{-\epsilon} e_2$ avec $\epsilon \in \{ \pm 1 \}$.  Maintenant \'ecrivons que l'action de $G_{K_e/{\BBq}_p}$ \'etendue par semi-lin\'earit\'e stabilise ${\cal L}(\beta)\!\otimes_{O_{L_e}}\!O_{K_e}$. C'est automatique en ce qui concerne $\omega$, puisque $\beta \in O_{L_e}$ (${\cal E}_{\beta}$ est d\'ej\`a d\'efinie sur $O_{L_e}$) ; pour $\tau_e$ cela \'equivaut \`a $\tau_e(\pi_e^{-2\epsilon +1}\beta) = \pi_e^{-2\epsilon +1}\beta$, c'est-\`a-dire $\pi_e^{-2\epsilon +1} \beta \in {\BBq}_{p^2} \cap L_e ={\BBq}_p$. Finalement, comme $\beta \in O_{L_e}$, on obtient : pour  $\epsilon = 1$, $\beta \in {\BBz}_p \pi_e$ ; pour $\epsilon = -1$, $\beta \in {\BBz}_p\pi_e^{e-3}$. \\
3) Cela d\'ecoule de l'assertion d'unicit\'e du th\'eor\`eme~\ref{descthm}. Si $e\in \{ 3,6\}$ alors ${\BBz}_p \pi_e \cap {\BBz}_p\pi_e^{e-3}= \{ 0 \}$ ; pour $\beta \in {\BBz}_p \pi_e \backslash \{ 0 \}$ (resp. $\beta \in {\BBz}_p\pi_e^{e-3} \backslash \{ 0 \}$) fix\'e, il y a, \`a ${\BBq}_p$-isomorphisme pr\`es, une seule courbe prolongeant ${\cal E}_{\beta}$ sur ${\BBq}_p$ avec un d\'efaut de semi-stabilit\'e \'egal \`a $e$ et elle correspond \`a une action prolong\'ee de $\tau_e$ avec $\epsilon =1$ (resp. $\epsilon = -1$). Si $e=4$ ou $\beta =0$ il y en a deux, l'une correspondant \`a une action prolong\'ee avec $\epsilon = 1$, l'autre avec $\epsilon = -1$.
\end{pf}

\begin{rem}
{\em Si $e=3$ (resp. $e=6$), le sch\'ema ${\cal E}_0$ se prolonge aussi en un sch\'ema elliptique ${\cal A}_0$ sur ${\BBz}_p$ ; comme $j({\cal A}_0)=0$ les deux courbes prolongeant ${\cal E}_0$ sur ${\BBq}_p$ avec un d\'efaut de semi-stabilit\'e $3$ (resp. $6$) sont les tordues de ${\cal A}_0\!\times_{\!{\BBz}_p}\!{\BBq}_p = A_0$ correspondant aux \'el\'ements $p^4$ et $p^2$ (resp. $(-p)^5$ et $(-p)$) de ${\BBq}_p^{\times}/({\BBq}_p^{\times})^6 \simeq \mbox{Twist}((A_0,{\bf 0}),{\BBq}_p)$. Si $e=4$, le sch\'ema ${\cal E}_0$ se prolonge aussi en un sch\'ema elliptique ${\cal B}_0$ sur ${\BBz}_p$ ; comme $j({\cal B}_0)=1728$ les deux courbes prolongeant ${\cal E}_0$ sur ${\BBq}_p$ avec un d\'efaut de semi-stabilit\'e $4$ sont les tordues de ${\cal B}_0\!\times_{\!{\BBz}_p}\!{\BBq}_p = B_0$ correspondant aux \'el\'ements $(-p)^3$ et $(-p)$ de ${\BBq}_p^{\times}/({\BBq}_p^{\times})^4 \simeq \mbox{Twist}((B_0,{\bf 0}),{\BBq}_p)$. Pour tous ces cas voir les exemples donn\'es en~\ref{sec:exsupp}.}
\end{rem}

\medskip
On note $E_{\alpha, \epsilon}$ les courbes elliptiques sur ${\BBq}_p$ correspondant aux \'el\'ements de $O_{L_e}''$, o\`u $\alpha \in {\BBz}_p$ et $\epsilon \in \{ \pm 1 \}$ sont tels que $E_{\alpha, \epsilon}\times_{{\BBq}_p} L_e \simeq {\cal E}_{\beta}\times_{O_{L_e}}L_e$ avec $\beta = \alpha \pi_e$ et une action \'etendue avec $\epsilon =1$ ou bien $\beta = \alpha \pi_e^{e-3}$ et une action \'etendue avec $\epsilon = -1$.

\begin{rem}
\label{supisogrmq}
{\em Si $\alpha \in {\BBz}_p^{\times}$ le morphisme $\varphi$ de $M$ envoie, apr\`es extension des scalaires, ${\cal L}(\alpha \pi_e^{e-3})$ dans ${\cal L}(\alpha^{-1} \pi_e)$ ; comme il provient d'un morphisme de $\widetilde{E}$, \`a savoir le Frobenius, on en d\'eduit que les sch\'emas ${\cal E}_{\alpha^{-1} \pi_e}$ et ${\cal E}_{\alpha \pi_e^{e-3}}$ sont $O_{L_e}$-isog\`enes. De plus, on v\'erifie que le morphisme $\varphi : {\bf D}^*_{cris,K_e}(V_p(E_{\alpha,-1})) = D_1 \rightarrow D_2 = {\bf D}^*_{cris,K_e}(V_p(E_{\alpha^{-1},1}))$ d'objets de ${\bf MF}_{K_e}(\varphi)$ commute \`a l'action de $G_{K_e/{\BBq}_p}$, de sorte que c'est un morphisme dans ${\bf MF}_{K_e/{\BBq}_p}(\varphi)$. L'injection 
 $$ \mbox{Hom}_{{\BBq}_p}( E_{\alpha^{-1},1} , E_{\alpha,-1} ) \hookrightarrow \mbox{Hom}_{{\BBq}_p[G]}( V_p(E_{\alpha^{-1},1}) , V_p(E_{\alpha,-1}) )  \simeq \mbox{Hom}_{{\bf MF}_{K_e/{\BBq}_p}(\varphi)} (D_1,D_2)$$
montre que cette $O_{L_e}$-isog\'enie est d\'efinie sur ${\BBq}_p$. Donc les courbes $E_{\alpha^{-1},1}$ et $E_{\alpha,-1}$ sont ${\BBq}_p$-isog\`enes.}
\end{rem}

\begin{rem}
\label{allsupsingrmq}
{\em (i) Soit $\widetilde{E}'/{\BBf}_p$ la tordue de $\widetilde{E}$ sur ${\BBf}_{p^2}$ et, pour tout $\beta \in O_{L_e}$, soit ${\cal E}_{\beta}'/O_{L_e}$ le sch\'ema obtenu en tordant ${\cal E}_{\beta}$ sur $K_e$. Alors les rel\`evements de $\widetilde{E}'$ en un sch\'ema elliptique sur $O_{L_e}$ sont, \`a isomorphisme pr\`es, les ${\cal E}_{\beta}'$. De plus, ${\cal E}_{\beta}'$ est d\'efinie sur ${\BBq}_p$ avec un d\'efaut de semi-stabilit\'e $e$ si et seulement si ${\cal E}_{\beta}$ l'est ; les courbes sur ${\BBq}_p$ obtenues \`a partir de $\widetilde{E}'$ sont les tordues sur ${\BBq}_{p^2}$ de celles obtenues \`a partir de $\widetilde{E}$. \\
(ii) Si $(e,p)=(6,5)$ on se ram\`ene \`a la situation $(e,p)=(3,5)$ en tordant sur l'extension quadratique ${\BBq}_5(\pi_2)/{\BBq}_5$. En effet, si $E/{\BBq}_p$ a potentiellement bonne r\'eduction avec $\mbox{dst}(E)=6$, alors sa tordue $E'/{\BBq}_p$ sur ${\BBq}_p(\pi_2)$ est telle que $\mbox{dst}(E')=3$, et vice versa (et, sur ${\BBq}_p(\pi_6)$, elles ont la m\^eme fibre sp\'eciale). \\
On obtient ainsi toutes les courbes elliptiques sur ${\BBq}_p$, \`a ${\BBq}_p$-isomorphisme pr\`es, qui sont potentiellement supersinguli\`eres avec un d\'efaut de semi-stabilit\'e $e\in \{ 3,4,6 \}$.}
\end{rem}

\begin{rem}
\label{twistsupisogrmq}
{\em Soit $\gamma : \widetilde{E} \rightarrow \widetilde{E}'$ un isomorphisme d\'efini sur ${\BBf}_{p^2}$ ; via l'isomorphisme $\mbox{End}_{{\BBf}_{\!p^2}}(\widetilde{E}) \rightarrow \mbox{Hom}_{{\BBf}_{\!p^2}}(\widetilde{E}, \widetilde{E}')$ donn\'e par $\psi \mapsto \gamma \circ \psi$, les \'el\'ements de $\mbox{Hom}_{{\BBf}_p}(\widetilde{E},\widetilde{E}')$ correspondent aux $\psi \in \mbox{End}_{{\BBf}_{\!p^2}}(\widetilde{E})$ tels que $\psi^{\sigma}=-\psi$, o\`u $\sigma$ est le Frobenius absolu. Soit $\psi_e = [\zeta_e] - [\zeta_e^{-1}] \in \mbox{End}_{{\BBf}_{\!p^2}}(\widetilde{E})$. Comme $\psi_e^{\sigma}=-\psi_e$ l'isog\'enie $\gamma \circ \psi_e$ est d\'efinie sur ${\BBf}_p$, d'o\`u un morphisme ${\bf M}(\gamma \circ \psi_e) ={\bf M}_{{\BBf}_{\!p^2}}(\psi_e) \circ {\bf M}_{{\BBf}_{\!p^2}}(\gamma) : {\bf M}(\widetilde{E}'(p)) \rightarrow {\bf M}(\widetilde{E}(p))$. Du fait que ${\bf M}_{{\BBf}_{\!p^2}}(\psi_e)$ envoie, apr\`es extension des scalaires, ${\cal L}(\beta)$ dans ${\cal L}(-\beta)$, on d\'eduit que $\gamma \circ \psi_e$ se rel\`eve en une $O_{L_e}$-isog\'enie $(\gamma \circ \psi_e)_{\beta} : {\cal E}_{\beta} \rightarrow {\cal E}_{-\beta}'$ pour tout $\beta \in O_{L_e}$. Lorsque $\beta \in O_{L_e}'$ et pour une action \'etendue avec le m\^eme invariant $\epsilon$, on v\'erifie que le morphisme associ\'e \`a $(\gamma \circ \psi_e)_{\beta}$ sur les objets de ${\bf MF}_{K_e}(\varphi)$ correspondants commute \`a l'action de $G_{K_e/{\BBq}_p}$. Donc, si pour tous $\alpha \in{\BBz}_p$ et $\epsilon \in \{ \pm 1 \}$ on note $E_{\alpha, \epsilon}'$ la tordue sur ${\BBq}_{p^2}$ de $E_{\alpha, \epsilon}$, on obtient que les courbes $E_{\alpha, \epsilon}$ et $E_{-\alpha, \epsilon}'$ sont ${\BBq}_p$-isog\`enes.}
\end{rem}

\subsection{Sur les cas ordinaires}
\label{sec:cpotord}

Si $\widetilde{E}/{\BBf}_p$ est ordinaire on a $a_p(\widetilde{E})\in {\BBz}_p^{\times}$ et il existe une base $(e_1,e_2)$ de ${\bf M}(\widetilde{E}(p))$ qui diagonalise $\varphi$. Si $e_1$ est un vecteur propre dont la valeur propre est une unit\'e, les rel\`evements de ${\bf M}(\widetilde{E}(p))$ en un module de Dieudonn\'e sur $O_{L_e}$ correspondent bijectivement aux filtrations ${\cal L}(\beta) = (\beta \cdot e_1 \otimes \pi_e^{1-e} + e_2 \otimes 1)O_{L_e}$ avec $\beta \in O_{L_e}$. On obtient avec des m\'ethodes tout \`a fait similaires les r\'esultats qui suivent.

\begin{prop}
\label{schellordprop}
Soit $e<p-1$. Soit $\widetilde{E}/{\BBf}_p$ une courbe elliptique ordinaire d'invariant modulaire $\tilde{ \mbox{\j}}$ ; on pose $m(\tilde{ \mbox{\j}})=1$ si $\tilde{ \mbox{\j}} \not\in \{ 0,1728 \}$, $m(1728)=2$ et $m(0)=3$. Via le choix d'une ${\BBz}_p$-base de diagonalisation de $\varphi$ dans ${\bf M}(\widetilde{E}(p))$, l'association 
 $$ \left \{
   \begin{array}{rll}
   {\bf {\cal C}}_{O_{L_e}}(\widetilde{E})  &  \rightarrow   &  O_{L_e}   \\
   (\widetilde{E} , \Gamma , \nu )          &  \mapsto       &  \beta \makebox[1.5cm]{tel que} {\cal L}(\beta)= {\bf M}(\nu)_{L_e}({\cal L}(\Gamma))
   \end{array}
   \right . $$
induit une bijection entre les classes d'isomorphisme dans ${\bf {\cal C}}_{O_{L_e}}(\widetilde{E})$ et l'ensemble $O_{L_e} / {\bf \sim}\,$, avec $x {\bf \sim} y$ si et seulement si $x^{m(\tilde{ \mbox{\j}})} = y^{m(\tilde{ \mbox{\j}})}$.
\end{prop}

En composant avec le foncteur ${\bf ST}$ on obtient une bijection entre les classes d'isomorphisme dans ${\bf {\cal SE}}_{O_{L_e}}(\widetilde{E})$ et $O_{L_e}/{\bf \sim}\,$. Le choix d'une autre ${\BBz}_p$-base de diagonalisation de $\varphi$ dans ${\bf M}(\widetilde{E}(p))$ change l'invariant $\beta$ en $\eta \beta$ avec $\eta \in {\BBz}_p^{\times}$. Quand $\tilde{ \mbox{\j}} \in \{ 0,1728 \}$ ces classes sont param\'etr\'ees par un quotient de $O_{L_e}$ parce que, contrairement au cas supersingulier, le groupe des ${\BBf}_p$-automorphismes de $\widetilde{E}$ est alors strictement plus grand que $\{ \pm 1\}$. 

Signalons que les rel\`evements d'une courbe elliptique ordinaire ont d\'ej\`a \'et\'e \'etudi\'es : voir par exemple \cite{Me}, Appendix (en particulier la prop.3.2.), o\`u les m\'ethodes utilis\'ees n'imposent pas de restriction sur la ramification.

Pour tout $\beta \in O_{L_e}/{\bf \sim}$ on note ${\cal E}_{\beta}$ le sch\'ema elliptique sur $O_{L_e}$ qui correspond par ${\bf ST}$ \`a un triplet isomorphe dans ${\cal C}_{O_{L_e}}$ \`a un $(\widetilde{E},J_{\beta},\nu_{\beta})$ avec ${\bf M}(\nu_{\beta})_{L_e}({\cal L}(J_{\beta}))={\cal L}(\beta) \subset {\cal M}$. Le lemme~\ref{homlem} implique que les assertions suivantes sont \'equivalentes :  
\begin{itemize}
\item[(i)] $\beta =0$
\item[(ii)] le Frobenius de $\widetilde{E}$ se rel\`eve en un morphisme de ${\cal E}_{\beta}$ 
\item[(iii)] si $j(\widetilde{E})=0$ ou $1728$ alors $j({\cal E}_{\beta})=0$ ou $1728$
\item[(iv)] la suite exacte $(*_{ord})$ associ\'ee \`a ${\cal E}_{\beta}$ est scind\'ee (cf.~\ref{sec:notsgeom}).
\end{itemize}
Donc ${\cal E}_0/O_{L_e}$ est le rel\`evement canonique de $\widetilde{E}$ (cf.~\cite{Me}, App., cor.1.2 et 1.3).

\begin{rem}
\label{canordisogrmq}
{\em Soient $\widetilde{E}/{\BBf}_p$ et $\widetilde{E}'/{\BBf}_p$ deux courbes elliptiques ordinaires ; soient ${\cal E}_0/O_{L_e}$ et ${\cal E}_0'/O_{L_e}$ les rel\`evements canoniques de $\widetilde{E}$ et $\widetilde{E}'$ respectivement. Alors $\mbox{Hom}_{O_{L_e}}({\cal E}_0,{\cal E}_0') \simeq \mbox{Hom}_{{\BBf}_p}(\widetilde{E},\widetilde{E}')$, i.e. toute ${\BBf}_p$-isog\'enie $\widetilde{E} \rightarrow \widetilde{E}'$ se rel\`eve en une $O_{L_e}$-isog\'enie ${\cal E}_0 \rightarrow {\cal E}_0'$.} 
\end{rem}

\begin{rem}
\label{ordnonisogrmq}
{\em Soient $\beta, \beta' \in O_{L_e}$ tels que $\beta \beta'\neq 0$ et soit $\psi : ({\bf M}(\widetilde{E}(p)),{\cal L}(\beta)) \rightarrow ({\bf M}(\widetilde{E}(p)),{\cal L}(\beta'))$ un morphisme de modules de Dieudonn\'e filtr\'es. Dans une base de ${\bf M}(\widetilde{E}(p))$ qui diagonalise $\varphi$, la matrice de $\psi$ est de la forme $\mbox{Diag}(a,d)$ avec $a,d \in {\BBz}_p$ tels que $a\beta = d\beta'$. Si $\psi$ provient d'un \'el\'ement de $\mbox{End}_{{\BBf}_p}(\widetilde{E})$ son polyn\^ome caract\'eristique est dans ${\BBq}[X]$, d'o\`u $[\,{\BBq}(a,d):{\BBq}\,] \leq 2$. Prenons $\beta = \alpha \pi_e^i$ et $\beta' = \alpha' \pi_e^i$ avec $i \in {\BBn}$, $\alpha' \in {\BBz} \backslash \{ 0 \}$ et $\alpha \in {\BBz}_p \backslash \{ 0 \}$ tel que $[\,{\BBq}(\alpha):{\BBq}\,] > 2$ ; alors les sch\'emas ${\cal E}_{\beta}$ et ${\cal E}_{\beta'}$ ne sont pas $O_{L_e}$-isog\`enes.}
\end{rem}

\medskip
Soit $e \in \{ 3,4,6 \}$ tel que $e<p-1$. On prend maintenant $\widetilde{E}/{\BBf}_p$ ordinaire d'invariant modulaire $\tilde{ \mbox{\j}}(e)$ avec $\tilde{ \mbox{\j}}(3)=\tilde{ \mbox{\j}}(6)=0$ et $\tilde{ \mbox{\j}}(4)=1728$ ; dans ce cas $e\mid p-1$, $[\zeta_e] \in \mbox{Aut}_{{\BBf}_p}(\widetilde{E})$ et cet automorphisme commute avec le Frobenius. Pour satisfaire la condition $e<p-1$ il faut \'ecarter les valeurs $(e,p)=(4,5)$ et $(e,p)=(6,7)$. En appliquant le th\'eor\`eme~\ref{descthm} on trouve que ${\cal E}_{\beta}$ est d\'efinie sur ${\BBq}_p$ avec un d\'efaut de semi-stabilit\'e $e$ si et seulement si $\beta \in {\BBz}_p \pi_e^{e-3}$ (correspondant \`a une action \'etendue par $\tau_e = \xi_e$) ou bien $\beta \in {\BBz}_p \pi_e$ (correspondant \`a une action \'etendue par $\tau_e = \xi_e^{-1}$).

\begin{rem}
\label{allordrmq}
{\em (i) Si $(e,p)=(6,7)$ on se ram\`ene \`a la situation $(e,p)=(3,7)$ en tordant sur ${\BBq}_7(\pi_2)$, cf. rmq.~\ref{allsupsingrmq}(ii). \\
(ii) Par contre, si le d\'efaut de semi-stabilit\'e d'une courbe elliptique est $4$, il reste inchang\'e par torsion quadratique. Donc si $(e,p)=(4,5)$ nos m\'ethodes ne s'appliquent pas (du moins pas avec l'utilisation des modules de Dieudonn\'e filtr\'es). Dans ce cas nous faisons appel \`a \cite{Me}, Appendix ; en particulier, la prop. 3.3 que l'on y trouve permet de d\'eduire les analogues des remarques \ref{canordisogrmq} et \ref{ordnonisogrmq} ci-dessus.}
\end{rem}

\begin{rem}
\label{potcanordisogrmq}
{\em Soient $\widetilde{E}/{\BBf}_p$ et $\widetilde{E}'/{\BBf}_p$ deux courbes elliptiques ordinaires qui sont ${\BBf}_p$-isog\`enes, c'est-\`a-dire telles que $a_p(\widetilde{E})=a_p(\widetilde{E}')$. Soient $\beta,\beta' \in O_{L_e}$ tels que $\beta,\beta' \in {\BBz}_p \pi_e^{e-3}$ ou bien $\beta,\beta' \in {\BBz}_p \pi_e$~; notons ${\cal E}_{\beta}/O_{L_e}$ et ${\cal E}_{\beta'}'/O_{L_e}$ les sch\'emas elliptiques relevant $\widetilde{E}$ et $\widetilde{E}'$ respectivement, ainsi que ${E}_{\beta,\epsilon}/{\BBq}_p$ et ${E}_{\beta',\epsilon}'/{\BBq}_p$ les courbes qui les prolongent avec le m\^eme invariant $\epsilon \in \{ \pm 1\}$. Du fait que tout morphisme commutant avec les Frobenii commute avec $\tau_e$ dans les $(\varphi,G_{K_e/{\BBq}_p})$-modules filtr\'es associ\'es, on d\'eduit que $\mbox{Hom}_{{\BBq}_p}({E}_{\beta,\epsilon},{E}_{\beta',\epsilon}') \simeq \mbox{Hom}_{O_{L_e}}({\cal E}_{\beta},{\cal E}_{\beta'}')$, i.e. toute $O_{L_e}$-isog\'enie ${\cal E}_{\beta} \rightarrow {\cal E}_{\beta'}'$ se prolonge en une ${\BBq}_p$-isog\'enie ${E}_{\beta,\epsilon} \rightarrow {E}_{\beta',\epsilon}'$. En particulier, les courbes ${E}_{0,\epsilon}$ et ${E}_{0,\epsilon}'$ sont ${\BBq}_p$-isog\`enes.}
\end{rem}

\section{Les ${\BBq}_p[G]$-modules provenant d'une courbe elliptique sur ${\BBq}_p$}
\label{sec:pmodellip}

\subsection{R\'esultat et cons\'equences}
\label{sec:resultp}

Nous allons maintenant donner des conditions n\'ecessaires et suffisantes pour qu'une repr\'esentation $p$-adique $V_p$ de $G$ de dimension $2$ provienne d'une courbe elliptique sur ${\BBq}_p$, i.e. pour qu'il existe  $E/{\BBq}_p$ telle que $V_p \simeq V_p(E)$ en tant que ${\BBq}_p[G]$-modules.

Soit $E/{\BBq}_p$ une courbe elliptique. Tout d'abord, on sait que la repr\'esentation $V_p(E)$ est potentiellement semi-stable. Ensuite, la repr\'esentation de Weil-Deligne ${\bf W}^*_p(V_p(E))$ associ\'ee \`a $V_p(E)$ v\'erifie les conditions $(1^{\circ})$, $(2^{\circ})$ et $(3^{\circ})$ du th\'eor\`eme~\ref{lmodellipthm} (compatibilit\'e). Enfin, on sait que ${\bf D}^*_{pst}(V_p(E))$ est de type Hodge-Tate $(0,1)$ ; on dira que $V_p(E)$ est de type Hodge-Tate $(0,1)$, le foncteur ${\bf D}^*_{pst}$ \'etant sous-entendu.

L\`a encore, le r\'esultat est que ces conditions n\'ecessaires sont aussi suffisantes : d'une part les $(\varphi,N,G)$-modules filtr\'es faiblement admissibles v\'erifiant ces conditions sont exactement ceux de la liste ${\bf D^*}$, et d'autre part tous les objets de cette liste proviennent d'une courbe elliptique sur ${\BBq}_p$.

\begin{thm}
\label{pmodellipthm}
Soit $V_p$ une repr\'esentation $p$-adique de $G$ de dimension $2$. Les assertions suivantes sont \'equivalentes : 
\begin{itemize}
\item[(1)] il existe une courbe elliptique $E$ sur ${\BBq}_p$ telle que $V_p(E)$ soit isomorphe \`a $V_p$, 
\item[(2)] $V_p$ est potentiellement semi-stable de type Hodge-Tate $(0,1)$ et ${\bf W}^*_p(V_p)$ v\'erifie les conditions $(1^{\circ})$, $(2^{\circ})$ et $(3^{\circ})$ du th\'eor\`eme~\ref{lmodellipthm}, 
\item[(3)] $V_p$ est potentiellement semi-stable et ${\bf D}^*_{pst}(V_p)$ est isomorphe \`a un objet de la liste ${\bf D^*}$.
\end{itemize}
\end{thm}

\begin{rem}
{\em La condition $(3^{\circ})$, qui porte sur la repr\'esentation de Weil-Deligne associ\'ee \`a un objet $D$ de ${\bf MF}_{K/{\BBq}_p}(\varphi,N)$, se lit en prenant le polyn\^ome caract\'eristique de $\varphi$ sur le ${\BBq}_p$-espace vectoriel form\'e des \'el\'ements de $D$ qui sont fixes par l'action d'un rel\`evement du Frobenius absolu dans $G_{K/{\BBq}_p}$.}
\end{rem}

\begin{pf}
Les objets $D={\bf D}^*_{pst}(V_p)$ obtenus avec les conditions de (2) sont exactement ceux de la liste ${\bf D^*}$ : la repr\'esentation de Weil-Deligne associ\'ee \`a $D$ se lit sur le $(\varphi,N,G_{K/{\BBq}_p})$-module $D^{(0)}$ obtenu en oubliant la filtration et la m\^eme preuve que celle du th\'eor\`eme~\ref{lmodellipthm} montre que $D^{(0)}$ est l'un des $(\varphi,N,G_{K/{\BBq}_p})$-modules d\'eduits de la liste ${\bf D^*}$ ; puis la filtration sur $D_K$ est obtenue en \'ecrivant que $D$ est de type Hodge-Tate $(0,1)$ et faiblement admissible, voir~\ref{sec:classp}.

Pour les cas ${\bf D^*_m(e;b;\alpha)}$ le r\'esultat provient du fait que l'application de $p{\BBz}_p \backslash \{ 0 \}$ dans ${\BBq}_p$ qui \`a $q$ associe $\alpha(q)= - \log(u_q)/v_p(q)$ est surjective (o\`u l'on a \'ecrit $q=u_q p^{v_p(q)}$, cf. rmq.~\ref{multrmq} ; $\log(u_q)$ parcourt $p{\BBz}_p$ et $v_p(q)$ parcourt les entiers $\geq 1$), ainsi que de la description des twists quadratiques donn\'ee en~\ref{sec:twistsp}.

Pour les cas ${\bf D^*_c(e;a_p;\alpha)}$ et ${\bf D^*_{pc}(e;a_p;\epsilon;\alpha)}$ le r\'esultat provient du th\'eor\`eme~\ref{lmodellipthm} et du fait que pour toute courbe elliptique ordinaire sur ${\BBf}_p$ il existe un rel\`evement tel que $(*_{ord})$ est scind\'ee (ce qui \'equivaut \`a $\alpha =0$) ainsi qu'un rel\`evement tel que $(*_{ord})$ n'est pas scind\'ee (voir~\ref{sec:cpotord} et les exemples donn\'es en~\ref{sec:exordp}). 

Pour les cas ${\bf D^*_c(e;0)}$ le r\'esultat provient du th\'eor\`eme~\ref{lmodellipthm} (ici la filtration n'apporte pas de donn\'ee suppl\'ementaire).

Pour les cas ${\bf D^*_{pc}(e;0;\alpha)}$, si $(e,p)=(6,5)$ on se ram\`ene \`a la situation $(e,p)=(3,5)$ en tordant sur ${\BBq}_5(\pi_2)$ (voir rmq.~\ref{allsupsingrmq}(ii)). Le r\'esultat provient alors de l'\'etude faite en~\ref{sec:cpotsup} dont on reprend les notations. Soient $E_{\alpha, \epsilon}$ les courbes elliptiques sur ${\BBq}_p$ correspondant aux \'el\'ements de $O_{L_e}''$ avec $\alpha \in {\BBz}_p$ et $\epsilon \in \{ \pm 1 \}$. On a une application 
 $$\delta_e : O_{L_e}'' \rightarrow  {\BBp}^1({\BBq}_p)$$
qui \`a $E_{\alpha, \epsilon}$ associe l'unique $\delta_e(E_{\alpha, \epsilon}) \in {\BBp}^1({\BBq}_p)$ tel que ${\bf D}^*_{cris,K_e/{\BBq}_p}(V_p(E_{\alpha, \epsilon})) \simeq {\bf D^*_{pc}(e;0;}\delta_e(E_{\alpha, \epsilon}){\bf )}$ en tant que $(\varphi,G_{K_e/{\BBq}_p})$-modules filtr\'es. Un calcul montre alors que l'on a $\delta_e(E_{\alpha, \epsilon})= -p \alpha^{-\epsilon}$~; en particulier $\delta_e$ est surjective.
\end{pf}

\begin{rem}
{\em L'\'etude faite en~\ref{sec:cpotsup} montre que les fibres de $\delta_e$ sont finies : pour $\alpha \in {\BBp}^1({\BBq}_p)$ le cardinal de $\delta_e^{-1}(\alpha)$ est $4$ si $v_p(\alpha)=1$ et $2$ sinon.}
\end{rem}

\begin{rem}
{\em De m\^eme que pour les cas $\ell \neq p$, un ${\BBz}_p[G]$-module $T_p$ provient d'une courbe elliptique sur ${\BBq}_p$ si et seulement si le ${\BBq}_p[G]$-module ${\BBq}_p\!\otimes_{{\BBz}_p}\!T_p$ provient d'une courbe elliptique sur ${\BBq}_p$.}
\end{rem}

\begin{cor}
\label{nbclasspcor}
Le nombre de classes d'isomorphisme d'objets de ${\bf Rep}_{{\BBq}_p}(G)$ provenant d'une courbe elliptique sur ${\BBq}_p$ ayant potentiellement bonne r\'eduction est fini si et seulement si $p \equiv 1 \bmod 12$ ; il vaut alors $\,8[2\sqrt p\,] + 66 $.
\end{cor}

La premi\`ere assertion provient du fait que les courbes $E$ ayant potentiellement bonne r\'eduction avec $\mbox{dst}(E)\geq 3$ sont toutes potentiellement ordinaires si et seulement si $p \equiv 1 \bmod 12$ (et aussi bien s\^ur du th\'eor\`eme~\ref{pmodellipthm}). Il y a alors : $2\,\mbox{Card}({\cal N}_p^{\times})=4 [2\sqrt{p}\,]$ classes provenant de courbes elliptiques ayant bonne r\'eduction ordinaire sur ${\BBq}_p$ et autant provenant d'un twist quadratique ramifi\'e de telles ; 1 classe provenant de courbes ayant bonne r\'eduction supersinguli\`ere sur ${\BBq}_p$ et 1 provenant d'un twist quadratique ramifi\'e de telles ; $4\,\mbox{Card}({\cal N}_{p,3}^{\times})=24$ classes provenant de $E/{\BBq}_p$ potentiellement ordinaires avec $\mbox{dst}(E)=3$~; $4\,\mbox{Card}({\cal N}_{p,4}^{\times})=16$ classes provenant de telles courbes avec $\mbox{dst}(E)=4$ ; et $4\,\mbox{Card}({\cal N}_{p,6}^{\times})=24$ classes provenant de telles courbes avec $\mbox{dst}(E)=6$.

\medskip
Enfin, la proposition qui suit est essentiellement une cons\'equence de l'\'etude des courbes sur ${\BBq}_p$ ayant potentiellement bonne r\'eduction faite en~\ref{sec:construct}.

\begin{prop}
\label{isogprop}
Soit $E/{\BBq}_p$ une courbe elliptique. L'assertion :
$$(*)\makebox[1.5cm]{}  V_p(E) \mbox{ et } V_p(E') \mbox{ sont } {\BBq}_p[G]\mbox{-isomorphes} \ \ \Leftrightarrow \ \ E \mbox{ et } E' \mbox{ sont }  {\BBq}_p\mbox{-isog\`enes} $$
est vraie si et seulement si on est dans l'un des trois cas suivants : \\
(1) $E$ a potentiellement mauvaise r\'eduction multiplicative  \\
(2) $E$ a potentiellement bonne r\'eduction supersinguli\`ere et $\mbox{\em dst}(E)\geq 3$  \\
(3) $E$ est le rel\`evement canonique sur une extension finie de ${\BBq}_p$ d'une courbe ordinaire.
\end{prop}

\begin{pf}
L'implication $(1) \Rightarrow (*)$ provient pour des courbes de Tate de la remarque~\ref{multisogrmq} et le cas g\'en\'eral s'en d\'eduit par torsion quadratique. L'implication $(2) \Rightarrow (*)$ provient des remarques~\ref{supisogrmq},~\ref{allsupsingrmq} et~\ref{twistsupisogrmq}. L'implication $(3) \Rightarrow (*)$ provient de la remarque~\ref{canordisogrmq} pour $\mbox{dst}(E)=1$ ou $2$ et des remarques~\ref{potcanordisogrmq} et~\ref{allordrmq} pour $\mbox{dst}(E)\geq 3$. Enfin, les remarques~\ref{supnonisogrmq} ainsi que~\ref{ordnonisogrmq} et~\ref{allordrmq} montrent que dans tous les autres cas l'assertion $(*)$ est fausse. 
\end{pf}

\subsection{Exemples}
\label{sec:exp}

\subsubsection{Courbes elliptiques potentiellement ordinaires}
\label{sec:exordp}

\underline{Si $4 \mid p-1$} : on reprend les courbes $\widetilde{E}_j/{\BBf}_p$ et $E_{i,j}/{\BBq}_p$ avec $1 \leq j \leq 4$ et $0\leq i \leq 3$ consid\'er\'ees en~\ref{sec:exordl}. On a alors pour chaque $j$ :
$$\begin{array}{rclcrcl}
{\bf D}^*_{pcris}(V_p(E_{0,j})) & \simeq & {\bf D^*_c(1;a_{p,\,j};0)} & \hspace{0.5cm} & {\bf D}^*_{pcris}(V_p(E_{1,j})) & \simeq & {\bf D^*_{pc}(4;a_{p,\,j};1;0)}       \\
{\bf D}^*_{pcris}(V_p(E_{2,j})) & \simeq  & {\bf D^*_c(2;a_{p,\,j};0)} & \hspace{0.5cm} & {\bf D}^*_{pcris}(V_p(E_{3,j})) & \simeq & {\bf D^*_{pc}(4;a_{p,\,j};-1;0)}
\end{array}$$

On pose $E_{i,j}' : y^2 = x^3 + [u_j](-p)^i x + (-p)^{n(i)}$ pour $1 \leq j \leq 4$, $0\leq i \leq 3$ et $n(i)= 1,2,4,5$ si $i=0,1,2,3$ respectivement. Ce sont des courbes sur ${\BBq}_p$ d'invariant modulaire entier congru \`a $1728$ modulo $p{\BBz}_p$ mais diff\'erent de $1728$ ; elles acqui\`erent bonne r\'eduction sur ${\BBq}_p(\pi_4)$ et la courbe r\'eduite de $E_{i,j}'$ est $\widetilde{E}_j$. On a alors pour chaque $j$~:
 $$\begin{array}{rclcrcl}
{\bf D}^*_{pcris}(V_p(E_{0,j}')) & \simeq & {\bf D^*_c(1;a_{p,\,j};1)} & \hspace{0.5cm} & {\bf D}^*_{pcris}(V_p(E_{1,j}')) & \simeq & {\bf D^*_{pc}(4;a_{p,\,j};1;1)}    \\
{\bf D}^*_{pcris}(V_p(E_{2,j}')) & \simeq & {\bf D^*_c(2;a_{p,\,j};1)} & \hspace{0.5cm} & {\bf D}^*_{pcris}(V_p(E_{3,j}')) & \simeq & {\bf D^*_{pc}(4;a_{p,\,j};-1;1)}
\end{array}$$

\medskip
\noindent \underline{Si $3 \mid p-1$} : on reprend les courbes $\widetilde{\cal E}_j/{\BBf}_p$ et ${\cal E}_{i,j}/{\BBq}_p$ avec $1 \leq j \leq 6$ et $0\leq i \leq 5$ consid\'er\'ees en~\ref{sec:exordl}. On a alors pour chaque $j$ :
 $$\begin{array}{rclcrcl}
{\bf D}^*_{pcris}(V_p({\cal E}_{0,j})) & \simeq & {\bf D^*_c(1;a_{p,\,j};0)} & \hspace{0.5cm} & {\bf D}^*_{pcris}(V_p({\cal E}_{1,j})) & \simeq &  {\bf D^*_{pc}(6;a_{p,\,j};1;0)}   \\
{\bf D}^*_{pcris}(V_p({\cal E}_{2,j})) & \simeq & {\bf D^*_{pc}(3;a_{p,\,j};1;0)} & \hspace{0.5cm} & {\bf D}^*_{pcris}(V_p({\cal E}_{3,j})) & \simeq & {\bf D^*_c(2;a_{p,\,j};0)} \\
{\bf D}^*_{pcris}(V_p({\cal E}_{4,j})) & \simeq & {\bf D^*_{pc}(3;a_{p,\,j};-1;0)} & \hspace{0.5cm} & {\bf D}^*_{pcris}(V_p({\cal E}_{5,j})) & \simeq & {\bf D^*_{pc}(6;a_{p,\,j};-1;0)}  
\end{array}$$

On pose ${\cal E}_{i,j}' : y^2 = x^3 +(-p)^{m(i)}x+ [v_j](-p)^i$ pour $1 \leq j \leq 6$, $0\leq i \leq 5$ et $m(i)= 1,1,2,3,3,4$ si $i=0,1,2,3,4,5$ respectivement. Ce sont des courbes sur ${\BBq}_p$ d'invariant modulaire entier congru \`a $0$ modulo $p{\BBz}_p$ mais non nul~; elles acqui\`erent bonne r\'eduction sur ${\BBq}_p(\pi_6)$ et la courbe r\'eduite de ${\cal E}_{i,j}'$ est $\widetilde{\cal E}_j$. On a alors pour chaque $j$ :
$$\begin{array}{rclcrcl}
{\bf D}^*_{pcris}(V_p({\cal E}_{0,j}')) & \simeq & {\bf D^*_c(1;a_{p,\,j};1)} & \hspace{0.5cm} & {\bf D}^*_{pcris}(V_p({\cal E}_{1,j}')) & \simeq & {\bf D^*_{pc}(6;a_{p,\,j};1;1)}  \\
{\bf D}^*_{pcris}(V_p({\cal E}_{2,j}')) & \simeq & {\bf D^*_{pc}(3;a_{p,\,j};1;1)} & \hspace{0.5cm} & {\bf D}^*_{pcris}(V_p({\cal E}_{3,j}')) & \simeq & {\bf D^*_c(2;a_{p,\,j};1)}      \\
{\bf D}^*_{pcris}(V_p({\cal E}_{4,j}')) & \simeq & {\bf D^*_{pc}(3;a_{p,\,j};-1;1)} & \hspace{0.5cm} & {\bf D}^*_{pcris}(V_p({\cal E}_{5,j}')) & \simeq & {\bf D^*_{pc}(6;a_{p,\,j};-1;1)}    
\end{array}$$

\subsubsection{Courbes elliptiques potentiellement supersinguli\`eres}
\label{sec:exsupp}

\underline{Si $4 \mid p+1$} : on reprend les courbes $E_i/{\BBq}_p$ avec $0\leq i\leq 3$ consid\'er\'ees en~\ref{sec:exsupl}. On a alors :
 $$\begin{array}{rclcrcl}
{\bf D}^*_{pcris}(V_p(E_0)) & \simeq & {\bf D^*_c(1;0)} & \hspace{0.5cm} & {\bf D}^*_{pcris}(V_p(E_1)) & \simeq & {\bf D^*_{pc}(4;0;\infty)}  \\
{\bf D}^*_{pcris}(V_p(E_2)) & \simeq & {\bf D^*_c(2;0)} & \hspace{0.5cm} & {\bf D}^*_{pcris}(V_p(E_3)) & \simeq & {\bf D^*_{pc}(4;0;0)}
\end{array}$$

\medskip
\noindent \underline{Si $3 \mid p+1$} : on reprend les courbes ${\cal E}_i/{\BBq}_p$ avec $0\leq i\leq 5$ consid\'er\'ees en~\ref{sec:exsupl}. On a alors :
$$\begin{array}{rclcrcl}
{\bf D}^*_{pcris}(V_p({\cal E}_0)) & \simeq & {\bf D^*_c(1;0)} & \hspace{0.5cm} & {\bf D}^*_{pcris}(V_p({\cal E}_1)) & \simeq & {\bf D^*_{pc}(6;0;\infty)}        \\
{\bf D}^*_{pcris}(V_p({\cal E}_2)) & \simeq & {\bf D^*_{pc}(3;0;\infty)} & \hspace{0.5cm} & {\bf D}^*_{pcris}(V_p({\cal E}_3)) & \simeq  & {\bf D^*_c(2;0)}             \\
{\bf D}^*_{pcris}(V_p({\cal E}_4)) & \simeq & {\bf D^*_{pc}(3;0;0)} & \hspace{0.5cm} & {\bf D}^*_{pcris}(V_p({\cal E}_5)) & \simeq & {\bf D^*_{pc}(6;0;0)} 
\end{array}$$

\bigskip

Department of Mathematical Sciences, University of Durham, Durham DH1 3LE, England

E-mail : maja.volkov@durham.ac.uk


\begin{thebibliography}{Wa-Mi}



\bibitem[Co-Fo]{Co-Fo} {\em P. Colmez et J.-M. Fontaine}, Construction des repr\'esentations $p$-adiques semi-stables, Invent. math. {\bf 140}, 1 (2000), 1-43.

\bibitem[C-D-T]{C-D-T} {\em B. Conrad, F. Diamond, and R. Taylor}, Modularity of certain potentially Barsotti-Tate Galois representations, J. of the Am. Math. Soc. {\bf 12}, 2 (1999), 521-567.

\bibitem[De 1]{De 1} {\em P. Deligne}, Les constantes des \'equations fonctionnelles des fonctions $L$, {\em in} Modular Functions of One Variable II, LNM {\bf 349}, Springer-Verlag (1973), 501-595.

\bibitem[De 2]{De 2} {\em P. Deligne}, La conjecture de Weil II, Publ. Math. IHES {\bf 52} (1980), 137-252.

\bibitem[Fo 1]{Fo 1} {\em J.-M. Fontaine}, Le corps des p\'eriodes $p$-adiques, expos\'e II, {\em in} P\'eriodes $p$-adiques, Ast\'erisque {\bf 223}, Soc. Math. de France (1994).

\bibitem[Fo 2]{Fo 2} {\em J.-M. Fontaine}, Repr\'esentations $p$-adiques semi-stables, expos\'e III, {\em in} P\'eriodes $p$-adiques, Ast\'erisque {\bf 223}, Soc. Math. de France (1994). 

\bibitem[Fo 3]{Fo 3} {\em J.-M. Fontaine}, Repr\'esentations $\ell$-adiques potentiellement semi-stables, expos\'e VIII, {\em in} P\'eriodes $p$-adiques, Ast\'erisque {\bf 223}, Soc. Math. de France (1994).

\bibitem[Fo 4]{Fo 4} {\em J.-M. Fontaine}, Groupes $p$-divisibles sur les corps locaux, Ast\'erisque {\bf 47-48}, Soc. Math. de France (1977).

\bibitem[Fo 5]{Fo 5} {\em J.-M. Fontaine}, Sur certains types de repr\'esentations $p$-adiques du groupe de Galois d'un corps local ; construction d'un anneau de Barsotti-Tate, Annals of Math. {\bf 115} (1982), 529-577.

\bibitem[Fo-Ma]{Fo-Ma} {\em J.-M. Fontaine and B. Mazur}, Geometric Galois Representations, {\em in} Conference on Elliptic Curves and Modular Forms, Hong Kong, December 18-21 (1995), 41-77.

\bibitem[Ho-Ta]{Ho-Ta} {\em J. Tate}, Classes d'isog\'enie des vari\'et\'es ab\'eliennes sur un corps fini (d'apr\`es {\em T. Honda}), S\'eminaire Bourbaki {\bf 352} (1968), 15p.

\bibitem[Ka]{Ka} {\em N. Katz}, Serre-Tate local moduli, {\em in} Surfaces alg\'ebriques, LNM {\bf 868}, Springer-Verlag (1981), 138-202.

\bibitem[Kr]{Kr} {\em A. Kraus}, D\'etermination du poids et du conducteur associ\'es aux repr\'esentations des points de $p$-torsion d'une courbe elliptique, Diss. Math. {\bf 364} (1997), 39p.

\bibitem[La]{La} {\em S. Lang}, Abelian Varieties, Interscience Publishers (1959).

\bibitem[LS]{LS} {\em B. Le Stum}, La structure de Hyodo-Kato pour les courbes, Rend. Sem. Mat. Univ. Padova {\bf 94} (1995), 279-301.

\bibitem[Ma]{Ma} {\em B. Mazur}, On monodromy invariants occuring in global arithmetic, and Fontaine's theory, {\em in} $p$-adic Monodromy and the Birch and Swinnerton-Dyer Conjecture (Boston, MA, 1991), Contemp. Math. {\bf 165} (1994), 1-20.

\bibitem[M-T-T]{M-T-T} {\em B. Mazur, J. Tate, and J. Teitelbaum}, On $p$-adic analogues of the conjectures of Birch and Swinnerton-Dyer, Invent. math. {\bf 84} (1986), 1-48.

\bibitem[Me]{Me} {\em W. Messing}, The Crystals associated to Barsotti-Tate Groups : with Applications to Abelian Schemes, LNM {\bf 264}, Springer-Verlag (1972). 

\bibitem[Ra]{Ra} {\em M. Raynaud}, 1-Motifs et monodromie g\'eom\'etrique, expos\'e VII, {\em in} P\'eriodes $p$-adiques, Ast\'erisque {\bf 223}, Soc. Math. de France (1994).

\bibitem[Ro]{Ro} {\em D.E. Rohrlich}, Elliptic Curves and the Weil-Deligne Group, CRM Proceedings and Lecture Notes {\bf 4} (1994). 

\bibitem[Se 1]{Se 1} {\em J.-P. Serre}, Abelian $\ell$-adic representations and elliptic curves (2nd ed.), Addison-Wesley (1989). 

\bibitem[Se 2]{Se 2} {\em J.-P. Serre}, Propri\'et\'es galoisiennes des points d'ordre fini des courbes elliptiques, Invent. math. {\bf 15} (1972), 259-331.

\bibitem[Se-Ta]{Se-Ta} {\em J.-P. Serre and J. Tate}, Good reduction of abelian varieties, Annals of Math. {\bf 88} (1968), 492-517. 

\bibitem[Si 1]{Si 1} {\em J.H. Silverman}, The Arithmetic of Elliptic Curves, GTM {\bf 106}, Springer-Verlag (1986).

\bibitem[Si 2]{Si 2} {\em J.H. Silverman}, Advanced Topics in the Arithmetic of Elliptic Curves, GTM {\bf 151}, Springer-Verlag (1994). 

\bibitem[Ta]{Ta} {\em J. Tate}, Endomorphisms of Abelian Varieties over Finite Fields, Invent. math. {\bf 2} (1966), 134-144. 

\bibitem[Wa-Mi]{Wa-Mi} {\em W.C. Waterhouse and J.S. Milne}, Abelian Varieties over Finite Fields, {\em in} AMS Proceedings of Symposia in Pure Mathematics {\bf XX} (1971), 53-64.

\bibitem[We]{We} {\em A. Weil}, The Field of Definition of a Variety, Am. J. of Math. {\bf 78} (1956), 509-524.


\addcontentsline{toc}{section}{R\'ef\'erences}


\end{thebibliography}
\end{document}